\pdfoutput=1
\PassOptionsToPackage{normalem,normalbf}{ulem}
\documentclass[VANCOUVER,STIX1COL]{WileyNJD-v2}

\usepackage{ifthen}

\usepackage{amssymb}

\usepackage[final]{changes}

\definechangesauthor[name=R1, color=blue]{R1}
\definechangesauthor[name=R2, color=red]{R2}
\definechangesauthor[name=R3, color=green]{R3}

\newboolean{arxiv}
\newboolean{sisc}
\newboolean{linearalgebra}
\newboolean{appendix}

\setboolean{arxiv}{false}
\setboolean{sisc}{false}
\setboolean{linearalgebra}{true}
\setboolean{appendix}{false}

\usepackage{graphicx}

\ifthenelse{\boolean{arxiv}}{
}{} 
\ifthenelse{\boolean{sisc}}{
	\usepackage[numbers,sort&compress]{natbib}
}{} 
\ifthenelse{\boolean{linearalgebra}}{
}{} 

\newcommand{\TheTitle}{\replaced[id=R3]{Stabilised Asynchronous Fast Adaptive
Composite Multigrid using Additive Damping}{Dynamically Adaptive FAS for an Additively Damped AFAC Variant}}

\newcommand{\TheAuthors}{C.~D.~Murray and T.~Weinzierl}

\usepackage{algorithm}
\usepackage{algpseudocode}

\usepackage{enumitem}

\newtheorem{idea}{Idea}

\ifthenelse{\boolean{arxiv}}{
}{} 
\ifthenelse{\boolean{sisc}}{
	\newtheorem{observe}{Observation}
}{} 
\ifthenelse{\boolean{linearalgebra}}{

}{} 

\ifthenelse{\boolean{arxiv}}{
  \usepackage{geometry}
\geometry{left=4.1cm,textwidth=14.2cm,top=2cm,textheight=24cm}

\title{\TheTitle}

}{} 
\ifthenelse{\boolean{sisc}}{
\headers{\TheTitle}{\TheAuthors}

\title{
 \TheTitle
 \thanks{
  This is an extended and revised version of a paper
  submitted to the Copper Mountain 2018 Student Paper competition.
  \funding{
    The work was funded by a Durham University/EPSRC DTA PhD scholarship.
    It made use of the facilities of the Hamilton HPC Service of Durham
    University.
  }
 }
}

\slugger{sisc}{xxxx}{xx}{x}{x--x}

\usepackage{amsopn}
\DeclareMathOperator{\diag}{diag}

\author{
  Charles D. Murray
  \thanks{
    Department of Computer Science, Durham University,
    Great Britain 
    (\email{c.d.murray@durham.ac.uk}).
  }
  \and
  Tobias Weinzierl
  \thanks{
    Department of Computer Science, Durham University,
    Great Britain 
    \url{http://www.dur.ac.uk/tobias.weinzierl}).
  }
}

}{} 
\ifthenelse{\boolean{linearalgebra}}{

\title{
	\TheTitle 
	\protect\thanks{
 		The work was funded by a Durham University/EPSRC DTA PhD scholarship.
 		Award reference 1764342.
  		It made use of the facilities of the Hamilton HPC Service of Durham
  		University.  		
  	}
  }
\articletype{Article Type}%

\received{26 April 2016}
\revised{6 June 2016}
\accepted{6 June 2016}

\usepackage{amsopn}
\DeclareMathOperator{\diag}{diag}

\author[1]{Charles D. Murray*}

\author[1]{Tobias Weinzierl}

\authormark{Murray and Weinzierl}

\address[1]{\orgdiv{Department of Computer Science}, \orgname{Durham University}, 
\orgaddress{\state{Durham}, \country{United Kingdom}}}

\corres{*Charles Murray, Department of Computer Science, Durham University, Durham, United Kingdom.
 \email{c.d.murray@durham.ac.uk}}


\jnlcitation{\cname{%
\author{Charles D Murray},  and 
\author{Tobias Weinzierl}} (\cyear{<year>}), 
\ctitle{<journal title>}, \cjournal{<journal name>} <year> <vol> Page <xxx>-<xxx>}

}{}

\begin{document}

\ifthenelse{\boolean{linearalgebra}}{
	\abstract[Abstract]{
Multigrid solvers face multiple challenges on parallel computers.
Two fundamental ones \replaced[id=R3]{read as
follows:}{\replaced[id=R2]{being}{are that}} 
\replaced[id=R3]{Multiplicative}{multiplicative} solvers issue coarse
grid solves which exhibit low concurrency and \deleted[id=R3]{that}many
multigrid implementations suffer from an expensive coarse grid identification
phase \replaced[id=R3]{plus}{\replaced[id=R2]{and}{as well as dynamic}} adaptive
mesh refinement \deleted[id=R3]{(AMR)} overhead.
We \deleted[id=R2]{therefore }propose a new additive
\added[id=R3]{multigrid} variant\replaced[id=R3]{ for spacetrees, i.e.~meshes
as they are constructed from octrees and quadtrees:}{ of the fast adaptive
composite (FAC) method which can be combined with Full Approximation Storage (FAS) plus BoxMG inter-grid transfer operators on spacetrees\replaced[id=R2]{:
 Additivity implies a high concurrency level,
 FAS implies that the major solver ingredients such as smoothers
 can work on logically regular grids,
 and BoxMG implies that the solver benefits from algebraic multigrid's accuracy
 even though it relies on a geometric grid hierarchy.}{.
This allows for a
straightforward realisation of arbitrary dynamic AMR on geometric
multiscale grids with algebraic operators.}}
\added[id=R3]{
 It is an additive scheme, i.e.~all multigrid resolution levels are updated
 concurrently.
 This ensures a high concurrency level, while the transfer operators between the
 mesh levels can still be constructed algebraically.
}
The novel flavour of the additive scheme is 
an augmentation of the solver with an additive, auxiliary damping \added[id=R2]{parameter }per grid
level \added[id=R3]{per vertex} that is in turn constructed through the next
coarser level---an idea which utilises smoothed aggregation principles or the motivation
behind AFACx:
\added[id=R3]{
 Per level, we solve an additional equation whose purpose is to damp too
 aggressive solution updates per vertex which would otherwise, in combination
 with all the other levels, yield an overcorrection and, eventually, oscillations.
 This additional equation is constructed additively as well, i.e.~is once more
 solved concurrently to all other equations.
} 
This yields improved stability\replaced[id=R3]{, closer to what is seen with}{ as we experience it}
with multiplicative schemes, while pipelining techniques help us to write down the additive solver with single-touch
semantics
\added[id=R3]{ for dynamically adaptive meshes}.
	}
}{}

\maketitle

\ifthenelse{\boolean{sisc}}{
	\begin{abstract}
  Multigrid solvers face multiple challenges on parallel computers.
  Two fundamental ones are that multiplicative solvers issue coarse grid solves
  in regular intervals which exhibit low concurrency by construction,
  while additive solvers tend to be less stable and efficient,    
  and that many multigrid implementations suffer from an expensive mesh setup/coarse grid
  identification phase as well as dynamic adaptive mesh refinement
  (AMR) overhead.
  We therefore propose to combine the additive variant of the fast adaptive
  composite (FAC) method with Full Approximation Storage (FAS) plus
  BoxMG inter-grid transfer operators on spacetrees.
  This allows for a
  straightforward realisation of arbitrary dynamic AMR on geometric
  multiscale grids with algebraic operators.
  Additive multigrid variants avoid isolated coarse grid solves.
  An augmentation of the solver with an additive damping per grid level that is
  in turn derived from the next coarser level---an idea
  which utilises smoothed aggregation principles---yields our additively damped
  asynchronous fast adaptive composite extended (adaFACx) scheme and exhibits
  improved stability.
  Once we throttle the propagation of the algebraic operators to one resolution
  level per iterate, once we tailor auxiliary stencils, and
  once we exploit pipelining techniques, the additive solver can be written with single-touch semantics.    
	\end{abstract}
}{}

\ifthenelse{\boolean{sisc}}{

 
}
{}

\section{Introduction}
\label{section:introduction}

%
%
The elliptic partial differential equation (PDE)
\begin{equation}
 - \nabla\added[id=R2]{\cdot } \left( \epsilon \nabla \right) u = f, \qquad
 \epsilon: \Omega \subset \mathbb{R}^d \mapsto \mathbb{R}^+ \mbox{\added[id=R1]{positive, bounded
 away from zero, and} either constant or varying}
 \label{equation:PDE}
\end{equation}

\noindent
serves as a building block in many applications.
Examples are chemical dispersion in sub-surface reservoirs, the
heat distribution in buildings, or the diffusion of oxygen in tissue.
It is also the starting point to construct 
more complex differential operators.
Solving this PDE quickly is important yet not trivial.
One reason is buried within the operator: 
any local modification\added[id=R1]{s} of
the solution propagate\deleted[id=R1]{s} \replaced[id=R1]{through}{among} the whole
computational domain, \replaced[id=R1]{though this effect can be}{unless} damped
\replaced[id=R1]{by}{out, i.e.~effectively stopped, due to} large
$\epsilon $ variations.
The operator exhibits multiscale behaviour.
A successful family of iterative techniques to solve (\ref{equation:PDE})
hence is multigrid.
It relies on representations of the
operator's behaviour on multiple scales.
It builds the operator's multiscale behaviour into the algorithm.

%
%
There are \replaced[id=R1]{numerical}{conceptional}, algorithmic and
implementational hurdles that must be tackled when we write multigrid codes.
In this paper, we focus on three
\replaced[id=R1]{algorithmic/implementation}{conceptual} challenges which should
be addressed before we scale up multigrid.
(i) State-of-the-art multigrid codes have to support dynamically adaptive mesh
refinement (AMR) without significant overhead.
While constructing coarser (geometric) representations from regular grids is
straightforward, it is non-trivial for adaptive meshes.
\added[id=R1]{Furthermore, setup phases become costly if the $\epsilon $
distribution is complex \cite{Lin:18:Performance}, or if the mesh changes
frequently \cite{Weinzierl:18:BoxMG}.
}
(ii) If an algorithm
solves problems on cascades of coarser and coarser, i.e.~smaller and smaller, problems, the smallest problems eventually do not exhibit
enough computational work to scale among larger core counts.
\added[id=R1]{For adaptive meshes, such low concurrency phases can---depending
on the implementation---arise for fine grid solves, too, if we run through the
multigrid levels resolution by resolution.}
(iii) If an algorithm projects a problem to multiple resolutions and then
constructs a solution from these resolutions, its implementation tends to read
and write data multiple times \added[id=R1]{using indirect or scattered data
accesses}.
Repeated data access however is poisonous on today's hardware which suffers
from a widening gap between what cores could compute and what throughput
\added[id=R1]{the} memory can provide \cite{Dongarra:14:ApplMathExascaleComputing}.

%
%
\added[id=R1]{
 Numerous concepts have been proposed to tackle this triad of
 challenges.
 We refrain from a comprehensive overview but sketch some particular popular 
 code design decisions:
 Many successful multigrid codes use a cascade
 of geometric grids which are embedded into each other, and \replaced[id=R2]{then}{they}
 run through the resolutions level by level, i.e.~embedding by embedding.
 This simplifies the coarse grid identification, and, given a sufficiently
 homogeneous refinement pattern, implies that a fine grid decomposition induces
 a fair partitioning of coarser resolutions. 
 Many codes actually start from a coarse resolution and make the refinement
 yield the finer mesh levels \cite{Gmeiner:14:Parallel}:
 Combining this with algebraic coarse grid identification once the
 grid is reasonably coarse adds additional flexibility to a geometric scheme.
 It allows codes to treat massive $\epsilon $ variations (or complex geometries)
 through a coarse mesh
 \cite{Lu:14:HybridMG,May:15:Scalable,Sundar:12:ParallelMultigrid}, whereas the geometric multigrid component handles the bulk of the compute work and exploits structured grids.
 Many successful codes furthermore use classic rediscretisation on finer grids
 and employ expensive algebraic operator computation only on coarser
 meshes.
 They avoid the algebraic overhead to assemble accurate fine grid operators
 \cite{Gholami:16:SolverComparison}.
 Many large-scale codes do not use the
 recursive multigrid paradigm across all possible mesh resolutions, but switch to alternative solvers such as Krylov schemes for systems that are still
 relatively big \cite{May:15:Scalable,Rudi:15:Extreme}.
 They avoid low concurrency phases arising from very small multigrid
 subproblem solves or exact inversions of very small systems.
 Finally, our own work \cite{Reps:17:Helmholtz,Weinzierl:18:BoxMG} has studied
 rewrites of multigrid with single-touch semantics\replaced[id=R2]{,
 i.e.~each unknown is, amortised, fetched from the main memory into the
 caches once per grid sweep or additive cycle, respectively.
 Similar to other}{to reduce disadvantageous memory access behaviour. Other}
 strategies such as pipelining or smoother and
 stencil optimisation 
 \cite{Ghysels:12:PolynomialSmoothers,Ghysels:13:ModelMG,Gmeiner:15:HHG}
 \added[id=R2]{our single-touch rewrites} reduce the data movement of solvers as
 well as indirect and scattered memory accesses.
 If meshes exhibit steep adaptivity, i.e.~refine
 specific subregions to be particularly detailed, or if problems have very strongly
 varying coefficients, all approaches will run into issues.
 These take the form of 
 scalability challenges (from low concurrency phases on the coarse mesh or
 situations where only small parts of the mesh are updated as we lack a global,
 uniform coarsening scheme), load balancing challenges (a geometric fine grid
 splitting does not map to coarser resolutions anymore and there's no inclusion property of multiscale domain partitions), 
 or materialise in the memory overhead and performance penalty of algebraic
 multigrid \cite{Gholami:16:SolverComparison}.
}

\replaced[id=R1]{
 The majority of multigrid papers focus on its multiplicative form,
 as it exhibits superior convergence rates compared to additive alternatives;
 consequently such formulations are more often used as preconditioner.
 Additive multigrid however remains an interesting solver alternative to
 multiplicative multigrid in}
{Due to its reduced synchronisation and the absence of
isolated small-scale solves, the additive mindset is an interesting alternative
for} the era of massive concurrency growth
\cite{Dongarra:14:ApplMathExascaleComputing}.
\added[id=R1]{In the present paper, we}
propose \replaced[id=R1]{an additive}{a} solver-implementation combination which
tackles \added[id=R1]{introductory} challenges \replaced[id=R1]{by combining
three key concepts}{in one rush}.
\deleted[id=R1]{It combines several state-of-the-art building blocks.}
\replaced[id=R1]{The first concept is a pure geometric construction through}{Our
first algorithmic block is } the spacetree paradigm, a generalisation of the classic octree/quadtree idea
\cite{Weinzierl:11:Peano,Weinzierl:19:Peano}.
Spacetrees yield adaptive Cartesian grids which are nested within
each other
\cite{Weinzierl:11:Peano,Reps:17:Helmholtz,Weinzierl:18:BoxMG,Weinzierl:19:Peano}.
Adaptivity decreases the cost of a solve by
reducing the degrees of freedom without adversely affecting the
accuracy.
\replaced[id=R1]{Additional computational effort is invested
where it improves the solution significantly.
}{It invests work where it pays off.} 
With complex boundary conditions or
non-trivial $\epsilon$---or even $\epsilon (u)$ which renders (\ref{equation:PDE}) nonlinear---the regions where to refine
\replaced[id=R1]{may not be}{are not} known a priori, \added[id=R1]{might be
expensive to determine, or depend on the right-hand side (if the PDE is employed within  a
time-stepping scheme, for example)}.
Schemes that allow for dynamic mesh refinement are therefore key for many
applications. 
\added[id=R1]{On top of this, the spacetree idea yields a natural
multiresolution cascade well-suited for multigrid (cmp.~(i) above).} 
\replaced[id=R1]{The second concept is the increase of asynchronicity and
concurrency through additive multigrid (cmp.~(ii) above).}
{Our second algorithmic building
block is additive multigrid.
Additive multigrid exhibits greater parallelism than
conventional multiplicative multigrid, as different levels are computed independently of each other.
There is no close-to-serial coarse grid solve.
} 
\replaced[id=R1]{Our final concept is the application of mature implementation
patterns to our algorithms such that we obtain a single-touch multigrid solver
with low memory footprint (cmp.~(iii) above):
We rely on}{Our third algorithmic building block is} the triad of fast adaptive composite
(FAC), hierarchical transformation multigrid (HTMG) \cite{Griebel:90:HTMG}
and full approximation storage (FAS) \cite{Trottenberg:01:Multigrid}.
These three techniques allow us to elegantly realise a multigrid scheme which
straightforwardly works for dynamically adaptive meshes.
\replaced[id=R1]{We}{Fourth, it unfolds its full potential once we} merge it
with quasi matrix-free multigrid relying on algebraic BoxMG inter-grid transfer operators
\cite{Dendy:82:BlackboxMG,Dendy:10:BoxMgBy3,Weinzierl:18:BoxMG}.
\replaced[id=R1]{We also utilise}{Our last algorithmic building block is} 
pipelining combined with
recursive element-wise grid traversals. 
We run through the spacetree depth-first which yields excellent cache hit rates
\cite{Weinzierl:11:Peano} and simple recursive implementations.
The approach equals a multiscale element-wise grid traversal.
Going from coarse levels to fine levels and backtracking however does not fit
straightforwardly to FAS with additive multigrid, where we restrict the residual
from fine to coarse, prolong the correction from coarse to fine and inject the
solution from fine to coarse again.
Yet, we know that some additional auxiliary variables allow us to write
additive solvers \added[id=R2]{as }single-touch \cite{Reps:17:Helmholtz}.
Each unknown is read into the chip's caches
once per cycle.

\replaced[id=R1]{None of the enlisted ingredients or their implementation
flavour as enlisted are new.
Our novel contribution is a modification of the additive formulation plus the
demonstration that this modification still fits with the other presented ideas:} 
{However,}
Plain additive approaches face a severe problem.
They are less robust \added[id=R1]{than their multiplicative counterparts}.
Na\"ively restricting residuals to
multiple levels and eliminating errors concurrently tends to make the iterative
scheme overshoot \cite{Bastian:98:Additive,Wolfson:19:AsynchronousMultigrid}. 
Multiple strategies exist to improve the stability without compromising on the
additivity.
In the simplest case, \deleted[id=R1]{we merely employ }additive multigrid
\added[id=R1]{is employed} as a preconditioner and \deleted[id=R1]{use }a more
robust solver \added[id=R1]{is used} thereon.
Our work goes down the ``multigrid as a solver'' route:
A well-known approach to mitigate overshooting in the solver is to more aggressively 
damp levels the coarser they are.
This reduces their impact and improves stability but
decreases the rate of convergence \cite{Reps:17:Helmholtz}.
We refrain from such resolution-parameterised damping and 
follow up on the idea behind AFACx
\cite{McCormick:86:FAC,McCormick:89:AFAC,Lee:04:AFAC,Phillip:00:Elliptic}:
By introducing an additional correction component per level, our approach
predicts additive overshooting from coarser levels.
Different to AFACx, we however do not make the additional auxiliary solves
preprocessing steps. 
We phrase them\added[id=R1]{ in a} completely parallel (additive)\added[id=R1]{
way} to the actual correction's solve.
To make the auxiliary contributions meaningful nevertheless, we tweak them
through ideas resembling smoothed aggregation
\cite{Tuminaro:00:Parallel,Vanvek:96:Algebraic,Vanvek:95:Fast}
which approximate the smoothing steps of multiplicative multigrid
\cite{Yang:14:Reducing}.
We end up with an \emph{additively damped Asynchronous FAC} (adAFAC).

Our adAFAC implementation merges the levels of the multigrid scheme plus their
traversal into each other, and thus provides a single-touch implementation.
Through this, we eliminate synchronisation between the solves on
different resolution levels and anticipate that FAC yields multigrid grid
sequences where work non-monotonously grows and shrinks upon each resolution
transition.
\replaced[id=R1]{
 Additive literature usually emphasises the advantage of ``additivity'' in that
 the individual levels can be processed independently.
 Recent}{Notably recent}
work on a further decoupling of both the individual levels'
solves as well as the solves within a level \cite{Wolfson:19:AsynchronousMultigrid} shows
great upscaling potential.
\added[id=R1]{
 Our strategy allows us to head in the other direction: 
}
We vertically integrate solves
\cite{Adams:16:SegmentalRefinement}\added[id=R1]{, i.e.~we partition the finest mesh where the residual is computed, apply this decomposition vertically---to all mesh 
resolutions---then merge the traversal of multiple meshes per
subpartition.
The multigrid challenge to balance not one mesh but a cascade of meshes becomes
a challenge of balancing a single set of jobs again.
}

%
%
We reiterate which algorithmic ingredients we use in Section
\ref{section:ingredients} before we introduce our new additive solver
adAFAC.
\added[id=R1]{
 This Section \ref{section:integration} is the main new contribution of the
 present text.
}
Section \ref{section:single-touch} then translates adAFAC into a 
single-touch algorithm blueprint.
Some numerical results \replaced[id=R1]{outline}{uncover} the solver's potential
(Section \ref{section:results}):
\added[id=R1]{
 An extensive comparison of solver variants, upscaling studies or the
 application to real-world problems are out of scope. 
 Yet, adAFAC adds an interesting novel solver variant to the known suite of
 multigrid techniques available to scientists and engineers.
}
We close the
discussion with a brief summary and sketch
future work.

\section{Related work and methodological ingredients}
\label{section:ingredients}


\subsection{Spacetrees}

\begin{figure}[htb]
 \begin{center}
  \includegraphics[width=0.23\textwidth]{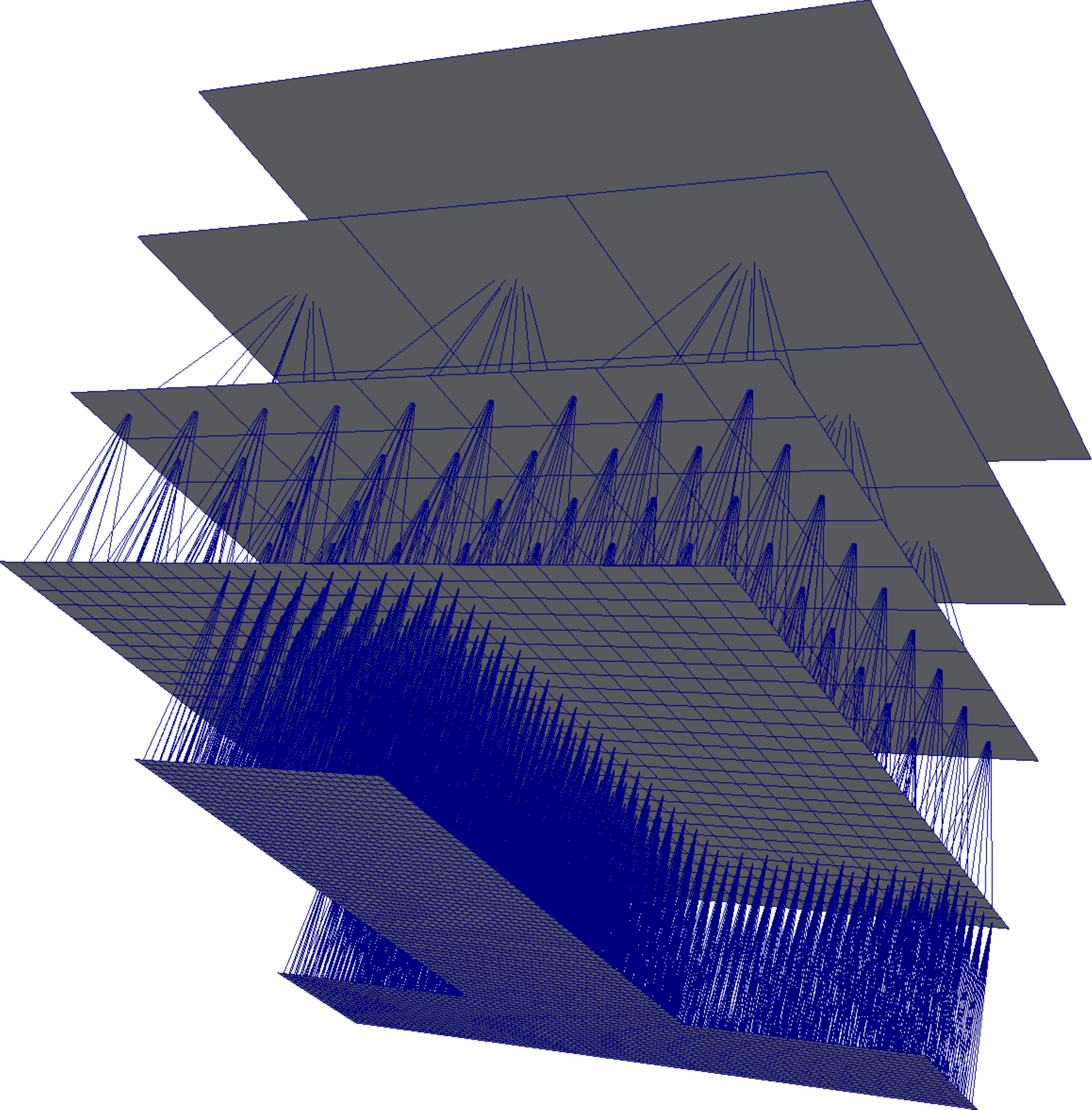}
  \hspace{0.4cm}
  \includegraphics[width=0.3\textwidth]{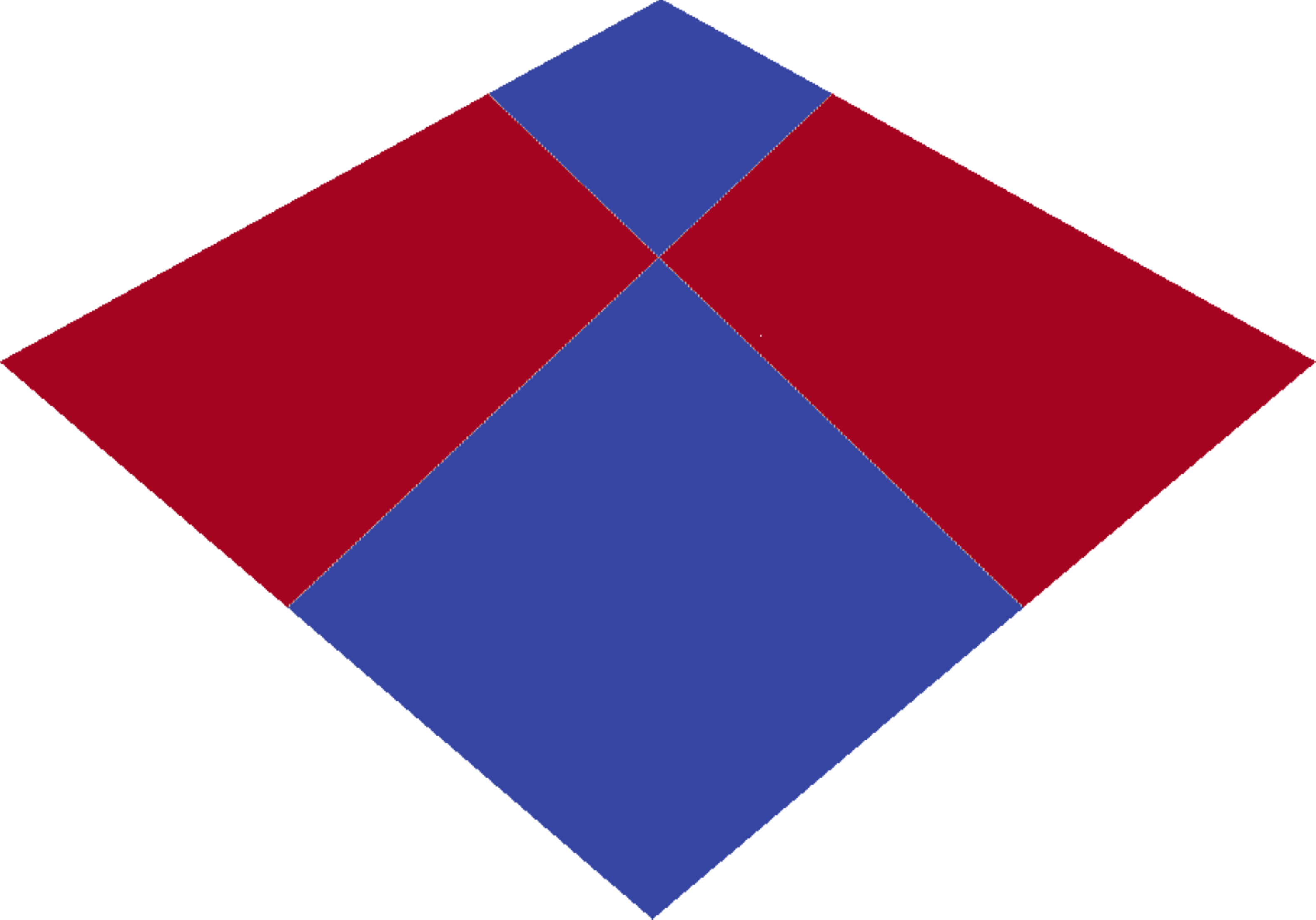}
  \hspace{0.4cm}
  \includegraphics[width=0.3\textwidth]{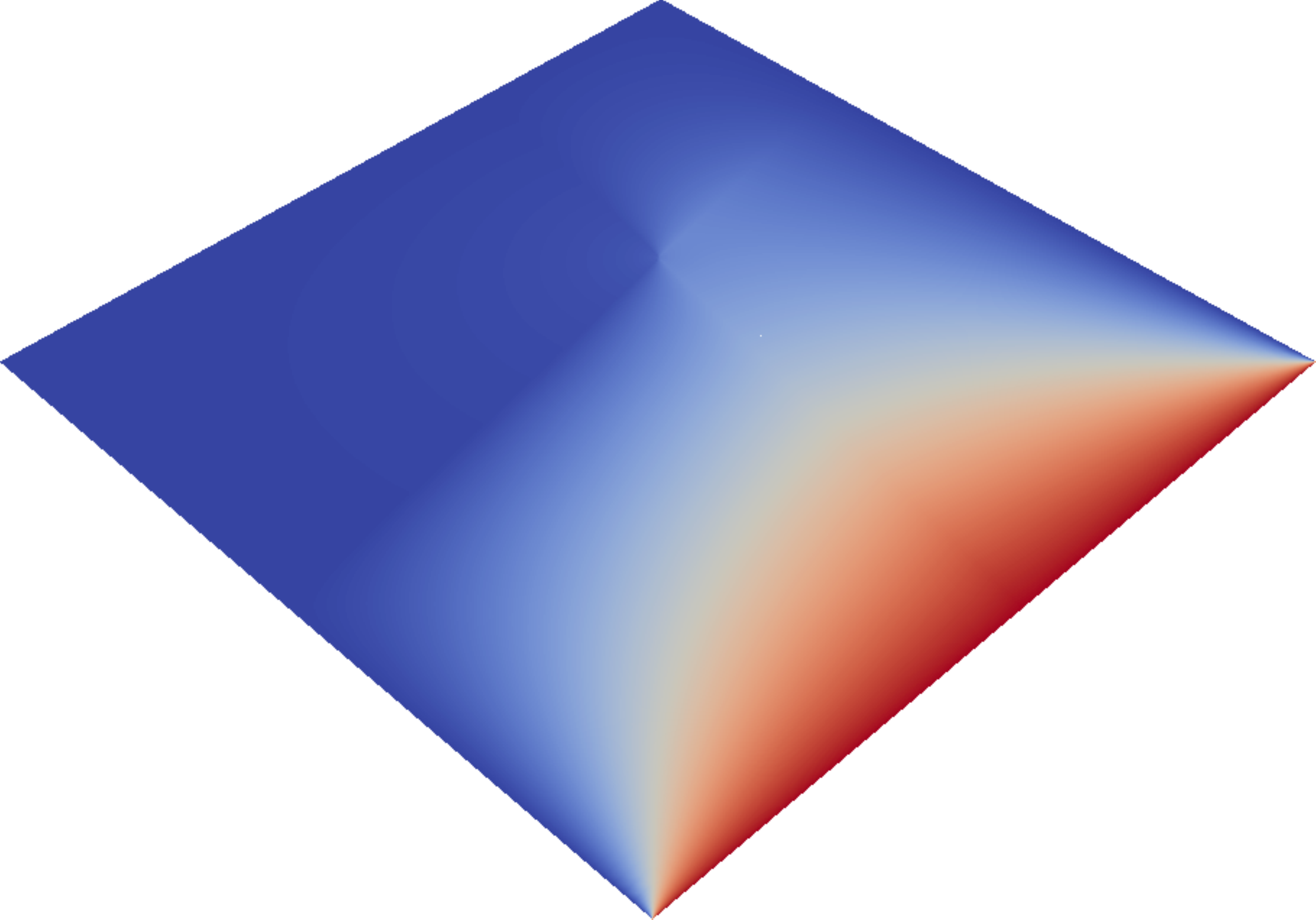}
 \end{center}
 \caption{
  Left:
  A $3 \times 3 $ mesh ($\ell =1$; top layer) serves as first refinement level.
  From here, we construct subsequent refinement levels by subdividing cells
  into $3 \times 3 $ patches.
  This yields a spacetree.
  Multiple Cartesian meshes are embedded into each other.
  Middle: Conductivity (material) parameter setup as used for a stationary heat
  equation solve where the right bottom side of the unit square is heated up.
  A high conductivity in two domain subregions makes the solution (right)
  asymmetric.
  \label{figure:spacetrees:grid-example}
 }
\end{figure}

Our meshing relies upon a spacetree \cite{Weinzierl:11:Peano,Weinzierl:19:Peano}
(Figure \ref{figure:spacetrees:grid-example}):
The computational domain is embedded into a square ($d=2$) or cube ($d=3$)
which yields a (degenerated) Cartesian mesh with one cell
and $2^d$ vertices.
We use cell as generic synonym for cube or square, respectively.
Let the bounding cell have level $\ell = 0$.
It is equidistantly cut into $k$ parts along
each coordinate axis.
We obtain $k^{d}$ child cells having level $\ell = 1$.
The construction continues recursively while we decide per cell individually
whether to refine further or not.
The process creates a cascade of Cartesian
grids $\Omega _{\ell=0}, \Omega _{\ell=1}, \Omega _{\ell=2}, \ldots $.
We count levels the other way round compared to most multigrid literature
\cite{Smith:96:DomainDecomposition,Trottenberg:01:Multigrid}
\replaced{assigning}{They make} the finest grid \deleted{hold} level $\ell =0$.
\deleted{ and assign increasing indices from fine to coarse.}
Our level grids might be ragged: $\Omega _\ell$ is a regular grid covering the
whole domain if and only if all cells on all levels
$\hat \ell < \ell$ are refined.
We use $k=3$.
Choosing three-partitioning is due to
\cite{Weinzierl:11:Peano,Weinzierl:19:Peano} acting as the implementation baseline. 
All of our concepts however apply to bipartitioning, too.

Our code discretises (\ref{equation:PDE}) with $d$-linear Finite Elements. 
Each vertex on each level $\ell $ that is surrounded by $2^d$ cells on level
$\ell $  carries one ``pagoda'', i.e.~bi- or tri-linear shape function.
The remaining vertices are hanging vertices.
Testing shape functions against other functions from the same level 
yields compact $3^d$ stencils.
For this, we make hanging and boundary vertices carry truncated
shape functions but no test functions.
A discussion of Neumann conditions is out of scope.
We therefore may assume that the scaling of the truncated shapes along
the boundary is known.
Due to the spacetree's construction pattern, stencils act on a nodal
generating system over an adaptive Cartesian grid $\Omega _h = \cup _\ell \Omega _\ell $.
If we study (\ref{equation:PDE})  only over the vertices from all levels
that carry a shape function and do not spatially coincide with any other vertex
of the grid from finer levels, we obtain a nodal shape space over an adaptive Cartesian grid
$\Omega _h$.

Let $\ell _\text{max}$ identify the finest mesh, i.e.~the maximum level, while $\ell
_\text{min} \geq 1$ is the coarsest level which holds
degrees of freedom.
$\ell _\text{max} \geq \ell _\text{min}$.
\replaced{For our benchmarking,}{Usually,} $\ell _\text{min}=1$ is appropriate.
\replaced{However,}{, though experience teaches us that} bigger $\ell _\text{min}$
might be reasonable if a problem's solution can't be accurately represented on the coarsest meshes anymore
\added{or performance arguments imply that it is not reasonable to continue to
use multigrid.
In these cases, most codes switch to Krylov methods, direct solvers or algebraic
multigrid with explicit assembly
\cite{Lu:14:HybridMG,May:15:Scalable,Rudi:15:Extreme,Sundar:12:ParallelMultigrid}.
We neglect such solver hybrids and emphasise that our new additive solver allows
us to use rather small $\ell _\text{min}$.
} 
Our subsequent discussion introduces the linear algebra
ingredients for a regular grid corresponding to $\ell _\text{max}$.
The elegant
handling of the adaptive grid is subject of a separate subsection where we
exploit the transition from a generating system into a basis.
Without loss of generality, (\ref{equation:PDE}) is thus discretised into 
\[
 A_{\ell _\text{max}} u_{\ell _\text{max}} = b_{\ell _\text{max}}.
\]

\subsection{Additive and multiplicative multigrid}

Additive multigrid reads as
\begin{equation}
  u_{\ell _\text{max}}  \gets u_{\ell _\text{max}}  +
  \left( \sum _{\ell =
  \ell_\text{min}}^{\ell _\text{max}} \omega _\text{add}( \ell ) P^{{\ell _\text{max}} -\ell }  
    M^{-1}_{\ell }
  R^{{\ell _\text{max}} -\ell } \right)
  \left(b _{\ell _\text{max}} - A _{\ell _\text{max}} u_{\ell _\text{max}} \right),
  \label{equation:additive-mg}
\end{equation}

\noindent
where $M _\ell$ is an approximation to $A _\ell$.
We use the Jacobi smoother $M^{-1}_\ell 
= \diag ^{-1} (A_\ell )$ \added[id=R2]{on all grid levels $\ell $.
No alternative (direct) solver or update scheme is employed on any level}.
The generic prolongation symbol $P$ accepts a solution on a particular
level $\ell -1$ and projects it onto the next finer level $\ell$. 
The exponent indicates repeated application of this inter-grid transfer
operator.
Restriction works the other way round, i.e.~projects from finer to coarser
meshes.
Ritz-Galerkin multigrid \cite{Trottenberg:01:Multigrid} finally yields
$A_\ell = R A_{\ell +1} P$ for $\ell < \ell _{max}$.

For an $\ell $-independent, constant $\omega _\text{add} (\ell )\in ]0,1]$, additive
multigrid tends to become unstable once $\ell _\text{max} - \ell _{min} $ becomes large
\cite{Bastian:98:Additive,Hart:89:FAC,Reps:17:Helmholtz}:
If the fine grid residual $b _{\ell _\text{max}} - A _{\ell _\text{max}} u_{\ell _\text{max}}
$ is homogeneously distributed, the residuals all \replaced{produce similar corrections---effectively attempting to reduce the same error multiple times}{push the solution
into the same direction}. 
Summation of all level contributions then moves the solution too
aggressively into this direction.
A straightforward fix is exponential damping $\omega _\text{add} (\ell ) = \hat
\omega _\text{add} ^{\ell _\text{max} - \ell}$ with a fixed $ \hat  \omega _\text{add} \in
]0,1[$.
If an adaptive mesh is used, $\ell _\text{max} - \ell $ is
ill-suited as there is no global $\ell _\text{max}$ hosting the solution.
We introduce an appropriate, adaptive damping in \cite{Reps:17:Helmholtz} where we make $\ell
_\text{max}$ a per-vertex property.
It is derived from a tree grammar \cite{Knuth:90:AttributeGrammar}.
Such exponential damping, while robust, 
struggles to track global solution effects efficiently once
many mesh levels are used:
The coarsest levels make close to no contribution to the solution.

Multiplicative multigrid is more robust than additive multigrid by construction.
Multiplicative multigrid does not make one residual feed into all level updates in one rush, but
updates the levels one after another. It starts with the finest level.
Before it transitions from a fine level to the next coarsest level, 
it runs some approximate solves (smoothing steps) on the current level \replaced{to}{, and}
\replaced{yield}{yields} a new residual.
We may
assume that the error represented by this residual is smooth.
Yet, the representation becomes rough again on the next \deleted{coarser }level, where we
\replaced{become able}{continue} to smooth it efficiently \added{again}.
Cascades of smoothers act on cascades of frequency bands. 
Multiplicative methods are characterised by the number of the pre- and
postsmoother steps $\mu _{pre}$ and $\mu _{post}$, i.e.~the number of relaxation
steps before we move to the next coarser level (pre) or next finer level (post), respectively.
The multiplicative multigrid solve closest to the additive scheme is a
$V(0,\mu _{post})$-cycle, i.e.~a scheme without any presmoothing and $\mu
_{post}$ postsmoothing steps.
\added{
 Different to real additive multigrid, the effect of smoothing on a level $\ell$
 here does feed into the subsequent smoothing on $\ell +1$.
} 
\replaced{Since}{However,} $\mu _{pre} =
0$ yields no classic multiplicative scheme\replaced{---}{, as}the resulting
solver does not smooth prior to the coarsening\replaced{---we }{. In
practice, it works nevertheless. We}
conclude that the $V(\mu _{pre}=1,0)$-cycle thus is the \replaced{(robust)
multiplicative scheme most similar}{closest cousin} to an additive scheme.
The multiplicative  
two-grid scheme with exact coarse grid solve reads
\begin{eqnarray}
 u_{\ell _\text{max}} & \gets & 
  PA_{\ell _\text{max}-1}^{-1}R(
   b_{\ell _\text{max}} - A_{\ell _\text{max}}
   \left[ 
    u_{\ell _\text{max}} + \omega _{{\ell_\text{max}}} M^{-1}_{{\ell_\text{max}}}  (
    b_{\ell _\text{max}} - A_{\ell _\text{max}} u_{\ell _\text{max}} ) \right]
  )
  \nonumber 
  \\
  && + 
   \left[ 
    u_{\ell _\text{max}} + \omega _{{\ell_\text{max}}} M^{-1}_{{\ell_\text{max}}}  (
    b_{\ell _\text{max}} - A_{\ell _\text{max}} u_{\ell _\text{max}} ) \right].
  \label{equation:multiplicative-two-grid}
\end{eqnarray}

\subsection{Multigrid on hierarchical generating systems}

%
%
Early work on locally adaptive multigrid
(see for example \cite{Brandt:77:Multi,McCormick:86:FAC} as well as the
historical notes in \cite{Lee:04:AFAC}) already relies on block-regular
Cartesian grids \cite{Dubey:16:SAMR} and nests the geometric grid resolutions
into each other.
The coarse grid vertices spatially coincide with
finer vertices where the domain is refined.
This yields a hierarchical generating system rather than a basis.

%
%
The fast adaptive composite (FAC) method \cite{Hart:86:FAC,Hart:89:FAC}
describes a multiplicative multigrid scheme over this hierarchical system: 
We start from the finest grid, determine the residual
equation there, smooth, re-compute the residual and
restrict it to the next coarser level. 
It continues recursively.
\replaced{As we compute corrections using residual equations,
this is a multigrid scheme.}{As we rely on residuals, this is a multigrid scheme.}
As we \added{sequentially} smooth and then recompute the residual, it is a
multiplicative scheme.
Early FAC papers \replaced{operate on the assumption of}{orbit around} a small set of reasonable fine grids and leave it
open to the implementation which iterative scheme to use.
Some explicitly speak of FAC-MG if a multigrid cycle is used as the iterative smoother per
level.
We may refrain from such details and consider \replaced{fast adaptive composite grid}{FAC} as a multiplicative scheme
overall which can be equipped with simple single-level smoothers.

%
%
The first \replaced{fast adaptive composite grid}{FAC} papers \cite{Hart:86:FAC}
acknowledge difficulties for operators along the resolution transitions.
While we discuss an elegant handling of these difficulties in
\ref{section:ingredients:HTMG}, \replaced{fast adaptive composite grid}{FAC} traditionally addresses them through a
top-down traversal \cite{Hart:89:FAC}:
The cycle starts with the coarsest grid, and then uses the updated solution to
impose Dirichlet boundary conditions on hanging nodes on the next finer level.
This inversion of the grid level order continues to yield a multiplicative
scheme as updates on coarser levels immediately propagate down and as all steps
are phrased as residual update equations.

%
%
%
FAC relies on spatial discretisations that are conceptually close to our
spacetrees.
Both approaches thus benefit from structural simplicity: 
As the grid segments per level are regular, solvers (smoothers) for regular
Cartesian grids can be (re-)used.
As the grid resolutions are aligned with each other, hanging nodes
can be assigned interpolated values from the next coarsest grid with a
geometrically inspired prolongation.
As all grid entities are cubes, squares or lines, all operators exhibit 
tensor-product structure.
\replaced{ 
FAC's hierarchical basis differs from textbook multigrid 
\cite{Trottenberg:01:Multigrid} for adaptive meshes:
}{
Through FAC's hierarchical basis approach, it does not share two properties with
standard, vanilla multigrid \cite{Trottenberg:01:Multigrid} if applied to
adaptive meshes:
}
The fine grid smoothers do not address the real fine grid, but only separate segments 
that have the same resolution.
The transition from fine to coarse grid does not imply that the number of
degrees of freedom decreases.
Rather, the number of degrees of freedom can increase if the
finer grid accommodates a very localised AMR region.
It is obvious that this poses challenges for parallelisation.

%
%
We can mechanically rewrite multiplicative FAC into an additive version.
The hierarchical generating system renders this endeavour straightforward.
However, plain additive multigrid on a FAC data structure again yields a
non-robust\deleted{, overshooting} solver \added{that tends to overcorrect}
\cite{Hart:89:FAC,Reps:17:Helmholtz}.
There are multiple approaches to tackle this:
Additive multigrid with exponential damping removes oscillations \added{from the
solution} at the cost of multigrid convergence behaviour.
\replaced{Bramble, Pasciak and Xu's scheme (BPX)}{BPX} \cite{Smith:96:DomainDecomposition} is the most popular variant where we
accept the non-robustness and use the additive scheme solely as a preconditioner.
To make this preconditioner cheap, BPX traditionally neglects (Ritz-Galerkin)
coarse grid operators. 
Instead, it replaces the $M_\ell^{-1}$ in (\ref{equation:additive-mg}) with a
diagonal matrix for the correction equations, where the diagonal matrix is
scaled such that it mimics the Laplacian.
The hierarchical basis approach starts from the observation that the
instabilities within the generating system are induced by spatially
coinciding vertices.
Therefore, it drops all vertices (and their shape functions) on one level that
coincide with coarser vertices.
The asynchronous \replaced{fast adaptive composite grid}{FAC} (AFAC) solver family
finally modifies the operators to anticipate overshooting.
We may read BPX as particular modification of additive multigrid and
AFAC as a generalisation of BPX \cite{Jimack:11:Asynchronous}.

\subsection{HTMG and FAS on spacetrees}
\label{section:ingredients:HTMG}

Though the implementation of multigrid on adaptive meshes is, in principle,
straightforward, implementational complexity arises along resolution
transitions.
Weights associated to the vertices
change \deleted{their} semantics once we compare vertices on a level $\ell $ which are surrounded by
refined spacetree cells to vertices on \replaced{that}{this} level which belong to the fine
grid:
The latter carry a nodal solution representation, i.e.~a scaling of the
Finite Element shape functions, while the former carry correction weights.
In classic multigrid starting from a fine grid and then traversing correction
levels, it is not straightforward how to handle the vertices \replaced{on the border }{in-}between a fine
grid region and a refined region within one level.
They carefully have to be separated
\cite{Reps:17:Helmholtz,Weinzierl:18:BoxMG}.

One elegant solution to address this ambiguity relies on full
approximation storage (FAS) \cite{Trottenberg:01:Multigrid}.
Every vertex holds a nodal solution representation.
If two vertices $v_{\ell}$ and $v_{\ell +1}$ from two levels spatially coincide,
the coarser vertex holds a copy of the finer vertex:
In areas where two grids overlap, the coarse grid \replaced{additionally holds}{is} the injection $u_{\ell} =
I u_{\ell +1}$ of the fine grid.
This definition exploits the regular construction pattern of spacetrees. 
Vertices in refined areas now carry a correction equation plus the
injected solution rather than a sole correction.
The injection couples the fine grid problem with its coarsened representation
and makes this representation consistent with the fine grid
problem on adjacent meshes which have not been refined further.
In the present paper, we use FAS exclusively to resolve the semantic
\deleted{dis}ambiguity that arises for vertices at the boundary between the fine grid and
a correction region on one level;
further potential such as $\tau $-extrapolation 
\cite{Rude:87:Multiple}
 or the application to nonlinear PDEs, i.e.~$\epsilon
= \epsilon (u)$ in (\ref{equation:PDE}),
is not exploited.

Our code relies on \replaced{hierarchical transformation multigrid (HTMG)}{HTMG}
\cite{Griebel:90:HTMG} for the implementation of \replaced{the full approximation storage scheme}{FAS}.
It also relies on the
assumption/approximation that all of our operators can be approximated by
Ritz-Galerkin multigrid $R A_{\ell + 1} P = A_{\ell }$.
Injection $u_{\ell} = I u_{\ell +1}$ allows us to rewrite each and every
nodal representation into its hierarchical representation $\hat u _{\ell} = (id
- PI) u _{\ell}$.
A hierarchical residual $\hat r$ is defined in the expected way.
This elegantly yields the modified \deleted{FAS }multigrid equation 
when we switch from the correction equation to
\begin{eqnarray}
  A_\ell \left( u_\ell + c_\ell \right) & = & 
  A_\ell u_\ell + A_\ell  c_\ell = 
  R \hat r_{\ell +1} 
  \nonumber \\
  & = &
  R \left( b_{\ell +1} - A_{\ell +1}(u_{\ell +1}-PI_{\ell +1}) \right) = 
  R \left( b_{\ell +1} - A_{\ell +1} \hat u_{\ell +1} \right),
  \label{equation:fas-with-htmg}
\end{eqnarray}

\noindent
i.e.~per-level equations
\[
  A_{\ell } u_{\ell } = \left\{
   \begin{array}{cl}
    b_{\ell } & \mbox{on the fine grid (regions)} \\
    b_{\ell } = R \hat r_{\ell +1} & \mbox{on the coarse grid (regions) with }
    \hat r_{\ell +1} = b_{\ell +1} - A_{\ell +1} \hat u_{\ell +1}.
   \end{array}
  \right.
\]

\noindent
To the smoother, $u_{\ell }$ resulting from the injection serves as 
the initial guess.
Subsequently it determines a correction
$c_\ell $.
This correction feeds into the multigrid prolongation.

Equation (\ref{equation:fas-with-htmg}) clarifies that the right-hand side of
\replaced{this full approximation storage}{FAS} does not require a complicated
calculation:
We ``simply'' have to determine the hierarchical representation $\hat u$ on the
finer level, compute a hierarchical residual $\hat r$ on this level (which uses
the smoother's operator), and restrict this value to the coarse grid's right-hand
side.

\subsection{BoxMG and algebraic-geometric multigrid hybrids}

BoxMG is a geometrically inspired algebraic 
technique \cite{Dendy:82:BlackboxMG,Dendy:10:BoxMgBy3,Yavneh:12:NonnsymBoxMg}
to determine inter-grid transfer operators
\replaced{on meshes that are embedded into each other.}{
In a spacetree context, each fine grid is split up along the next coarser
level's grid lines.}
We assume that the prolongation from a coarse vertex maps onto the nullspace of
the fine grid operator.
However, BoxMG does not examine the ``real'' operator.
Instead, it studies an operator which is collapsed along the coarse grid
\replaced{cell}{level} boundaries.

All fine grid points are classified into $c$-points (coinciding
spatially with coarse grid points of the next coarser level), 
$\gamma $-points which coincide with the faces of the next coarser levels
and $f$-points.
Prolongation and restriction are \replaced{defined as}{} the identity on $c$-points.
Along $\gamma $-points, we collapse the stencil:
If $\gamma $ members reside on a face with normal $n$, the stencil is
accumulated (lumped) along the $n$ direction.
The result contains only entries along the non-$n$ directions.
Higher dimensional collapsing can be constructed iteratively.
We solve $\tilde APe=0|_\gamma $---$\tilde A$ stems from the collapsed
operators---along these $\gamma $-points where $e$ is the characteristic 
\replaced[id=R2]{vector for a vertex}{vertex vector} on the coarse grid, i.e.~holds one entry $1$ and zeroes everywhere else.
Finally, we solve $APe=0 |_f$ for the remaining points. 
No two $f$-points separated by a coarse grid line are coupled to each other
anymore.

In our previous work~\cite{Weinzierl:18:BoxMG}, we detail how to store BoxMG's operators as
well as all Ritz-Galerkin operators which typically supplement BoxMG within the
spacetree.
This yields a hybrid scheme in two ways:
\added{On the one hand,}
BoxMG itself is a geometrically inspired way to construct algebraic inter-grid
transfer operators.
Storing the entries within the spacetree \added{on the other hand} allows for a
``matrix-free''
\deleted{(i.~e. one that does not globally assemble the full matrix)}
implementation \replaced{where no explicit matrix structure is held but
all matrix entries are embedded into the mesh.
With an on-the-fly compression of entries relative to rediscretised
values \cite{Weinzierl:18:BoxMG}, we effectively obtain the total memory
footprint of a matrix-free scheme.}
{with explicit matrix storage.}

\section{Additively damped AFAC solvers with FAS and HTMG}
\label{section:integration}

With our ingredients and observations at hand, our research agenda reads as
follows: 
We first introduce our additive multigrid scheme which avoids oscillations
without compromising on the convergence speed.
Secondly, we discuss two operators suited to realise our scheme.
Finally, we contextualise this idea and show that the new solver actually
belongs into the family of AFAC solvers.

\subsection{An additively damped additive multigrid solver}

Both additive and multiplicative multigrid sum up all the levels' corrections.
Multiplicative multigrid is more stable than additive---it does not
overshoot---as each level eliminates error modes tied to its resolution.
In practice, we cannot totally separate error modes, and we cannot assume
that a correction on level $\ell $ does not introduce a new error on level $\ell +1$.
Multigrid solvers thus often use postsmoothing.
Once we ignore this multiplicative lesson, the simplest class of
multiplicative solvers is $V(\mu _{pre}=1,0)$.

We start with our recast of
the multiplicative $V(1,0)$ two-grid cycle
(\ref{equation:multiplicative-two-grid}) into an additive formulation
(\ref{equation:additive-mg}).
Our objective is to quantify additive multigrid's \replaced{over-correction}{``too much'' of a correction}
relative to its multiplicative cousin.
For this, we compare the multiplicative two-grid scheme \added[id=R3]{$u_{\ell _\text{max},\text{mult}}$}
(\ref{equation:multiplicative-two-grid}) to the two level additive scheme with an exact solve on the coarse level

%
%
\[
 u_{\ell _\text{max},\text{add}}^{(n+1)} =  
  PA_{\ell _\text{max}-1}^{-1}R(
   b_{\ell _\text{max}} - A_{\ell _\text{max}}
    u_{\ell _\text{max}} ^{(n)}
  )
  + 
   \left[ 
    u_{\ell _\text{max}}^{(n)} + \omega _{\ell _\text{max}} M^{-1}_{\ell_\text{max}}  (
    b_{\ell _\text{max}} - A_{\ell _\text{max}} u_{\ell _\text{max}}^{(n)} ) \right].
\]

\noindent
The difference is 
\begin{eqnarray}
  u_{\ell _\text{max}, \text{mult}} ^{(n+1)} - u_{\ell _\text{max}, \text{add}} ^{(n+1)}  
  & = & 
  PA_{\ell _\text{max}-1}^{-1}R(
   b_{\ell _\text{max}} - A_{\ell _\text{max}}
   \left[ 
    u_{\ell _\text{max}} ^{(n)} + \omega _{{\ell_\text{max}}} M^{-1}_{{\ell_\text{max}}}  (
    b_{\ell _\text{max}} - A_{\ell _\text{max}} u_{\ell _\text{max}} ^{(n)} ) \right]
  )
  \nonumber
  \\
  && -
  PA_{\ell _\text{max}-1}^{-1}R(
   b_{\ell _\text{max}} - A_{\ell _\text{max}}
    u_{\ell _\text{max}} ^{(n)}
  ).
  \label{equation:integration:difference-add-mult}
\end{eqnarray}

\noindent
The superscripts $^{(n)}$ and $^{(n+1)}$ denote old and respectively new iterates of a
vector. 
We continue to omit it from here where possible.

Starting from the additive rewrite of the $V(1,0)$ multiplicative two-level
scheme, we intend to express multiplicative multigrid as an additive scheme.
This is a
popular endeavour as additive multigrid tends to \replaced{more readily show\deleted[id=R2]{s} improved performance on large scale parallel implementation}{scale better}.
There is no close-to-serial coarse grid solve.
There is no coarse grid bottleneck in an Amdahl sense.
Multiplicative multigrid however tends to converge faster and is more robust.
Different to popular approaches such as Mult-additive
\cite{Yang:14:Reducing,Wolfson:19:AsynchronousMultigrid}, our approach 
\replaced{does not aim to achieve the exact}{gives up on the idea to achieve the}  convergence rate of multiplicative multigrid.
\replaced{Instead, we aim to mimic the robustness of multiplicative multigrid in
an additive regime---i.~e. allow additive multigrid to successfully converge across
a wider range of setups.} {Instead, we save ``solely'' the robustness over into
the additive regime.} Our hypothesis is that any gain in concurrency will eventually outperform efficiency improvements on future machines.
A few ideas guide our agenda:

\begin{idea}
  We add an additional one-level term to our additive scheme
  \replaced{which compensates for additives overly aggressive updates compared to}{which
  imitates} multiplicative \added{$V(1,0)$ }multigrid.
\end{idea}

\noindent
This idea \replaced{describes the rationale behind}{circumscribes}
(\ref{equation:integration:difference-add-mult}) where we \deleted{already}
stick to a two-grid formalism.
Our strategy next is to find an approximation to 
\begin{eqnarray}
-PA_{\ell _\text{max}-1}^{-1}RA_{\ell _\text{max}}\omega _{\ell _\text{max}}
M^{-1}_{{\ell_\text{max}}}  ( b_{\ell _\text{max}} - A_{\ell _\text{max}} u_{\ell _\text{max}} )
  \label{equation:damping-term}
\end{eqnarray}
from (\ref{equation:integration:difference-add-mult}) 
such that we obtain a modified additive two-grid scheme which, on the one hand,
mimics multiplicative stability and, on the other hand, is cheap. 
For this, we read the difference term as an auxiliary solve.

\begin{idea}
  We approximate the auxiliary term (\ref{equation:damping-term}) with a single
  smoothing step.
\end{idea}

\noindent
The approach yields a per-level correction
\begin{eqnarray}
-P
\omega _{\ell-1} \tilde M_{\ell _\text{max}-1}^{-1}
RA_{\ell _\text{max}}\omega _{\ell _\text{max}}
M^{-1}_{{\ell_\text{max}}}  ( b_{\ell _\text{max}} - A_{\ell _\text{max}} u_{\ell _\text{max}}
).
\label{equation:adaFACx-auxiliary-equation}
\end{eqnarray}

\noindent
We use the tilde to denote the auxiliary solves.
Following Idea 1, this is a per-level correction:
When we re-generalise the scheme from two grids to multigrid 
(by a recursive expansion of
$A_{\ell _\text{max}-1}^{-1}$ within the original additive formulation),
we do not further expand the correction (\ref{equation:damping-term}) or
(\ref{equation:adaFACx-auxiliary-equation}).
This implies another error which we accept in return for a simplistic correction
term without additional synchronisation or data flow between levels.

\begin{idea}
  The damping runs asynchronously to the actual solve. It is another additive
  term computed concurrently to each correction equation.
\end{idea}

\noindent
Using $A_{\ell _\text{max}} \omega _{\ell _\text{max}} M^{-1}_{{\ell_\text{max}}}$
adds a
sequential ingredient to the damping term.
A fine grid solve must be finished before it can enter the auxiliary equation.
This reduces concurrency.
Therefore, we propose to merge this preamble smoothing step into the
restriction.
typically uses a simple aggregation/restriction
operator and then improves it by applying a smoother.
It is also similar to Mult-additive \cite{Yang:14:Reducing}, which constructs
inter-grid transfer operators that pick up multiplicative pre- or postsmoothing
behaviour.
We apply the smoothed operator concept to the
restriction $\tilde{R}= \omega R A M^{-1}$, and end up with a wholly additive correction term
\begin{equation}
 -\tilde \omega P \tilde M_{\ell _\text{max}-1}^{-1} \tilde R.
 \label{equation:integration:damping-equation}
\end{equation}

\begin{idea}
  We geometrically identify the auxiliary coarse grid levels with the actual multilevel grid
  hierarchy.
  All resolution levels integrate into the spacetree.
\end{idea}

\noindent
$\tilde M$ and $\tilde A$ are auxiliary operators but act on mesh levels which
we hold anyway.
With the spacetree at hand, we finally unfold the two-grid scheme into
\begin{eqnarray}
  u_{\ell _\text{max}} & \gets & u_{\ell _\text{max}}  
  + 
  \left( \sum _{\ell =
  \ell_\text{min}}^{\ell _\text{max}} \omega _\text{add}( \ell ) P^{{\ell _\text{max}} -\ell }  
    M^{-1}_{\ell }
  R^{{\ell _\text{max}} -\ell } \right)
  \left(b _{\ell _\text{max}} - A _{\ell _\text{max}} u_{\ell _\text{max}} \right) 
  \nonumber
  \\
  && 
  \phantom{u_{\ell _\text{max}}}
  -
  \left( \sum _{\ell =
  \ell_\text{min}}^{\ell _\text{max}}\tilde{ \omega} _\text{add}( \ell ) P^{{\ell _\text{max}} -\ell }  
    \tilde M^{-1}_{\ell }
  \widetilde{R}^{{\ell _\text{max}} -\ell } \right)
  \left(b _{\ell _\text{max}} - A _{\ell _\text{max}} u_{\ell _\text{max}} \right),
  \label{equation:adaFACx}
\end{eqnarray}

\noindent
where we set, without loss of generality, $M^{-1}_{\ell _\text{max}-1}=0$.
This assumes that no level coarser than $\ell _\text{min}$ hosts any degree of
freedom.

Algorithms in standard AFAC literature present all levels as correction levels.
That is, a global residual is computed on the composite grid
and then restricted to construct the right hand side
of error equations on all grid resolutions.
This includes the finest grid level.
Here instead we use standard multigrid convention
and directly smooth the finest grid level
(Algorithms \ref{algorithm:adAFAC:Jac} and \ref{algorithm:adAFAC:PI}).
We only restrict the residual to coarse grid levels.

\added{
 The four ideas align with the three key concepts from the introduction: 
 We stick to a geometric grid hierarchy and then also reuse this hierarchy for
 additional equation terms.
 We stick to an additive paradigm and then also make additional equation terms
 additive.
 We stick to a geometric-algebraic mindset.
}

\begin{algorithm}
  \caption{Blueprint of one cycle of the our adAFAC-Jac without AMR. 
  $R^{i}$ or $P^{i}$ denote the recursive application of the restriction or
  prolongation, respectively.
  $\widetilde{R}^{i}$ applies $R$ $i-1$ times,
  followed by an application of one smoothed operator.
  \label{algorithm:adAFAC:Jac}
  }
    \begin{algorithmic}
    \Function{adAFAC-Jac}{}
    	\State $r _{\ell _\text{max}} \gets b_{\ell _\text{max}} - A_{\ell _\text{max}} u_{\ell
    	_{max}} $
    	\ForAll {$\ell_\text{min} \leq \ell < \ell _\text{max}$}
    		\Comment{Restrict fine grid residual to grid levels $\ell $} 
    		\State ${b_\ell } \gets R^{\ell _\text{max} - \ell} r _{\ell _\text{max}}$
    		\State ${\widetilde{b}_{\ell }} \gets \widetilde{R}^{\ell _\text{max} - \ell} r _{\ell _\text{max}}$ 
    		\Comment{Additional restriction residual into additional grid space}
    	\EndFor	
    	\ForAll {$\ell_\text{min} < \ell < \ell _\text{max}$}
    		\State ${{c_\ell } \gets 0};$ ${{\widetilde{c}_{\ell-1} } \gets 0}$
    		\Comment Initial ``guess'' of correction and damping
    		\State $c_\ell \gets $ \Call{jacobi}{$A_\ell c_\ell = b_\ell, \omega $} 
    		\Comment Iterate of correction equation stored in $c_\ell $
    		\State $\hat c_{\ell -1} \gets $ \Call{jacobi}{$A_{\ell -1}
    		\widetilde{c}_{\ell -1} = \widetilde{b}_{\ell -1} , \widetilde{\omega} $} 
    		\Comment Iterate of corresponding damping
    		equation (one level coarser)
    	\EndFor
        \State ${{c_{\ell _\text{min}} } \gets 0};$ ${{\widetilde{c}_{\ell _\text{min}-1} } \gets 0}$
        \Comment Initial ``guesses'' on coarsest level
        \State $c_{\ell _\text{min}} \gets $ \Call{jacobi}{$A_{\ell _\text{min}} c_{\ell _\text{min}} = b_{\ell _\text{min}}, \omega $} 
   		\Comment Iterate of correction equation. No damping active on coarsest level
        \State ${{c_{\ell _\text{max}} } \gets 0};$ ${{\widetilde{c}_{{\ell _\text{max}}-1} } \gets 0}$
        \State $c_{\ell _\text{max}} \gets $ \Call{jacobi}{$A_{\ell _\text{max}} u_{\ell _\text{max}} = b_f, \omega $} 
   		\Comment Finegrid correction
   		\State $\hat c_{{\ell _\text{max}} -1} \gets $ \Call{jacobi}{$A_{{\ell _\text{max}} -1}
    		\widetilde{c}_{{\ell _\text{max}} -1} = \widetilde{b}_{{\ell _\text{max}} -1} , \widetilde{\omega} $} 
        \Comment Damping of finest grid correction (one level coarser)
  		\State $ u_{\ell _\text{max}} \gets u_{\ell _\text{max}} + c_{\ell
  		_\text{min}} + \sum^{{\ell _\text{max}}}_{\ell = {\ell _{\text{min} -1}}} {P^{\ell _\text{max} - \ell}{c_\ell }
  		- P^{\ell _\text{max} - (\ell-1)}{\widetilde{c}_{\ell-1}}}$ 
    \EndFunction
    \end{algorithmic}
\end{algorithm}

\subsection{Two damping operator choices}

\begin{figure}
 \begin{center}
  \includegraphics[width=0.8\textwidth]{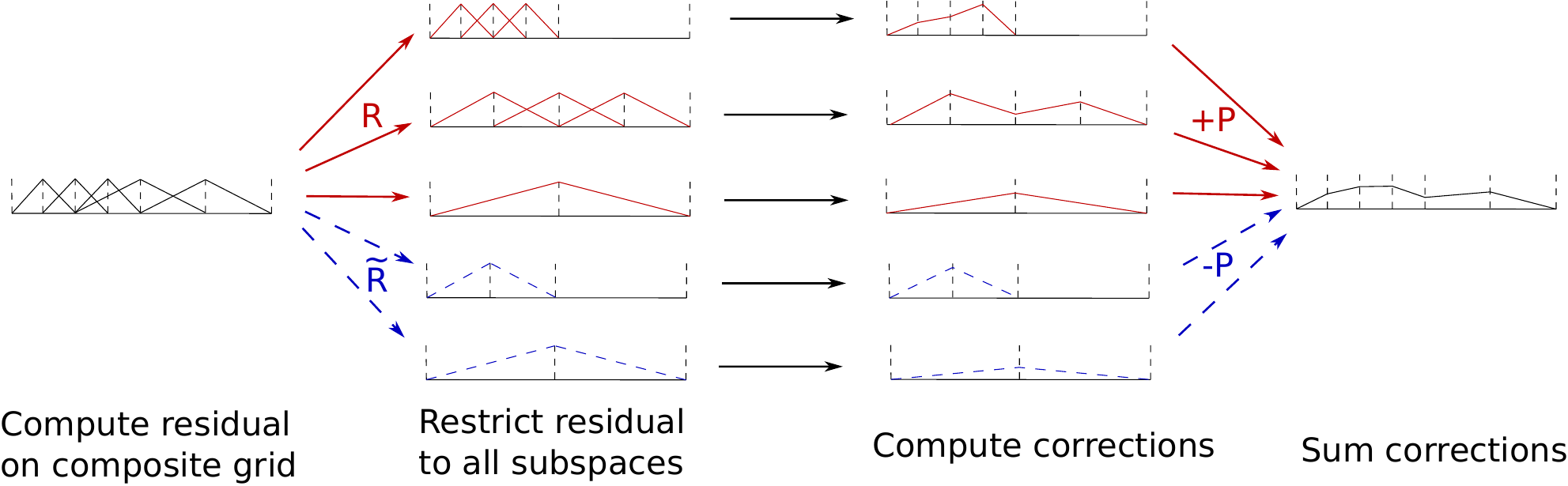}
 \end{center}
 \caption{
  Data flow overview of adAFAC-JAC. 
  Solid red lines denote
  traditional subspaces within additive correction equations, dashed blue lines correspond to
  auxiliary equations that damp the existing correction equations.
  \label{figure:ingredients:adaFAC-JAC-flow}
 }
\end{figure}

\begin{figure}
 \begin{center}
  \includegraphics[width=0.8\textwidth]{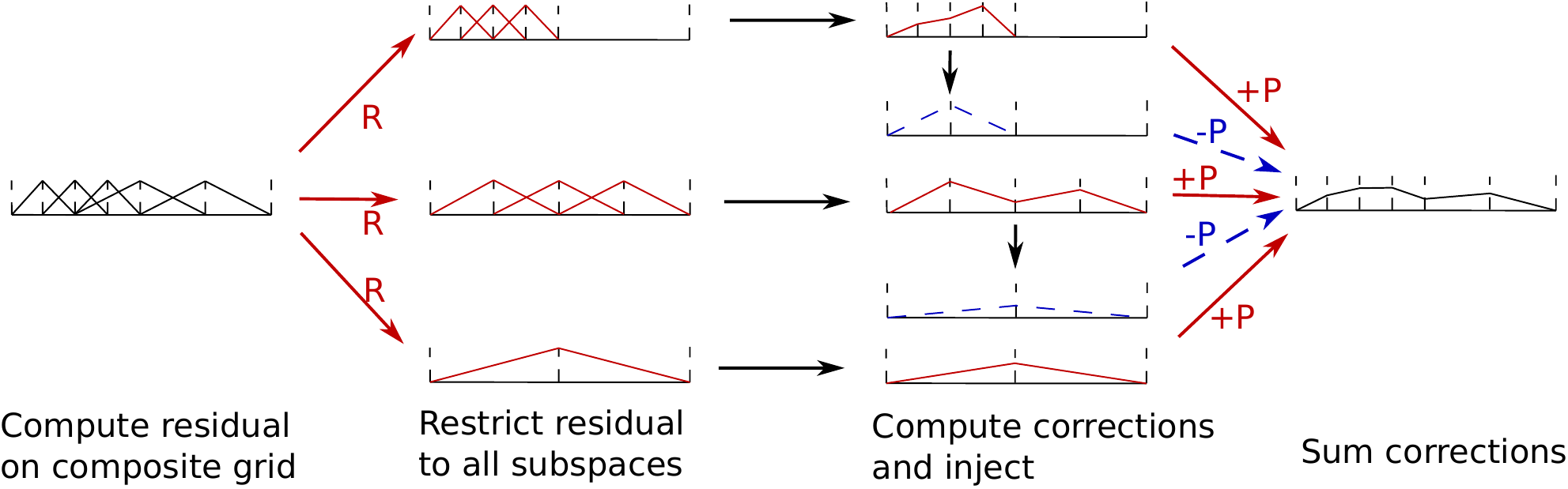}
 \end{center}
 \caption{
  Data flow overview of adAFAC-PI.  
  Solid red lines denote
  traditional subspaces within additive correction equations, dashed blue lines correspond to
  auxiliary equations that damp the existing correction equations.
  \label{figure:ingredients:adaFAC-PI-flow}
 }
\end{figure}

It is obvious that the effectiveness of the approach depends on a proper 
construction of (\ref{equation:integration:damping-equation}). 
We propose two variants.
Both \replaced{are based upon the assumption}{accept} that 
smoothed inter-grid transfer operators yield 
better operators than standard bi- and trilinear operators
(and obviously naive injection or piecewise constant interpolation) 
\cite{Tuminaro:00:Parallel,Vanvek:96:Algebraic,Vanvek:95:Fast}.
Simple geometric transfer operators fail to capture complex solution behaviour 
\cite{Press:91:Multigrid,Kouatchou:00:Optimal,Bjorgen:17:Numerical} for
non-trivial $\epsilon $ choices in (\ref{equation:PDE}).

Let $\epsilon $ in (\ref{equation:PDE}) be one.
We observe that a smoothed operator \added{$M^{-1}A_{\epsilon = 1}P$}
derived from bilinear interpolation \added{$P$ using a Jacobi smoother $M^{-1} = diag(A)^{-1}$} for
three-partitioning \replaced{corresponds to the stencil}{equals}

{\footnotesize
\[
\begin{bmatrix}
-0.0139 & -0.0417   & -0.0833   & -0.0972   & -0.083    & -0.0417   & -0.0139 \\
-0.0417 & 0 & 0 & 0.0833    & 0 & 0 & -0.0417 \\
-0.0833 & 0 & 0 & 0.167 & 0 & 0 & -0.0833 \\
-0.0972 & 0.0833    & 0.167 & 0.444444444   & 0.167 & 0.0833    & -0.0972 \\
-0.0833 & 0 & 0 & 0.167 & 0 & 0 & -0.0833 \\
-0.0417 & 0 & 0 & 0.0833    & 0 & 0 & -0.0417 \\
-0.0139 & -0.0417   & -0.0833   & -0.0972   & -0.0833   & -0.0417   & -0.0139
\end{bmatrix}.
\]
}

\noindent
\added{Here the stencil is a restructured row of the full operator.}
Assuming $\epsilon =1$ is reasonable as the term $A_{\ell} M_{\ell}^{-1}$
or $M_{\ell-1}^{-1}RA_{\ell}$, respectively,
enters the auxiliary restriction.
Such an expression removes the impact of $\epsilon$\deleted{---it yields the
Laplacian---} on all elements with non-variable $\epsilon$. 
Assuming $\epsilon $ is reasonably smooth, we 
neglect only small perturbations in the off-diagonals of the system matrix.

\begin{algorithm}
  \caption{Blueprint of our adAFAC-PI without AMR.
  $R^{i}$ or $P^{i}$ denote the recursive application of the single level restriction or
  prolongation, $R$ or $P$, respectively.
  $I$ is the injection operator.
  \label{algorithm:adAFAC:PI}
  }
    \begin{algorithmic}
    \Function{adAFAC-PI}{}
    	\State $r _{\ell _\text{max}} \gets b_{\ell _\text{max}} - A_{\ell _\text{max}} u_{\ell
    	_{max}} $
    	\ForAll {$\ell_\text{min} \leq \ell < \ell _\text{max}$}
    		\State ${b_\ell } \gets R^{\ell _\text{max} - \ell} r _{\ell _\text{max}}$
    		\Comment{Restrict fine grid residual to coarser levels. $b _{\ell _\text{max}}$ remains untouched} 
    	\EndFor
    	\ForAll {$\ell_\text{min} < \ell < \ell _\text{max}$}
    		\State ${{c_\ell } \gets 0};$ ${{\widetilde{c}_\ell } \gets 0}$
    		\Comment Initial ``guesses'' for corrections
    		\State $c_\ell \gets $ \Call{jacobi}{$A_\ell c_\ell = b_\ell, \omega $} 
    		\Comment Iterate of correction equation stored in $c_\ell $
    		\State ${\widetilde{c}_\ell } \gets  PI{c_\ell }  $
    		\Comment Computation of localised damping for $c_\ell $
    	\EndFor
        \State ${{c_{\ell _\text{min}}} \gets 0};$ ${{\widetilde{c}_{\ell _\text{min}}} \gets 0}$
   		\State $c_{\ell _\text{min}} \gets $ \Call{jacobi}{$A_{\ell _\text{min}} c_{\ell _\text{min}} = b_{\ell _\text{min}},
   		\omega $}
   		\Comment No auxiliary damping for coarsest level
        \State ${{c_{\ell _\text{max}} } \gets 0};$ ${{\widetilde{c}_{\ell _\text{max}} } \gets 0}$
          \Comment Initial ``guess'' for correction on finest grid
        \State $c_\ell \gets $ \Call{jacobi}{$A_\ell u_{\ell _\text{max}} = b_f, \omega $} 
    		\Comment Finegrid update
    		\State ${\widetilde{c}_\ell } \gets  PI{c_\ell }  $
  		\State $ u_{\ell _\text{max}} \gets u_{\ell _\text{max}} + c_{\ell
  		_\text{min}} + \sum^{{\ell _\text{max }}}_{\ell = {\ell _{\text{min}-1}}} {P^{\ell _\text{max} - \ell}{c_\ell }
  		- P^{\ell _\text{max} - \ell}{\widetilde{c}_\ell}}$ 
    \EndFunction
    \end{algorithmic}
\end{algorithm}

This motivates us to introduce two modified, i.e.~smoothed restriction operators
$\tilde R$:
\begin{enumerate}
  \item A ``smoothed'' $\tilde{R} \approx \omega R A M^{-1}_{\ell}$\replaced{
  where we take the transpose of the full operator above and}{.
  Implementations may} truncate the support, i.e.~throw away the
  \replaced{small negative entries}{smallish negative entries}
  by which the stencil support grows\replaced{. Furthermore, we}{, and}
  approximate $\tilde M_{\ell -1} = M_{\ell -1}$\added{, i.e.~reuse multigrid's
  correction operator within the damping term}.
  For this choice, memory requirements are slightly increased (we have to track
  one more ``unknown'') and two solves on all grid level besides the finest
  mesh are required (Algorithm \ref{algorithm:adAFAC:Jac}).
  The flow of data between grids can be seen in \replaced[id=R3]{(Figure \ref{figure:ingredients:adaFAC-JAC-flow})}{(Figure \ref{figure:ingredients:adaFAC-PI-flow})}.
  \item Sole injection \replaced{where we}{which can be seen as a mass-lumped
  smoothed stencil but using mass-lumping from finite elements as a presmoother. In
  this case, we} collapse $\tilde M_{\ell -1}IA_\ell$ into the identity. The
  overall damping reduces to $-\omega PIM_\ell^{-1}$. 
  We evaluate the original additive solution update. While we
  perform this update, we identify updates within $c$-points, i.e.~for vertices
  spatially coinciding with the next coarser mesh, inject these, immediately
  prolongate them down again, and damp the overall solution with the result.
  The damping equation is $PI$ (Algorithm \ref{algorithm:adAFAC:PI}).
  A schematic representation is shown in  \replaced[id=R3]{(Figure \ref{figure:ingredients:adaFAC-PI-flow})}{(Figure \ref{figure:ingredients:adaFAC-JAC-flow})}.
\end{enumerate}

\noindent
Both choices are  motivated through empirical observations. 
Our results study them for jumping coefficients in complicated domains, while 
our previous work demonstrates the suitability for 
Helmholtz-type setups \cite{Reps:17:Helmholtz}.
Though the outcome of both approaches is promising for our tests, we hypothesise
that more complicated setups such as convection-dominated phenomena require more
care in the choice of $\tilde R$, as $R$ has to be chosen more carefully \cite{Yavneh:12:NonnsymBoxMg}.

Both approaches can be combined with multigrid with geometric transfer operators
where $P$ is $d$-linear
everywhere or with algebraic approaches where $P$ stems from BoxMG. 
Both approaches inherit Ritz-Galerkin operators if they are used in the
baseline additive scheme.
Otherwise, they exploit redisretisation.

\begin{figure}
 \begin{center}
  \includegraphics[width=0.42\textwidth]{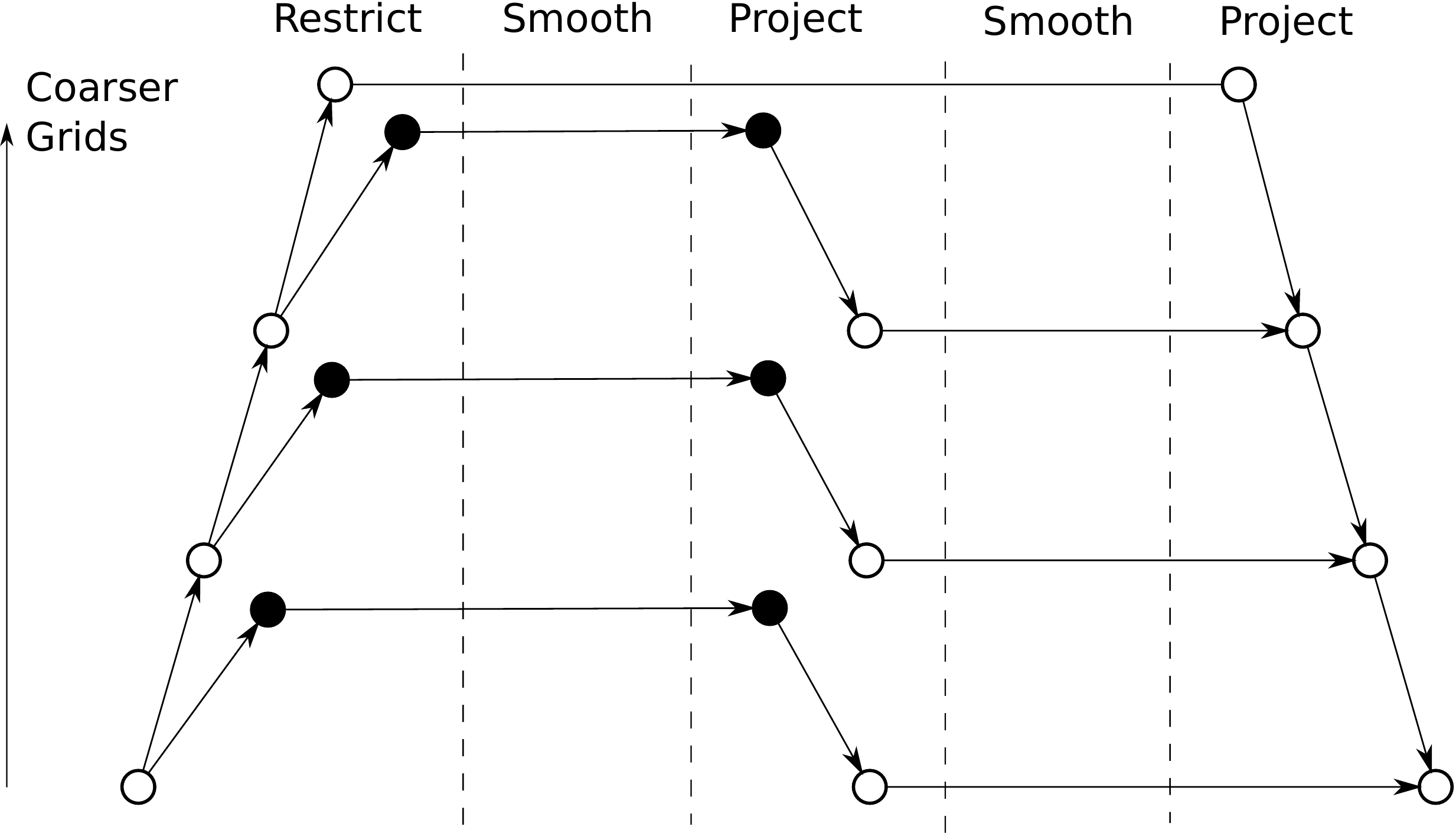}
  \hspace{1.4cm}
  \includegraphics[width=0.30\textwidth]{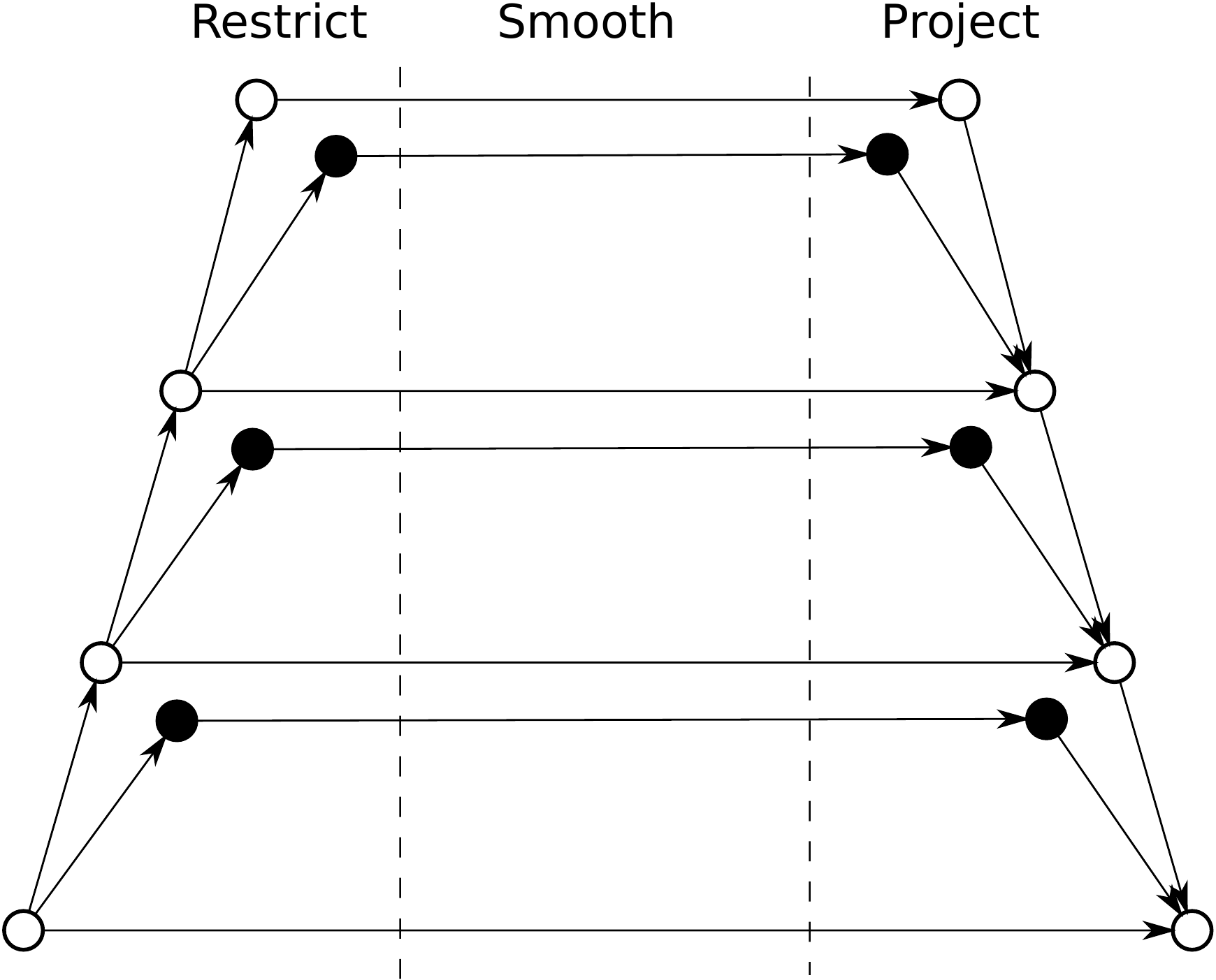}
 \end{center}
 \caption{
  Schematic overview of 
  AFACx (left) and our adAFAC (right). Black markers denote
  smoothing steps on the auxiliary equations, white markers correspond to
  traditional additive multigrid.
  \label{figure:ingredients:adaFAC-scheme}
 }
\end{figure}

\subsection{The AFAC solver family and further related approaches}

It is not a new idea to damp \replaced{the additive formulation of fast adaptive composite grids (FAC)}{FAC's additive formulation} such that additive
multigrid remains stable. 
Among the earliest endeavours is
FAC's asynchronous variant
\replaced{referred to as asynchronous fast adaptive composite grids (AFAC)}{AFAC} 
\cite{Lee:04:AFAC,McCormick:89:AFAC} which decouples the individual grid levels
to yield higher concurrency, too.
To remove oscillations, AFAC is traditionally served in two variants \cite{Hart:89:FAC}:

AFACc simultaneously determines the right-hand side for 
all grid levels $\ell $.
Before it restricts the fine grid residual to a particular level $\ell$, any
residuals on vertices spatially coinciding with vertices on the level $\ell $
are instead set to zero.
They are masked out on the fine grid.
This effectively damps the correction equation's right-hand side.
If we applied this residual masking recursively---a discussion explicitly not
found in the original AFACc paper where only the points are masked
which coincide with the target grid---i.e.~if we constructed
the masking restriction recursively over the levels instead of in one
rush, then AFACc would become a hybrid solver between additive multigrid and the
hierarchical basis approach \cite{Hart:89:FAC}.
 
AFACf goes down a different route:
The individual levels are treated independently from each other, but each
level's right-hand side is damped by an additional coarse grid contribution.
This coarse grid contribution is an approximate solve of the
correction term for the particular grid.
AFACf solves all meshes in parallel and sums up their contributions, but each mesh
has its contribution reduced by the local additional coarse grid cycle. 
The resulting scheme is similar to the combination technique as introduced
for sparse grids \cite{Bungartz:04:SparseGrids}:
We determine all solution updates additively but remove the intersection of
their coarser meshes.

Since multiplicative methods are superior to additive 
in terms of stability and simplicity, the transition from AFAC into AFACx
\cite{Lee:04:AFAC} is seminal:
Its inventors retain one auxiliary coarser level for each multigrid level, and 
split the additive scheme's solve into two phases
(Figure~\ref{figure:ingredients:adaFAC-scheme}):
A first phase determines per
level which modi might be eliminated by coarser grid solves.
For this, they employ the auxiliary helper level.
Each level keeps its additive right-hand side in the second phase, but it starts
with a projection from this auxiliary level as an initial guess.
The projection damps the correction after smoothing.
Only the resultant damped corrections derived in the second phase 
are eventually propagated to finer grids.
%
%
%
%

AFAC and FAC solvers
traditionally remain vague which solvers are to be used for particular substeps.
They are meta algorithms and describe a family of solvers. 
AFACx publications allow a free, independent choice of multigrid hierarchy
and auxiliary levels.
Our approach is different.
We stick to the spacetree construction paradigm.
As a result, real and auxiliary grid levels coincide.
Furthermore, we do not follow AFACx's multiplicative per-level update
(anticipate first the corrections made on coarser grids and then
determine own grid's contribution).
Instead, we run two computations in parallel (additively).
One is the classic additive correction computation.
The other term imitates the \added[id=R2]{effective} reduction of this update as compared
to multiplicative multigrid.
This additional, auxiliary term is subject to a single
smoothing step on one single auxiliary level which is the same as the next
additive resolution.

Our approach shares ideas with the Mult-additive approach
\cite{Yang:14:Reducing} where smoothed transfer operators are used in
the approximation of a $V(1,1)$ cycle.
Mult-additive yields faster convergence as it effectively yields stronger
smoothers.
We stick to the simple presmoothing approach and solely hijack the 
additional term to circumvent overshooting,
while the
asynchronicity of the individual levels is preserved.


We finally observe that our solver variant with a $PI$-term exhibits some
similarity with BPX.
BPX builds up its correction solely through inter-grid transfer operators while
the actual fine grid system matrix does not directly enter the
correction equations.
Though not delivering an explanation why the solver converges, the introduction
of the $PI$-scheme in \cite{Reps:17:Helmholtz} thus refers to this
solver as BPX-like.

\begin{idea}
 As our solver variants are close to AFAC, we call them adaptively damped AFAC
 and use the postfix $PI$ or $Jac$ to identify which damping equations we
 employ. Our manuscript thus introduces adAFAC-PI and adAFAC-Jac.
\end{idea}

\section{An element-wise, single-touch implementation}
\label{section:single-touch}

\begin{algorithm}[htb]
 \caption{
  Outline of single-touch adAFAC-Jac.
  A tilde identifies variables related to the auxiliary adAFAC grid.
  We invoke the cycle passing in the coarsest grid $\ell _{min}$\added[id=R2]{,
  i.e.~all helper variables are initially set to zero.}
    \label{algorithm:single-touch-adaFACx}
  }
  \begin{algorithmic}
    \Function{adaFAC-Jac}{$ \ell $}
    	\State $u_\ell \gets u_\ell + P^{\ell}_{\ell-1} sc_{\ell-1} - \widetilde{sl}$
    	  \Comment{Prolong contributions from both grids. \added[id=R2]{$sc$ holds
    	  the coarse grid corrections,}} 
    	\State $sc_\ell \gets sl_\ell - \widetilde{sl_\ell} + P^{\ell}_{\ell-1} sc_{\ell-1}- P ^{\ell}_{\ell-1}	\widetilde{sl}_{\ell-1}$
    	  \Comment{\replaced[id=R2]{i.e.~is recursively prolongated.}{Prepare
    	  for further prolongation}}
    	\State $u_\ell \gets u_\ell + sf_{\ell} - \widetilde{sf}_{\ell}$
    	  \Comment{Anticipate fine grid smoothing effects
    	  \added[id=R2]{deposited in $sf$.}} 
    	\State $\hat u_\ell \gets u_\ell -
    	  P^{\ell}_{\ell-1} I u_{\ell-1}$
    	  \Comment{Determine hierarchical residual.}  
    	\State $b_{\ell} \gets 0$;
    		$\widetilde{b}_{\ell} \gets 0$
    	  \Comment{Reset \replaced[id=R2]{right-hand side}{RHS} of correction
    	  equations.}

      \If {$\ell \neq {\ell _\text{max}}$}
    	\State \Call{adaFAC-Jac}{$l+1$}
      \EndIf

      \State $d _\ell = $ \Call{jacobi}{$A_\ell u_\ell = b_\ell, \omega $} 
        \Comment Determine update through Jacobi step.
      \State $\widetilde{d}_\ell = $
    	\Call{jacobi}{$\widetilde{A}_\ell
        0_\ell = \widetilde{b}_\ell , \omega $} 
        \Comment \replaced[id=R2]{Auxiliary}{Auxilliary} smooth with initial
        guess 0.
      \State $r_\ell \gets $ \Call{residual}{$A_\ell u_\ell = b_\ell$} 
        \Comment Bookmark residual from Jacobi update.
      \State $\hat{r}_\ell \gets $ \Call{residual}{$\widetilde{A}_\ell u_\ell = \widetilde{b}_\ell$} 
        \Comment Compute hierarchical residual.
      \State $sl _\ell \gets d _\ell $; 
             $u _\ell \gets u _\ell + d _\ell $; 
             $\widetilde{sl _\ell} \gets \widetilde{d _\ell} $
          \Comment{Bookmark \added[id=R2]{updates in $sl$} and apply
          \replaced[id=R2]{them, too.}{updates}}
        \State $b_\ell \gets R^{\ell-1}_{\ell}\hat r_{\ell +1}$;
        	$\widetilde{b}_\ell \gets \widetilde{R}^{\ell-1}_{\ell} r_{\ell +1}$
    	  \Comment Restrict \replaced[id=R2]{right-hand sides}{RHS} to coarse
    	  equation systems.
        \State $sf_{\ell -1} \gets I\left( sf_\ell + sl_\ell
          \right)$
        \Comment{Inform \replaced[id=R2]{coarser}{all} levels about
        updates\added[id=R2]{, but do not}} 
        \State $\widetilde{sf}_{\ell -1}
        \gets I\left( \widetilde{sf}_\ell + \widetilde{sl_\ell} \right)$
        \Comment{\replaced[id=R2]{apply them there. Information propagates
        bottom-up.}{do not apply them.}}

    \EndFunction

\end{algorithmic}
\end{algorithm}\vspace{-2.5mm}

adAFAC \replaced{works seamlessly with}{fits seamlessly to} our algorithmic building
blocks.
It solves up to three equations of the same type  per level.
\replaced{We}{Let's} distinguish the unknowns of these equations as follows:
$u$ is the solution in a \replaced{full approximation storage (FAS)}{FAS
multiscale} sense, while $\hat u$ is the hierarchical solution\deleted{ required
for HTMG}.
adAFAC solves
a correction equation, but no solution equation in the FAS sense. 
We do not to store need an additional $\widetilde u$ adAFAC unknown.
A complicated multi-scale representation along resolution boundaries is 
thus not required for the auxiliary damping equation: 
No semantic distinction between solution and correction
areas is required.
Let $d$ and $\widetilde d$ encode the iterative updates
of the unknowns through the additive \replaced{full approximation storage}{FAS
scheme} or the auxiliary adAFAC equation, respectively.

\subsection{Operator storage}
To make adAFAC stable and efficient for non-trivial $\epsilon $, each vertex
stores its element-wise operator parts from $A$.
Vertices hold the stencils.
For vertex members of the finest grid, the stiffness matrix entries result from the discretisation of
(\ref{equation:PDE}).
If we use $d$-linear inter-grid transfer operators this
storage scheme is applied to all levels.
Otherwise, we augment each vertex by further stencils for $P$ and $R$ and
proceed as follows:
For a vertex on a particular level which overlaps with finer resolutions, this
vertex belongs to a correction equation.
Its stencil results from the Ritz-Galerkin coarse grid
operator definition, whereas the inter-grid transfer operators $P$ and $R$
result from Dendy's BoxMG \cite{Dendy:82:BlackboxMG}.
BoxMG is well-suited for three-partitioning
\cite{Dendy:10:BoxMgBy3,Yavneh:12:NonnsymBoxMg,Weinzierl:18:BoxMG}.
We refer to \cite{Weinzierl:18:BoxMG} for remarks how to make the scheme
effectively matrix-free, i.e.~memory saving, nevertheless.
Each coarse grid vertex carries its prolongation and restriction
operator plus its stencil.
We are also required to store the auxiliary
$\widetilde R$ for adAFAC-Jac.
All further adAFAC terms use operators already held.
%

\subsection{Grid traversal}

For the realisation of the (dynamically) adaptive scheme, we follow
\cite{Mehl:06:MG,Reps:17:Helmholtz,Weinzierl:18:BoxMG,Weinzierl:19:Peano} and propose to run
through the spacetree in a depth-first (DFS) manner while each level's cells are
organised along a space-filling curve \deleted{(SFC) }\cite{Weinzierl:11:Peano}.
We write the code as recursive function where each cell has access to its $2^d$
adjacent vertices, its parent cell, and the parent cell's $2^d$ adjacent
vertices.
The latter ingredients are implicitly stored on the call stack of the recursive
function.

As we combine DFS with space-filling curves, our tree traversal is
weakly single-touch w.r.t.~the vertices: 
Vertices are loaded when an adjacent cell from the spacetree is first entered.
They are ``touched'' for the last time \replaced[id=R2]{once all $2^d$}{when the $2^d$th}
adjacent cell\added[id=R2]{s} within the
spacetree \replaced[id=R2]{have been}{is} left due to recursion backtracking.
In-between, they reside either on the call stack or can be temporarily stored in
stacks \cite{Weinzierl:11:Peano}.
The call stack is bounded by the depth of the spacetree---it is small---while
all temporary stacks are bounded by the time in-between the traversal of two
face-connected cells.
The latter is short due to the H\"older continuity of the underlying
\replaced{space-filling curve}{SFC}.
Hanging vertices per grid level, i.e.~vertices surrounded by less than $2^d$
cells, are created on-demand on-the-fly.
They are not held persistently.
We may assume that all data remains in the caches
\cite{Mehl:06:MG,Weinzierl:11:Peano,Weinzierl:19:Peano}.

As we extract element-wise operators for $A,P,R$ from the stencils stored within
the vertices or hard-code these element-wise operators,
we end up with a strict element-wise traversal in a multiscale sense.
All matrix-vector products (mat-vecs) are accumulated.
The realisation of the element-wise mat-vecs reads as
follows:
Once we have loaded the vertices adjacent to a cell,
we can derive the element-wise stiffness matrix or inter-grid transfer operator
for the cell.
To evaluate $r=Au$, we set one variable $r$ per vertex to zero, and
then accumulate the matrix-vector (mat-vec) contributions in each of the
vertex's adjacent cells.
Since the hierarchical $\hat u$ can be determined on-the-fly while running DFS
from coarse grids into fine grids, the evaluation of $\hat
r$ follows exactly $r$'s pattern.
So does the realisation of $\widetilde{r}$.
adAFAC's mat-vecs can be realised within a single spacetree traversal.
The mat-vecs are single-touch.

\subsection{Logical iterate shifts and pipelining}
%
%
A \replaced{full approximation storage}{(FAS)} sweep \deleted{however} can not
straightforwardly be realised within a single DFS grid sweep
\cite{Reps:17:Helmholtz}:
The residual computation propagates information bottom-up,
the corrections propagate information top-down, and the final
\deleted{FAS} injection propagates information bottom-up again.
This yields a cycle of causal dependencies.
We thus offset the additive cycle's smoothing steps by half a grid sweep:
Each grid sweep, i.e.~DFS traversal, evaluates all three mat-vecs---of FAS, of
\replaced{the hierarchical transformation multigrid}{HTMG}, of
adAFAC---but does not perform the actual updates.
Instead, correction quantities $sl$, $sc$, $sf$, $\widetilde{sc}$, $\widetilde{sf}$, and $\widetilde{sl}$ are bookmarked
as additional attributes within the vertices while the grid traversal
backtracks, i.e.~returns from the fine grids to the coarser ones.
Their impact is added to the solution throughout the downstepping of the
subsequent tree sweep.
Here, we also evaluate the prolongation.
Restriction of the residual to the auxiliary right-hand side and hierarchical residual 
continue to be the last
action on the vertices at the end of the sweep when variables are last
written/accessed.
As we plug into the recursive function's backtracking, we know that all
right-hand sides are accumulated from finer grid levels when we touch a vertex
for the last time throughout a multiscale grid traversal.
We can thus compute the unknown updates though we do not directly apply them.

As we use helper variables to store intermediate results
throughout the solve and apply them the next time, we need one tree
traversal per V-cycle plus one kick-off traversal.
Our helper variables pick up ideas behind pipelining and are a direct
translation of optimisation techniques proposed in \cite{Reps:17:Helmholtz}
to our scheme.
Per traversal, each unknown is read into memory/caches only once.
We obtain a single-touch implementation.
adAFAC's auxiliary equations do not harm its \replaced[id=R2]{suitability}{suitedness} to architectures with
a widening memory-compute facilities gap.

Dynamic mesh refinement integrates seamlessly into the single-touch traversal:
We rely on a top-down tree traversal which adds additional tree levels on-demand
throughout the steps down in the grid hierarchy.
The top-down traversal's backtracking drops parts of the tree if a coarsening
criterion demands so
it erases mesh parts.
Though the erasing feature is not required for the present test cases,
both refinement and coarsening integrate into Algorithm
\ref{algorithm:single-touch-adaFACx}.
We inherit FAC's straightforward handling of dynamic adaptivity,
simply the treatment of resolution transitions through \replaced{full
approximation storage}{FAS}, and provide an implementation which reads/writes
each unknown only once from the main memory.

%
%
\added{
 In a parallel environment, 
 the logical offset of the computations by half a grid sweep allows us to
 send out the partial (element-wise) residuals along the domain boundaries non-blockingly and to
 accumulate them on the receiver side.
 As the restriction is a linear operator, residuals along the domain
 boundaries can be restricted partially, too, before they are exchanged on the
 next level.
 Prior to the subsequent multiscale mesh traversal, all residuals and
 restricted right-hand sides are received and can be merged into the local
 residual contributions.
 If a rank is responsible for a cell, it also holds a copy of its parent cell.
 As we apply this definition recursively, a rank can restrict its partial
 residuals all the way through to the coarsest mesh.
 As we apply this definition recursively, all ranks hold the spacetree's root.
 The two ingredients imply that parallel adAFAC can be realised with one
 spacetree traversal per cycle.
 The two ingredients also imply that we do not run the individual grid levels in
 parallel even though we work with additive multigrid.
 Instead, we vertically integrate the levels, i.e.~if a tree traversal is
 responsible for a certain fine grid partition of the mesh, it is also
 responsible for its coarser representations, and the traversal of all of the
 grid resolution is collapsed (integrated) into one mesh run-through.
}

\section{Numerical results}
\label{section:results}

\begin{figure}
 \begin{center}
   \includegraphics[width=0.62\textwidth]{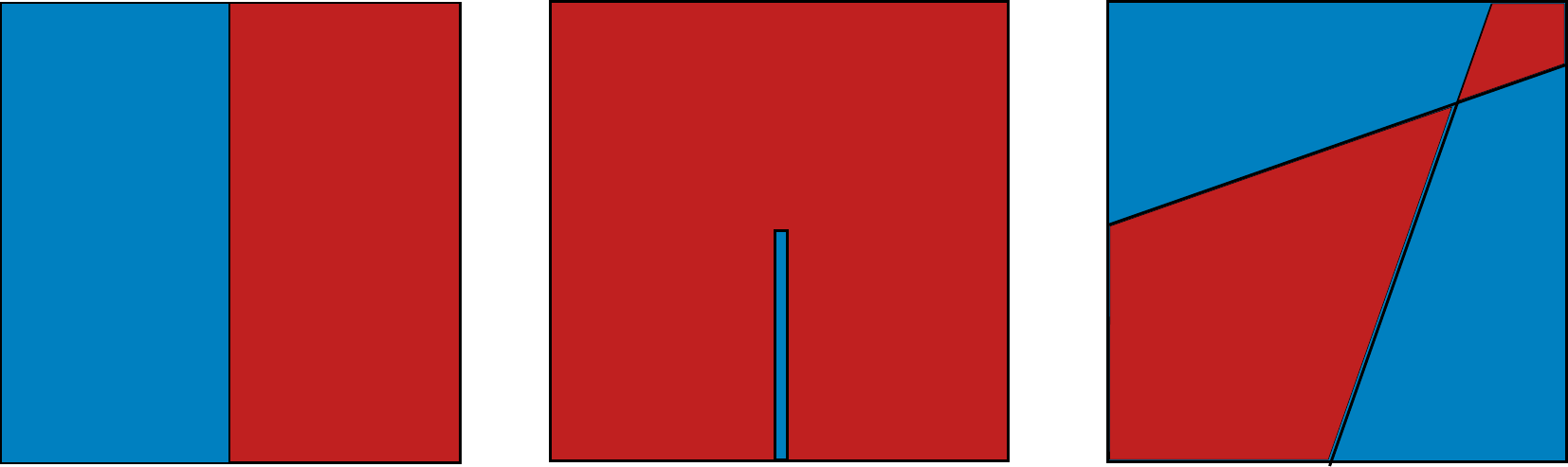}
 \caption{
   The three $\epsilon $ distributions studies throughout the tests.
   The blue area (left; inside of the inclusion;or top left and bottom right respectively) holds $\epsilon =1 $. The remaining domain $\epsilon =
   10^{-k}, \ k \in \{1,2,\ldots , 5\}$.
   \label{figure:results:grid-regions}
 }
 \end{center}
\end{figure}

To assess the potential of adAFAC, we study three test setups of type
(\ref{equation:PDE}) on the unit square.
They are simplistic yet already challenging for multigrid. 
All setups use $f=0$ and set $u|_{\partial \Omega } = 1$ for $x_2=0$ and
$u|_{\partial \Omega } = 0$ otherwise.
\added{This discontinuity in the boundary conditions and the reduced global smoothness
 are beyond the scope here.
 However, they stress the adaptivity criteria and pose a 
 challenge for the multigrid corrections. 
 The criteria have to refine towards the corners, while a sole geometric
 prolongation close to the corner introduces errors.
 Similar effects are well-known among lid-driven cavity setups for example.
}

A first test is the sole Poisson equation with a homogeneous
material parameter. 
The three other setups (Figure \ref{figure:results:grid-regions}) use regions
with $\epsilon =1$ and regions with $\epsilon = 10^{-k}$.
Per run, the respective $k \in \{1,2,\ldots , 5\}$ is fixed.
The second setup splits up the parameter domain into two equally sized
sections.
We emphasise that the split is axis-aligned but does not coincide with the mesh
as we employ three-partitioning of the unit square.
The third setup penetrates the area with a thin protruding line of width $0.02$. 
This line is parameterised with $\epsilon $.
It extends from the $x_0$
axis---$(x_1,x_2)^T \in \mathbb{R}^2$ are the coordinate axes---and
terminates half-way into the domain.
Such small material inhomogeneities cannot be represented
explicitly on coarse meshes.
The last setup makes the lines $x_2 = 5x_1 - 2.5$ and $x_2 = 0.2x_1 + 0.5$
separate domains which hold different $\epsilon $ in a checkerboard fashion.
No parameter split is axis-aligned.

If tests are labelled as regular grid runs, each grid level is regular and we
consequently end up with a mesh holding $(3^7-1)^d=4,778,596$ degrees of freedom
for $\ell_{max}=7$.
If not labelled as regular grid run, our tests rely on 
dynamic mesh refinement.

Otherwise our experiments focus on $d=2$ and start with a 2-grid algorithm
($\ell_{max}=2$) where the coarser level has $(3-1)^d=4$ degrees of freedom and
the finer level hosts $(3^2-1)^d=64$ vertices carrying degrees of freedom.
From hereon, we add further grid levels and
build up to an 8-grid scheme ($\ell_{max}=8$).

In every other cycle, our code manually refines the cells along the bottom
boundary, i.e.~the cells where one face carries
$u|_{\partial \Omega}\not = 0$.
We stop with this refinement when the grid meets $\ell _{max}$.
Our manual mesh construction ensures that we kick off with a low total vertex
count, while the solver does not suffer from pollution effects:
The scheme kickstarts further feature-based refinement.
Parallel to the manual refinement along the boundary, our implementation
measures the absolute second derivatives of the solution along both coordinate
axes in every single unknown.
A bin sorting algorithm is used to identify the vertices carrying the
(approximately) 10 percent biggest directional derivatives.
These are refined unless they already meet $\ell _{max}$. 
The overall approach is similar to full multigrid where coarse grid solutions
serve as initial guesses for subsequent cycles on finer meshes,
though our implementation lacks higher-order operators.
All interpolation from coarse to fine meshes both for hanging vertices and for
newly created vertices is $d$-linear.

Our runs employ a damped Jacobi smoother with damping $\omega = 0.6$ and
report the normalised residuals
\begin{equation}
 \frac{\|r^{(n)}\|_h}{\|r^{(0)}\|_h}
 \quad \mbox{where} \quad  
	\|r^{(n)}\|^2_h := \sum_{i} h_i^d (r_i^{(n)})^2,
 \label{equation:results:normalised-residual}
\end{equation}

\noindent
with $n$ being the cycle count.
$r_i^{(n)}$ is the residual in vertex $i$ and
$h_i$ is the local mesh spacing around vertex $i$.
Dynamic mesh refinement inserts additional vertices and
thus might increase the residual vectors between two subsequent iterations.
As a consequence, residuals under an Eulerian norm may temporarily grow due to
mesh expansion.
This effect is amplified by the lack of higher order interpolation for new
vertices. 
The normalised residual (\ref{equation:results:normalised-residual}) 
enables us to quantify how much the residual has decreased compared to the residual fed
into the very first cycle.

Where appropriate, we also display the normalised maximum residual
\[
 \frac{max _i \{|r_i^{(n)}|\}}{max _i\{|r_i^{(0)}|\}}.
\]
This metric identifies localised errors,
while (\ref{equation:results:normalised-residual}) weights errors with the mesh
size.

\subsection{Consistency study: the Poisson equation}

%
%
Our first set of experiments focuses on the Poisson equation, i.e.~$\epsilon =
1$ everywhere.
Multigrid is expected to yield a perfect solver for this setup:
Each cycle (multiscale grid sweep) has to reduce the residual by a constant
factor which is independent of the degrees of freedom, i.e.~number of vertices.
Ritz-Galerkin multigrid yields the same operators as rediscretisation, since
BoxMG gives bilinear inter-grid transfer operators.
The setup is a natural choice to validate the consistency and correctness
of the adaFACx ingredients.
All grids are regular.

\begin{figure}
 \begin{center}
   \includegraphics[width=0.49\textwidth]{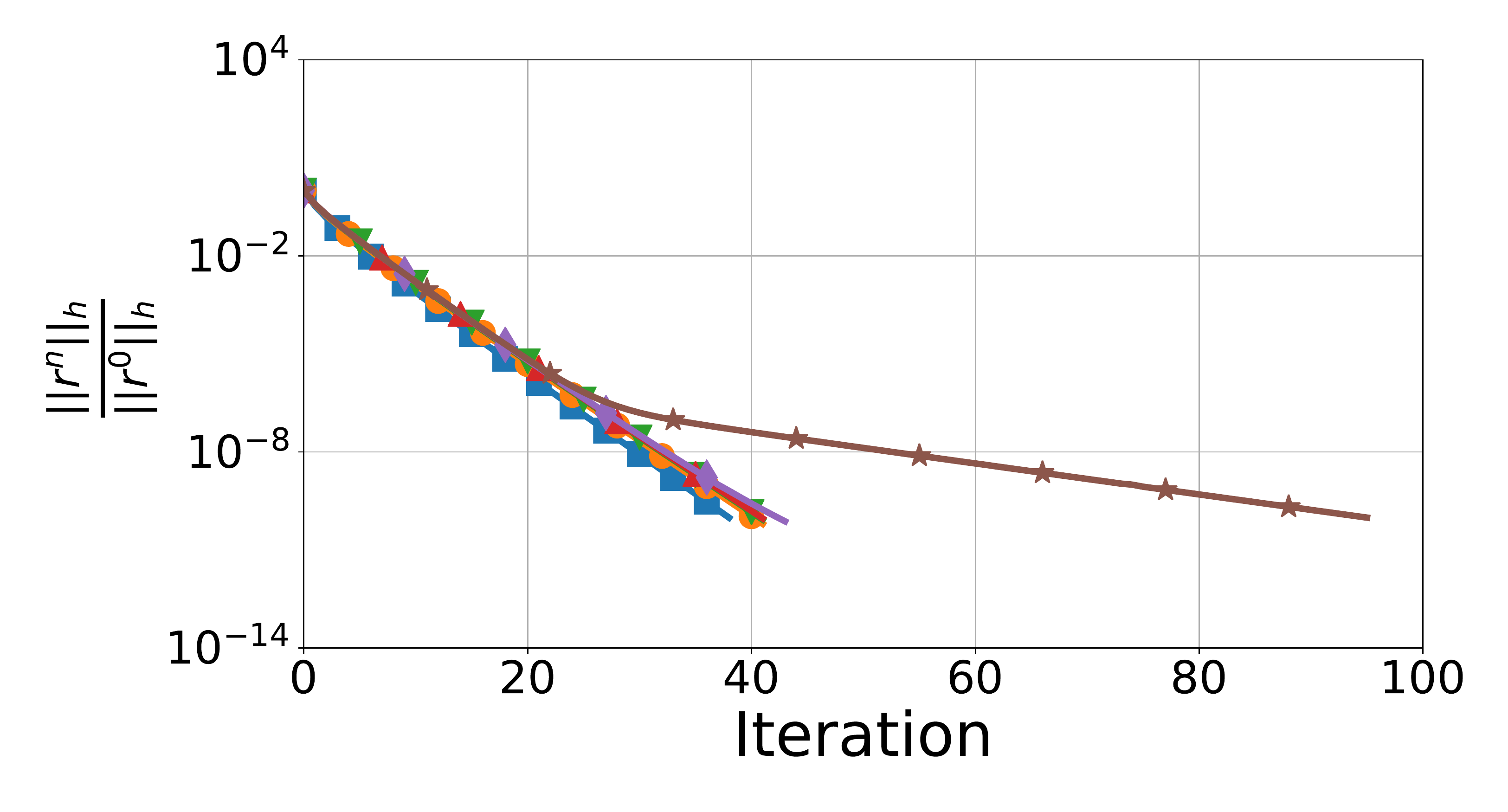}
   \includegraphics[width=0.49\textwidth]{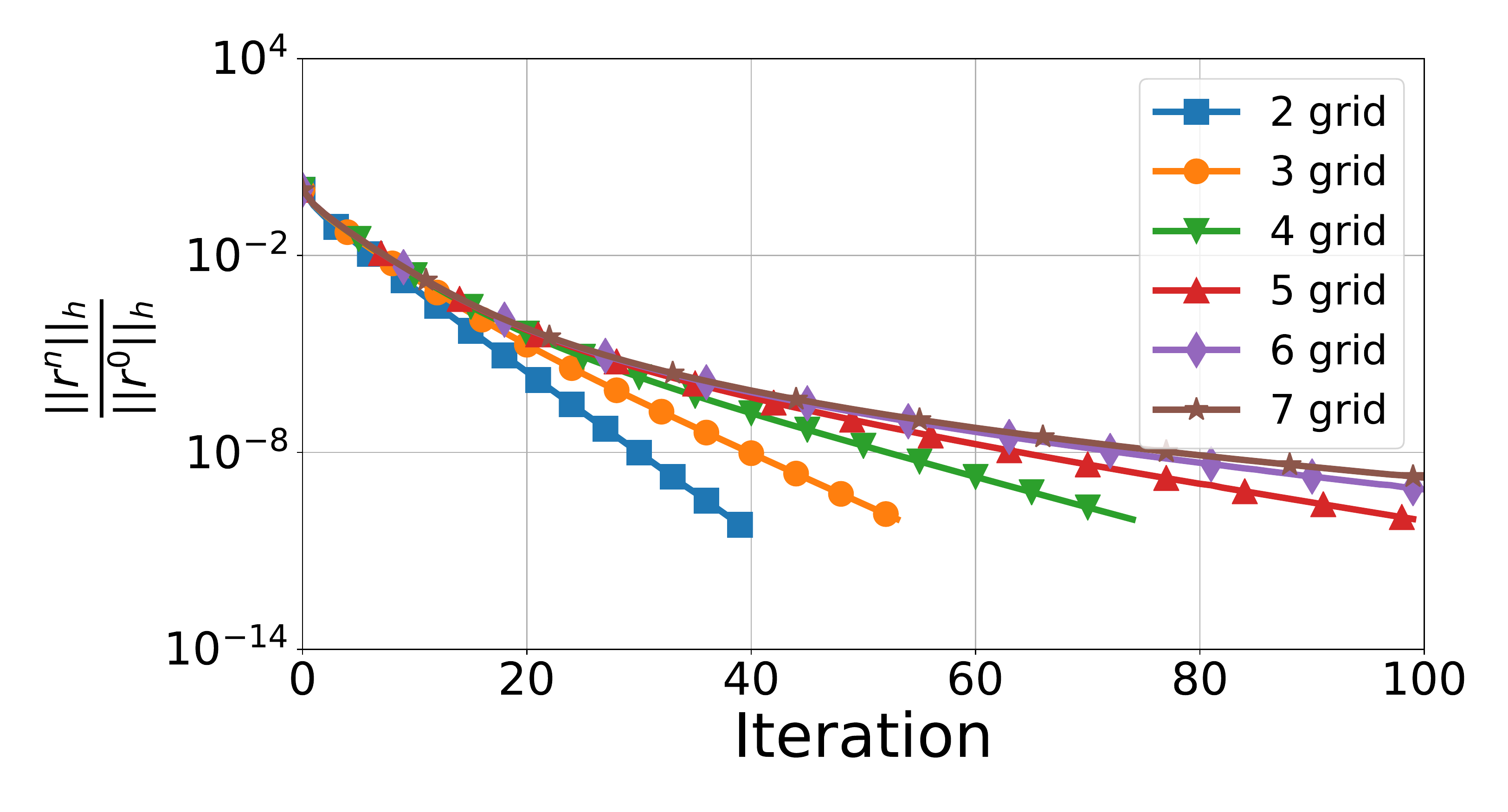}
   \\
   \includegraphics[width=0.49\textwidth]{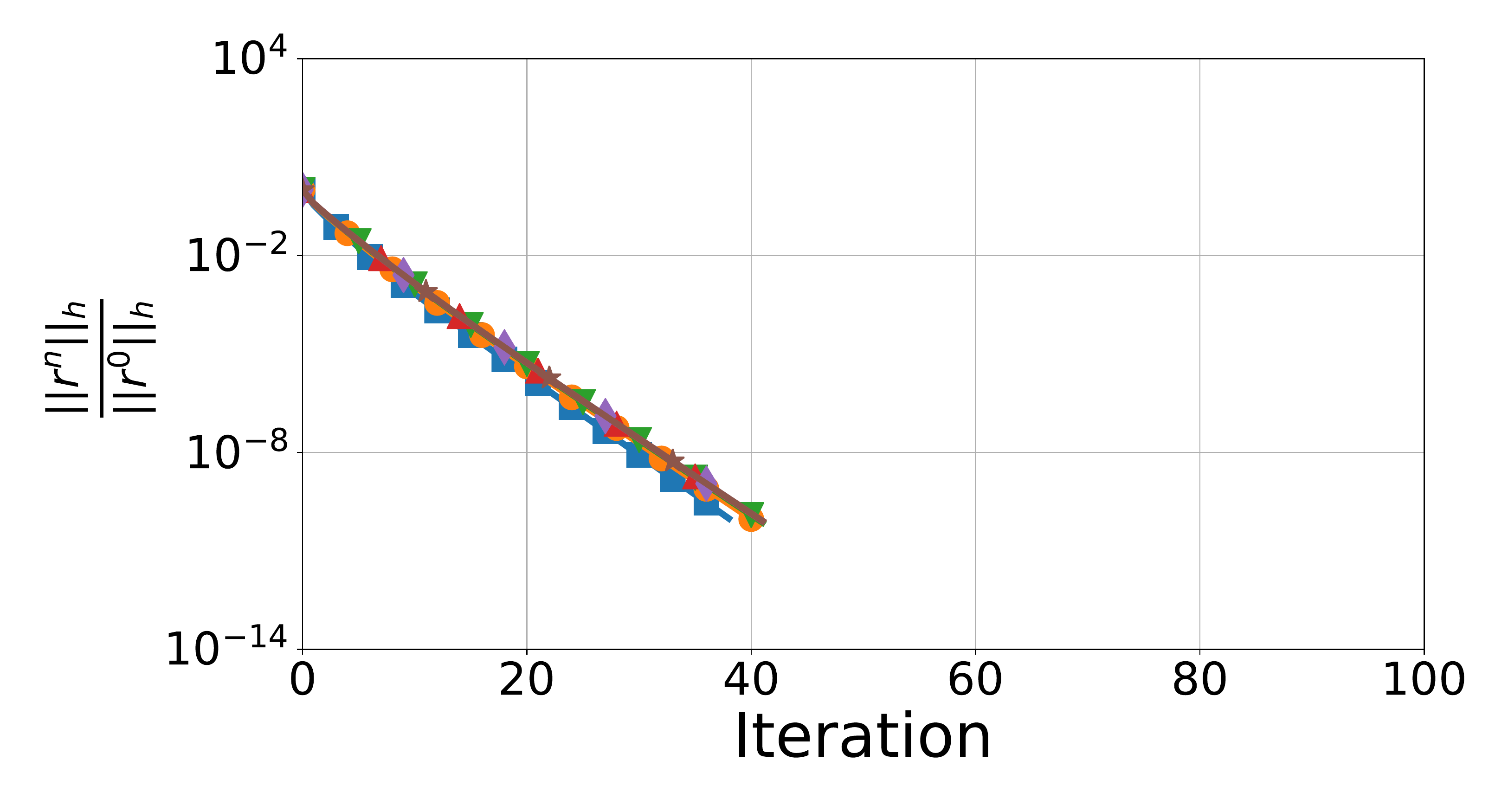}
   \includegraphics[width=0.49\textwidth]{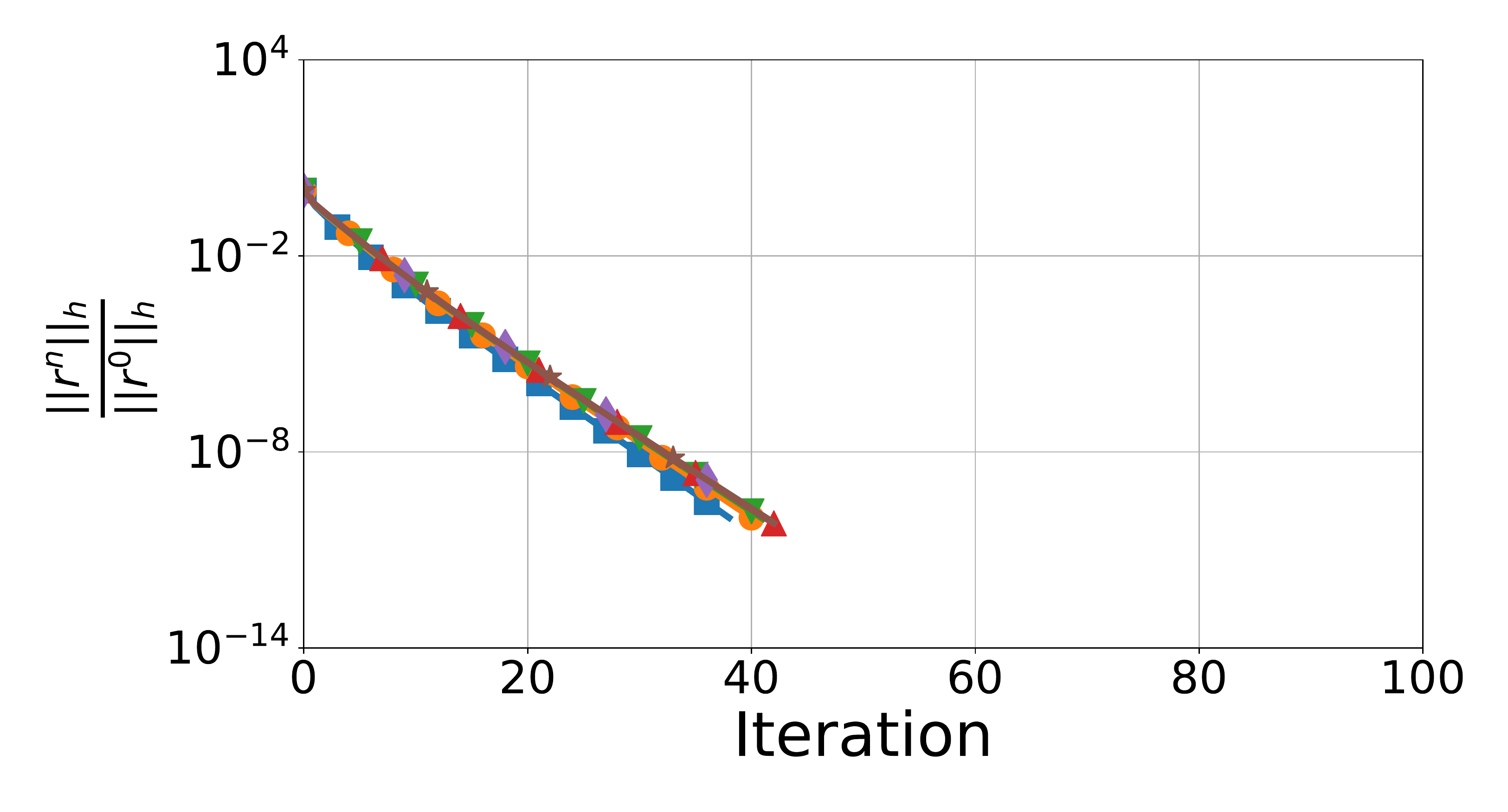}
 \end{center}
 \caption{
   Solves of the Poisson equation on regular grids of different levels.
   We compare plain additive multigrid (top, left), multigrid using
   exponential damping (top, right), adAFAC-PI (bottom, left) and adAFAC-Jac
   (bottom, right).
   \label{figure:results:Poisson:regular-grid:solvers}
 }
\end{figure}

%
%
Our experiments (Figure~\ref{figure:results:Poisson:regular-grid:solvers})
confirm that additive multigrid is insignificantly faster than the other
alternatives if it is stable.
The more grid levels are added, the more we overshoot per multilevel relaxation.
When we start to add a seventh level, this suddenly makes the plain additive
code's performance deteriorate.
With an eighth level added, the solver would diverge (not shown).
Exponentially damped multigrid does not suffer from the instability for
lots of levels, but its damping of coarse grid influence leads to the situation
that long-range solution updates don't propagate quickly.
The convergence speed suffers from additional degrees of freedom.
Both of our adAFAC variants are stable, but they do not suffer from a speed
deterioration.
Their \deleted{``clever'',} localised damping \replaced{makes both schemes
effective and stable}{pays off}.  
\replaced{We see similar rates of convergence in our damped implementations
to that of undamped additive multigrid.}{We play in the same league additive multigrid in terms of speed.}
adAFAC-PI and adAFAC-Jac are almost indistinguishable.

\begin{figure}
 \begin{center}
   \includegraphics[width=0.21\textwidth]{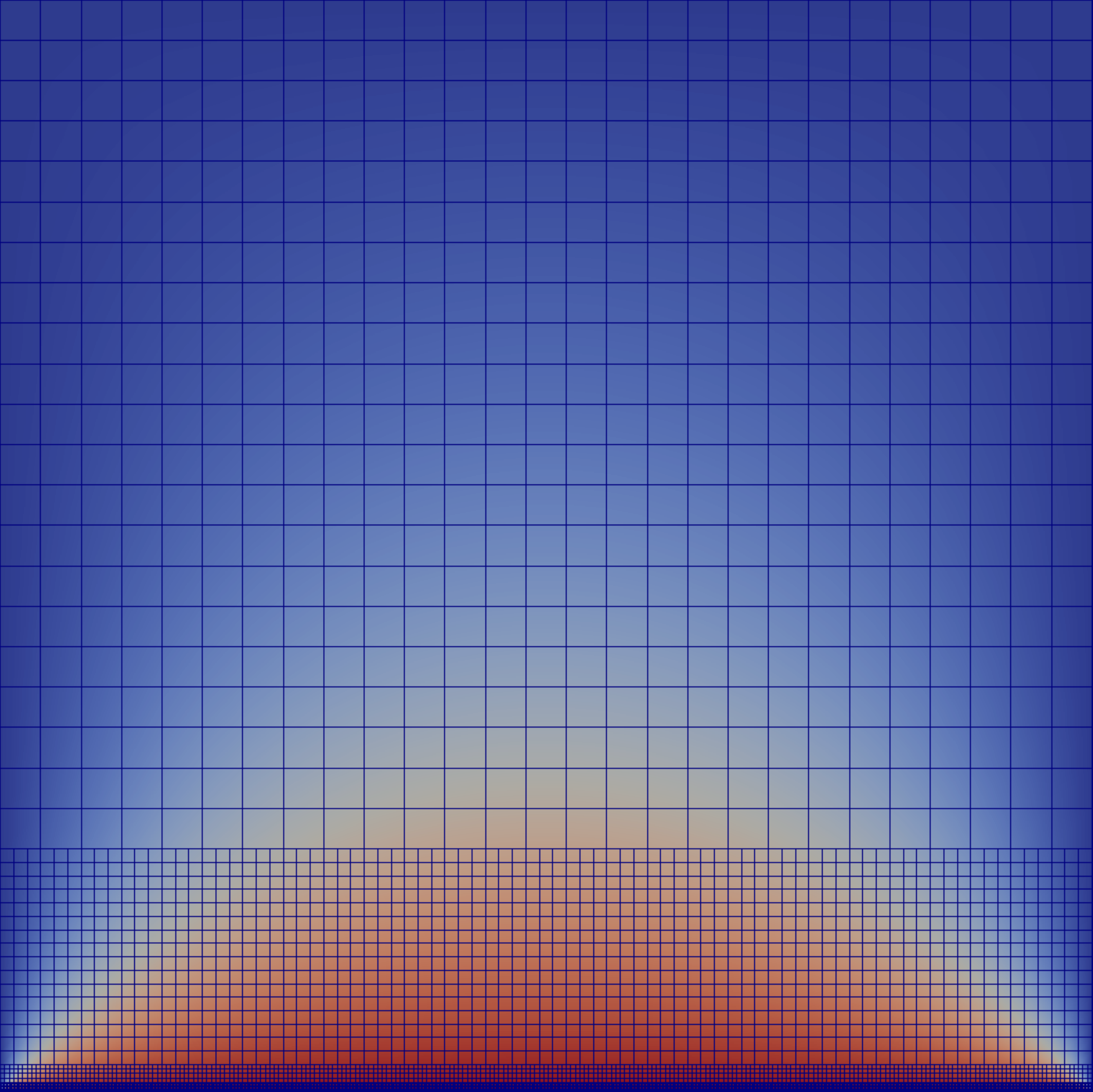}
   \hspace{1.4cm}
   \includegraphics[width=0.39\textwidth]{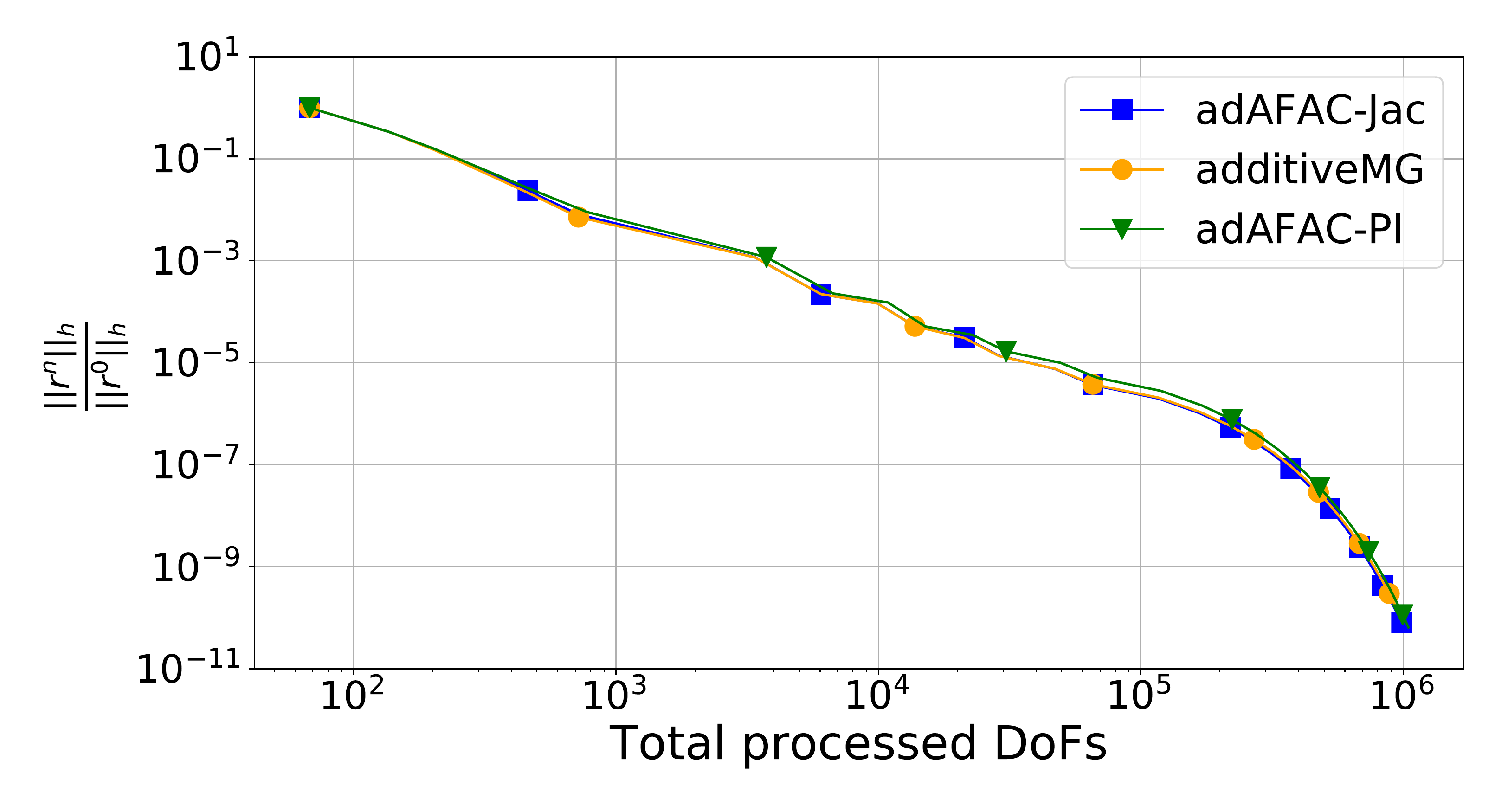}
 \end{center}
 \caption{
   Left: Typical adaptive mesh for pure Poisson (constant material parameter)
   once the refinement criterion has stopped adding further elements.
   Right: We compare different solvers on the pure Poisson equation using a
   hybrid FMG-AMR approach starting at a two grid scheme and stopping at an eight grid
   scheme. $\ell _{max}=8$.
   \label{figure:Poisson:adaptive-grid:solvers}
 }
\end{figure}

%
%
Despite the instability of plain additive multigrid, we continue to benchmark
against the undamped additive scheme, as exponential damping is not
competitive. 
All experiments from hereon are reasonable irregular/coarse to circumnavigate
the instabilities.
Feature-based dynamic refinement criterion makes the mesh spread out from
the bottom edge where
$u|_{\partial \Omega}=1$ (Figure \ref{figure:Poisson:adaptive-grid:solvers}).
To assess its impact on cost, we count the number of required degrees of freedom
updates plus the updates on coarser levels.
\added{
These degree of freedom updates correlate directly to runtime.}
We do not neglect the coarse grid costs.

One smoothing step on a regular mesh of level eight yields 
$4.3 \cdot 10^7$ updates plus the updates on the correction levels.
If the solver terminated in 40 cycles, we would have to
invest more than $10^9$ updates.
Dynamic mesh unfolding reduces the cost to reduce the residual by up to three
orders of magnitude.
For Poisson, this saving applies to both our adAFAC variants and plain
additive multigrid, while the latter remains stable.

%
%
If ran with BoxMG, our \replaced[id=R2]{codebase}{code base} uses Ritz-Galerkin coarse operator construction for both the
correction terms and the auxiliary adAFAC operators in adAFAC-Jac.
We validated that both the algebraic inter-grid transfer operators and 
geometric operators yield exactly the same outcome.
This is correct for Poisson as
BoxMG yields geometric operators here
and Ritz-Galerkin coarse operator construction for the correction terms thus
yields the same result as rediscretisation.

\begin{observation}
Our code is consistent.
For very simple, homogeneous setups, it however makes only limited sense to use
adAFAC---unless there are many grid levels.
If adAFAC is to be used, adAFAC-PI is sufficient.
There's no need to really solve an additional auxiliary equation.
\end{observation}

\added{
 \noindent
 We conclude with the observation that all of our solvers, if stable, converge
 to the same solution.
 They are real solvers, not mere preconditioners that only yield approximate solutions.
}

\subsection{One material jump}

\begin{figure}[htb]
 \begin{center}
   \includegraphics[width=0.21\textwidth]{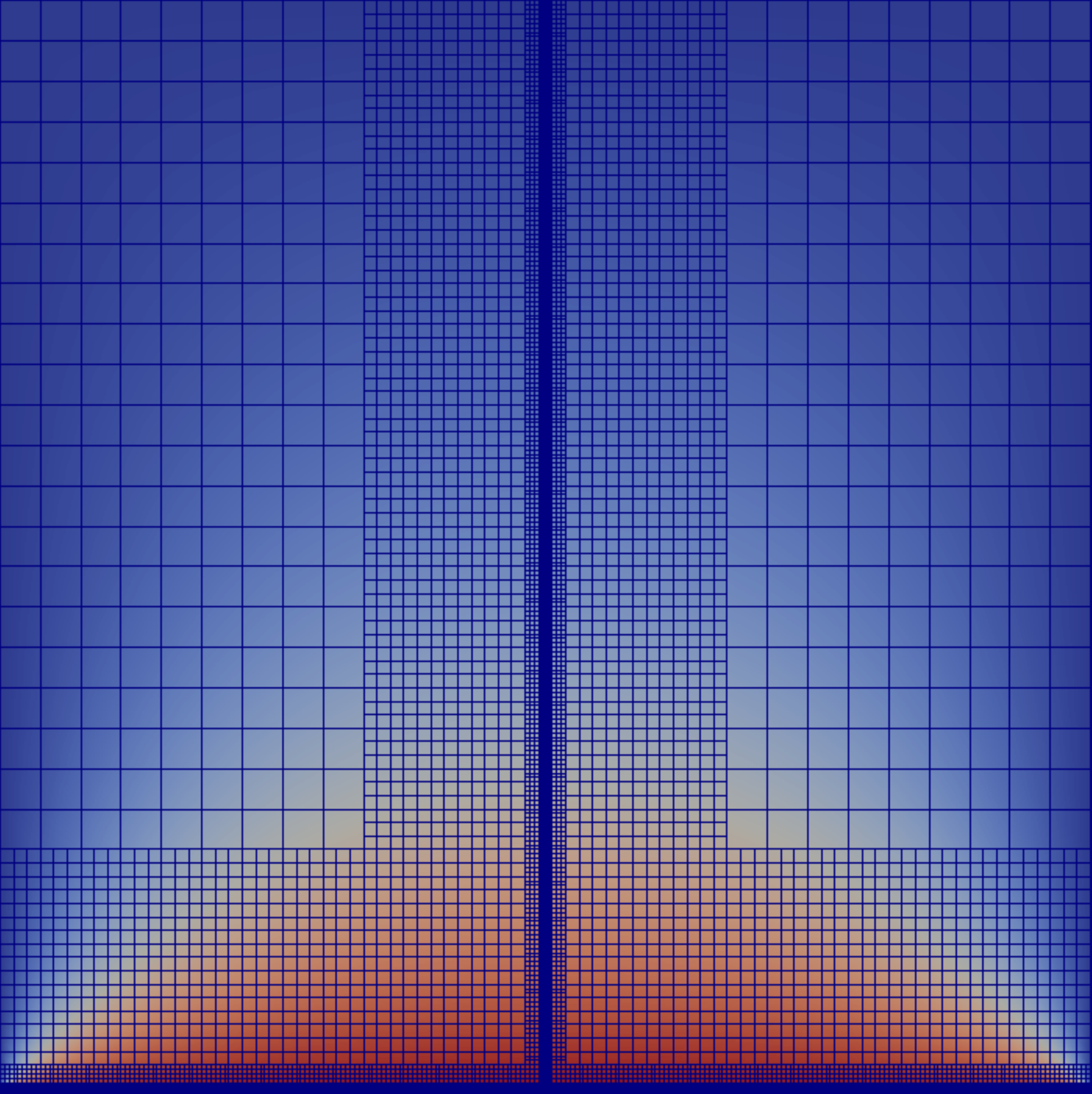}
   \includegraphics[width=0.39\textwidth]{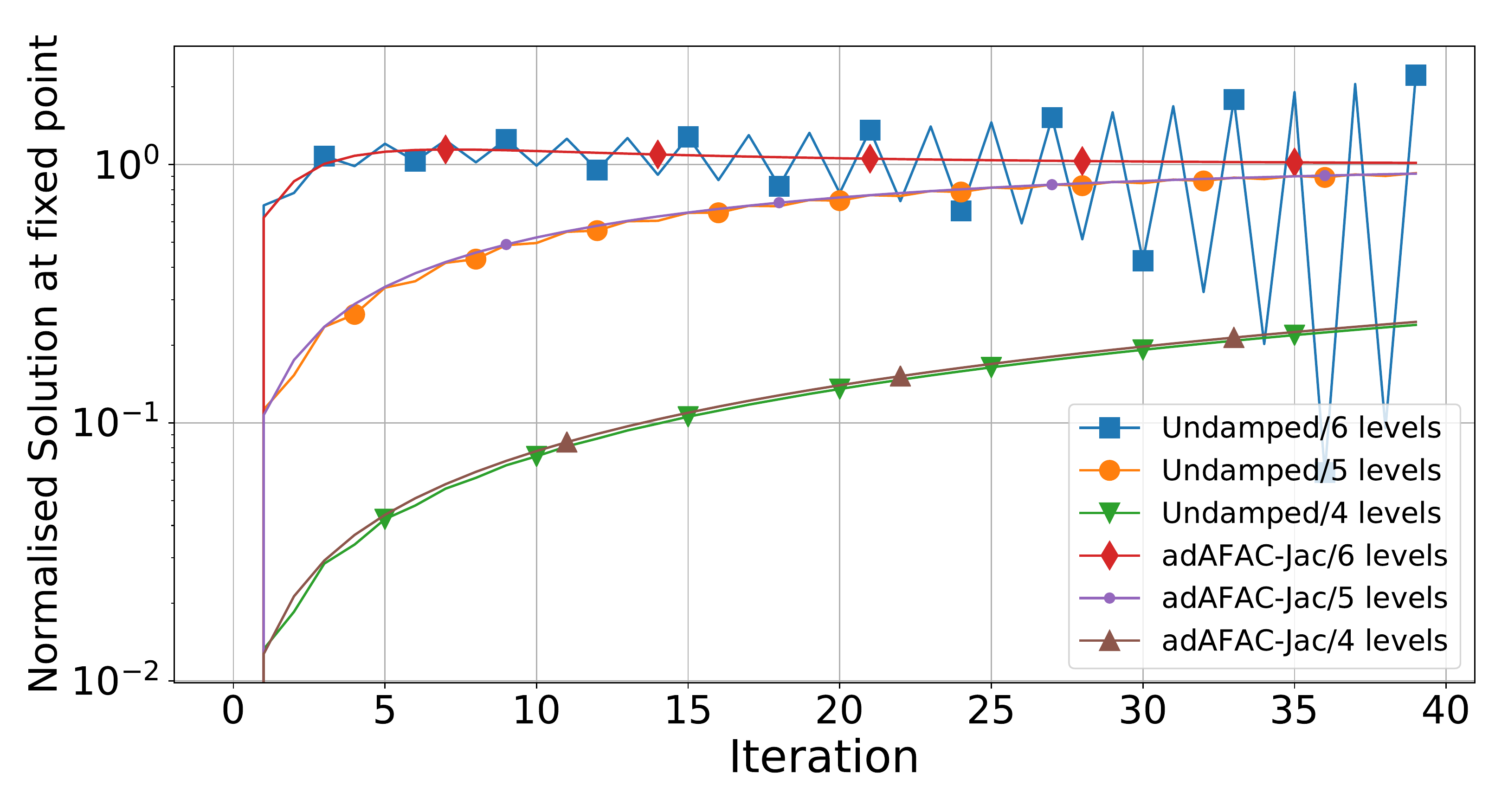}
   \includegraphics[width=0.39\textwidth]{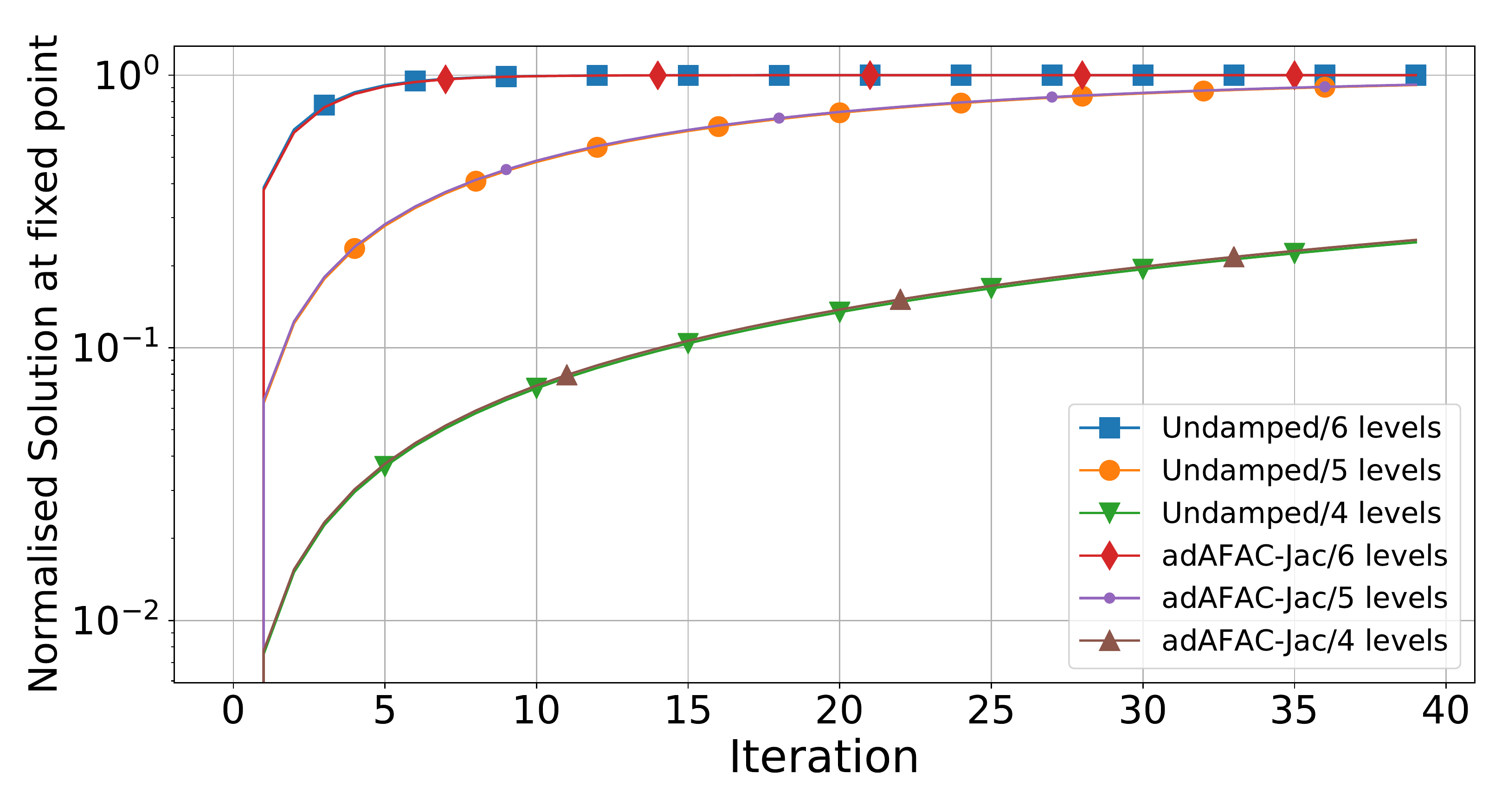}
 \caption{
  \added{
   The domain material is split into two halves with an $\epsilon $ jump
   from $\epsilon = 1$ to \replaced[id=R2]{$\epsilon = 10^{-7}$}{$\epsilon =
   0.1$}.
   Typical adaptive mesh for single discontinuity setup once the
   refinement criterion has stopped adding further elements (left). 
   Solution development in sample point next to a
   discontinuity\replaced[id=R2]{, normalised}{. Normalised} by
   \added[id=R2]{the} true solution value at that point, i.e.~one means the
   correct value.
   We compare $d$-linear inter-grid transfer (middle) to BoxMG operators
   (right).
  }
   \label{figure:results:solution-values}
 }
 \end{center}
\end{figure}

%
%
We next study a setup where the material ``jumps'' in the middle of the domain.
The stronger the material transition is the more important it is to pick up the
$\epsilon $ changes in the inter-grid transfer operators.
Otherwise, \added{a} prolongation of coarse grid corrections yields errors close
to $x_1=0.5$.
As no grid in the present setup has degrees of freedom exactly on the material
transition, the inter-grid transfer operators are never able to mirror the material transition exactly.
\deleted{It is the dynamic adaptivity which counterbalances this shortcoming:
Large errors in $P$ are compensated by many corrections on finer grids.}

\added{
Without dynamic adaptivity, multigrid runs the risk of deteriorating in the 
multiplicative case and becoming unstable in the additive case.
To document this phenomenon, we monitor the solution in one sample point
coinciding with the real degree of freedom next to $x=(0.5,0.5)^T$, and 
employ a jump in $\epsilon $ of seven orders of magnitude.
A regular grid corresponding to $\ell = 6$ is used.
We start from a single grid algorithm, and add an increasing number of correction levels.
Not all grid level setups are shown herein.
Without algebraic inter-grid transfer operators, oscillations 
arise if we do not use our additional damping parameter 
\replaced[id=R3]{(Figure~\ref{figure:results:solution-values})}{(Fig.~\ref{figure:results:solution-values})}.
The oscillations increase with the number of coarse grid
levels used.
Our damping parameter eliminates these oscillations and does not harm
the rate of convergence.
Algebraic \replaced[id=R3]{inter-grid}{intergrid} transfer operators eliminate
these oscillations, too.
The results show why codes without algebraic operators and without
damping usually require a reasonably coarse mesh to align with $\epsilon $
transitions.
}

\begin{figure}[htb]
 \begin{center}
   \includegraphics[width=0.39\textwidth]{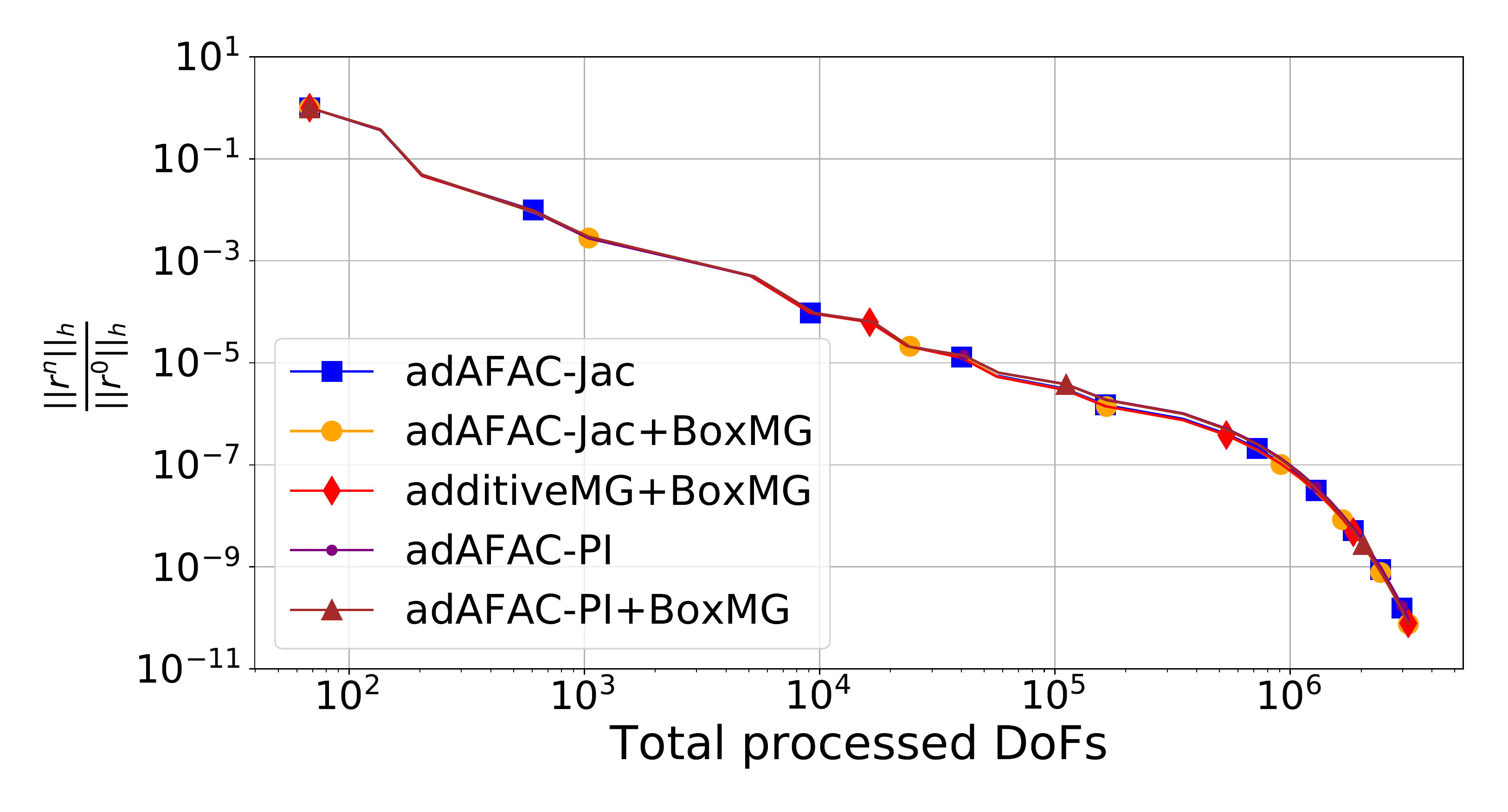}
   \includegraphics[width=0.39\textwidth]{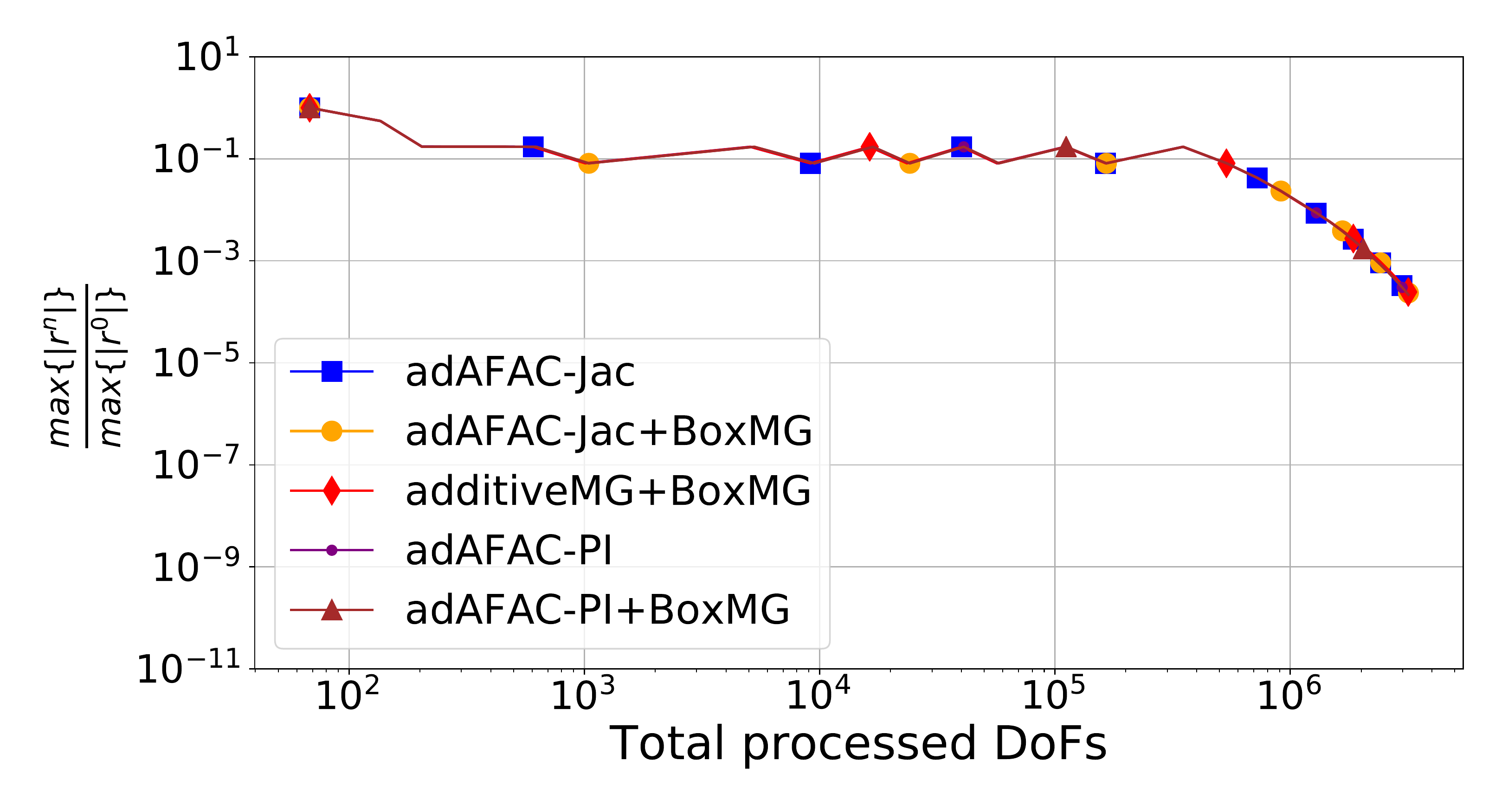}
 \end{center}
 \caption{
   The \replaced[id=R2]{left}{centre} plot shows the normalised residual and the
   right shows the normalised $L^{\infty}$-norm of the residual.
   $\ell _{max}=8$.
   \added[id=R2]{$\epsilon$ in $\{1,10^{-1}\}$, i.e.~the material parameter
   changes by one order of magnitude.
   We present only data for converging solver flavours.
   }
   \label{figure:Poisson:adaptive-grid:jump-solvers-1}
 }
\end{figure}

\begin{figure}[htb]
 \begin{center}
   \includegraphics[width=0.39\textwidth]{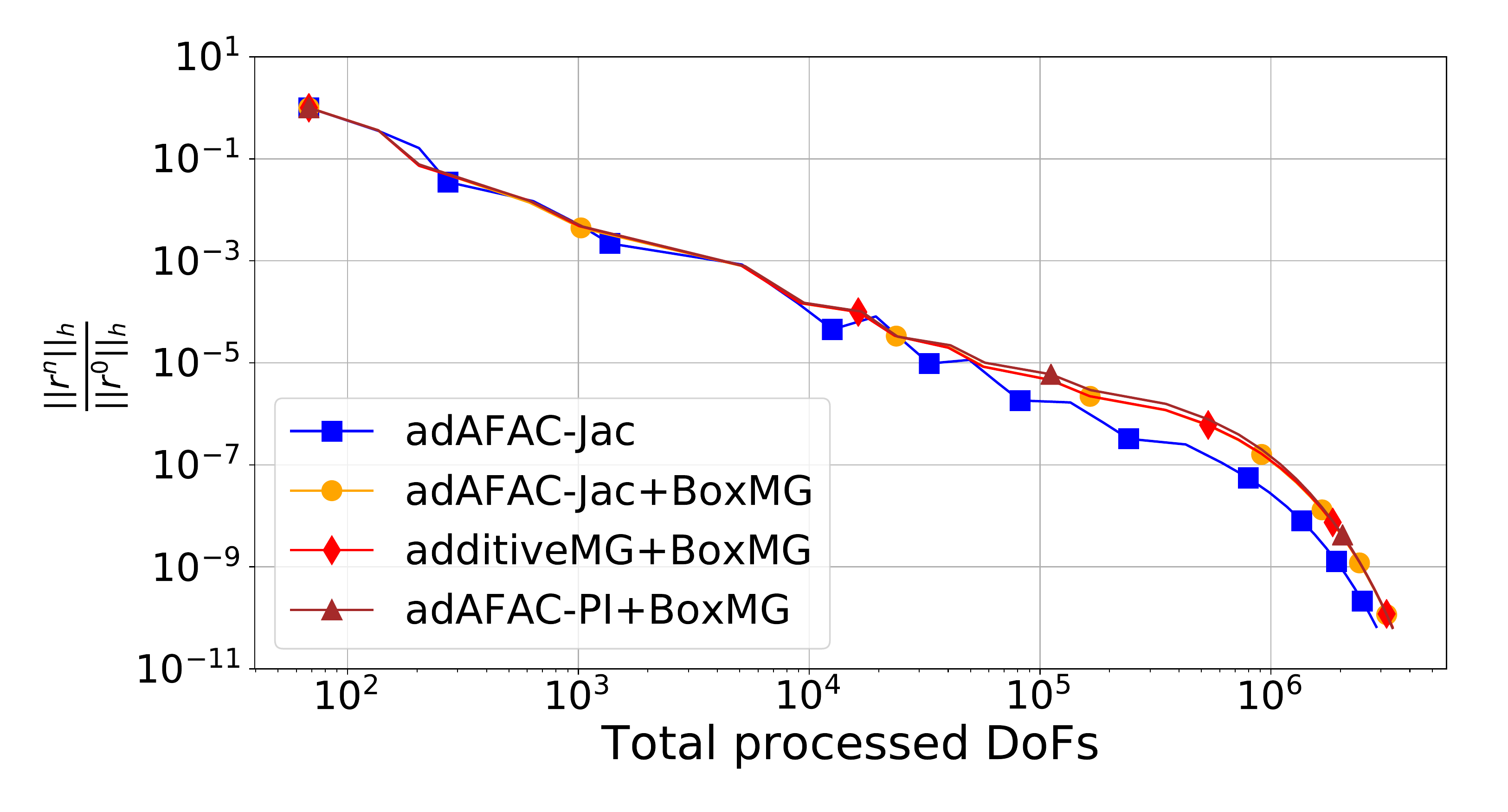}
   \includegraphics[width=0.39\textwidth]{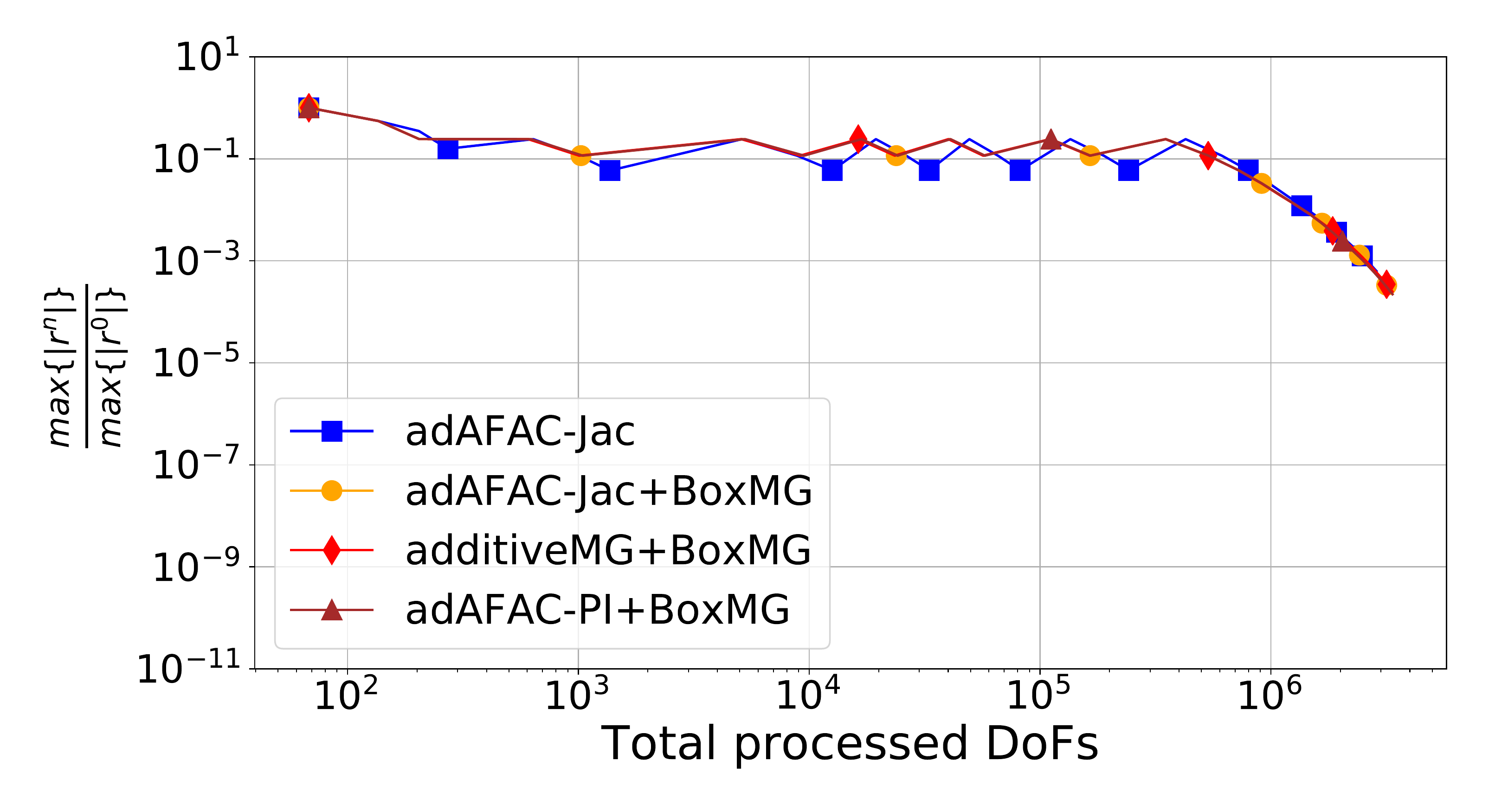}
 \end{center}
 \caption{
   Setup of Figure \ref{figure:Poisson:adaptive-grid:jump-solvers-1} but with a
   five orders of magnitude jump in the material parameter.
   We present only data for converging runs\added[id=R2]{ and observe that
   fewer solver ingredient combinations converge}.
   \label{figure:Poisson:adaptive-grid:jump-solvers-2}
 }
\end{figure}

\added{
 We continue with dynamically adaptive meshes.
}
All experiments use the AMR/FMG setup, i.e.~start from a
coarse mesh and then dynamically adapt the grid.
We observe that the hard-coded grid refinement refines along the stimulus
boundary at the bottom, while the dynamic refinement criterion unfolds it along
the material transition
(Figure~\ref{figure:results:solution-values}).

%
%
Starting from reasonably small changes in $\epsilon $ (Figure
\ref{figure:Poisson:adaptive-grid:jump-solvers-1}), additive multigrid
\added{with geometric inter-grid transfer operators again} fails to converge.\deleted{without the addition of BoxMG.}
Once BoxMG is used, it \deleted{however }becomes stable.
\replaced{The}{This is reasonable, as the} residual plot in the maximum norm
validates our statement that large errors arise along the material transition when we insert
new degrees of freedom.
We need an algebraic interpolation routine. 
Our adAFAC variants in contrast all converge.
The absence of a higher-order interpolation for new degrees of
freedom hurts, but it does not destroy the overall stability.
Once the dynamic AMR stops inserting new vertices---this happens after around
$10^6$ degrees of freedom have been processed---the residual
drops under both norms.

The picture changes when we increase the variation in $\epsilon $.
adAFAC-Jac with bilinear transfer operators converges for all
$\epsilon = 10^{-k}$ values tested,
whereas additive multigrid and adAFAC-PI diverge without BoxMG
\added[id=R2]{if the $\epsilon$-transition is too harsh}
(Figure \ref{figure:Poisson:adaptive-grid:jump-solvers-2}).
The geometric inter-grid transfer approach suffers 
from oscillations around the material transition.
All stable solvers play in the same league.


\begin{observation}
If we face reasonably small jumping materials, adAFAC-PI is superior to plain
additive multigrid, adAFAC-Jac or any algebraic-geometric extension, as it is
both stable and simple to compute.
Once the jump grows, adAFAC-Jac becomes the method of \replaced[id=R2]{choice}{joice}.
Its auxiliary damping equations compensates for the lack of algebraic inter-grid
transfer operators which are typically not cheap to compute.
\end{observation}

\subsection{A material inclusion}

\begin{figure}
 \begin{center}
   \includegraphics[width=0.21\textwidth]{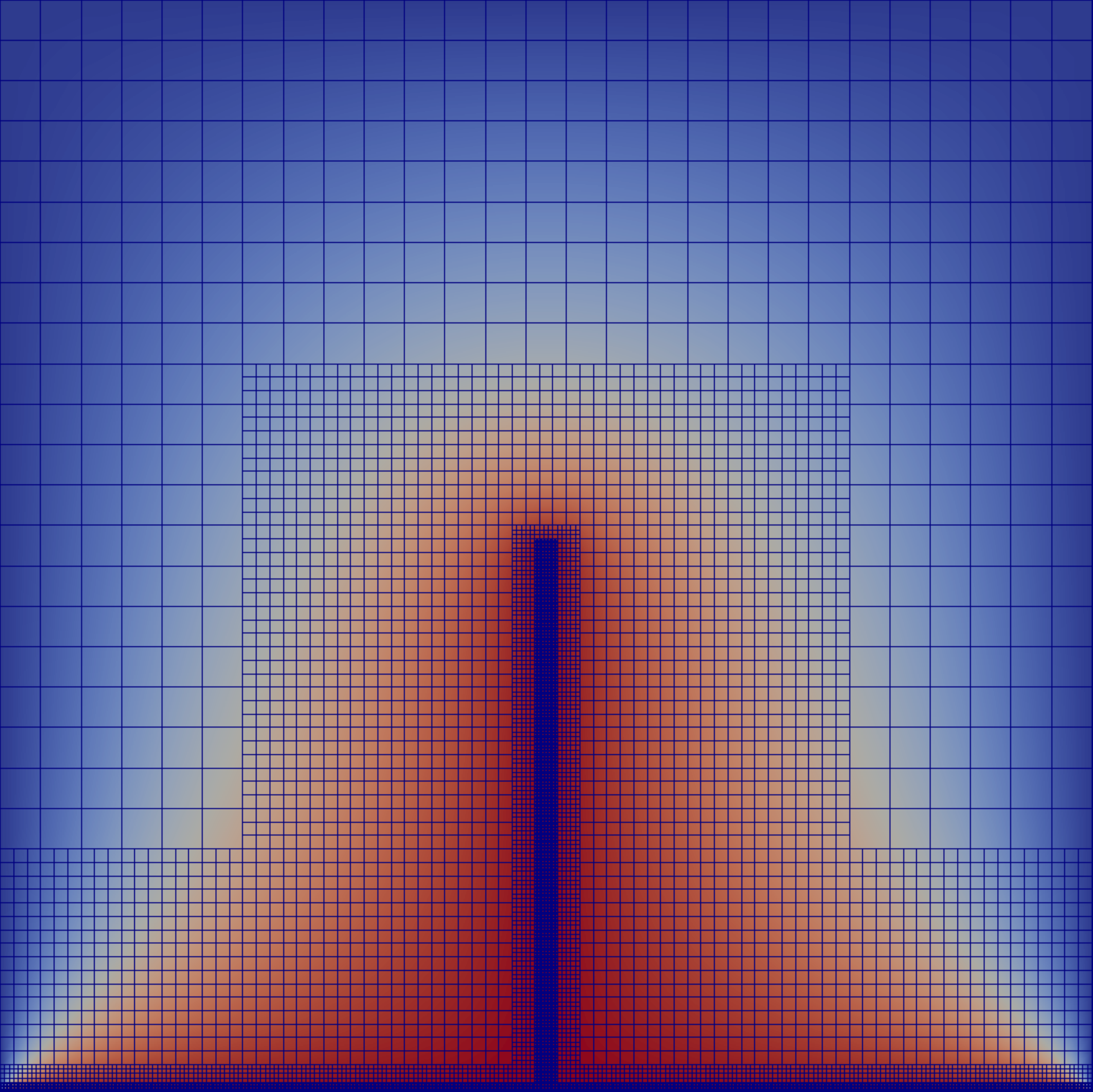}
   \includegraphics[width=0.39\textwidth]{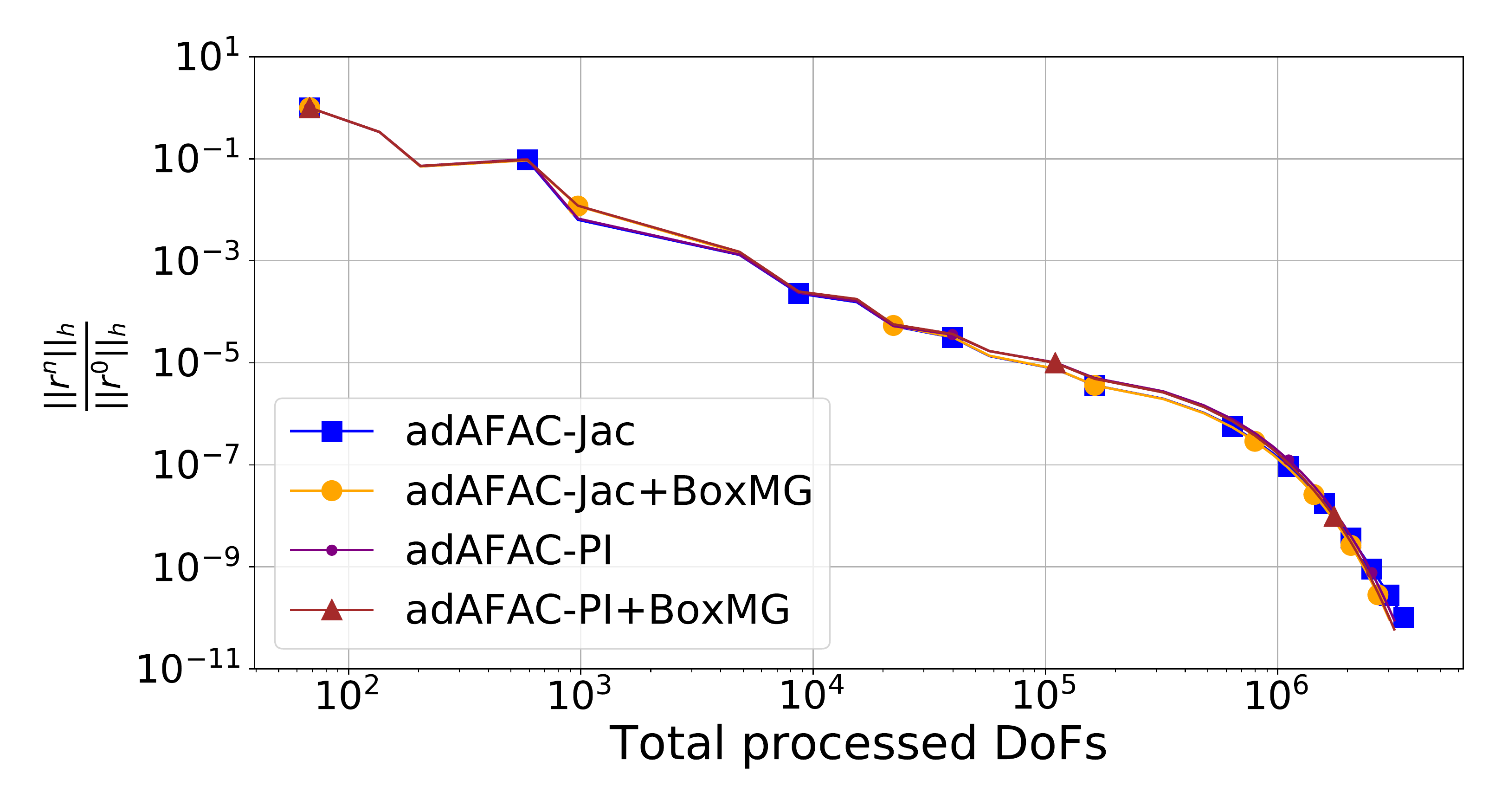}
   \includegraphics[width=0.39\textwidth]{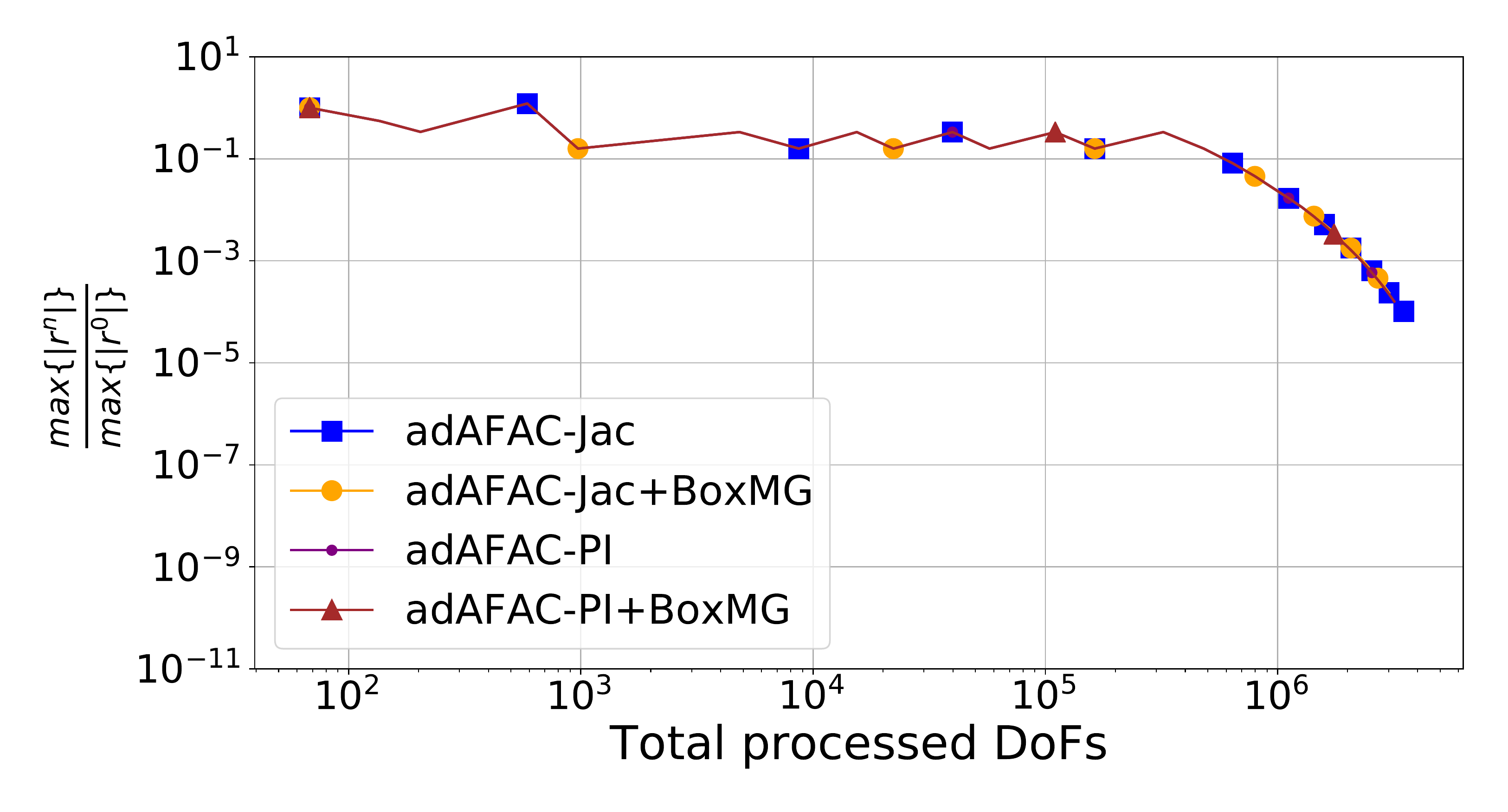}
   \\
   \includegraphics[width=0.39\textwidth]{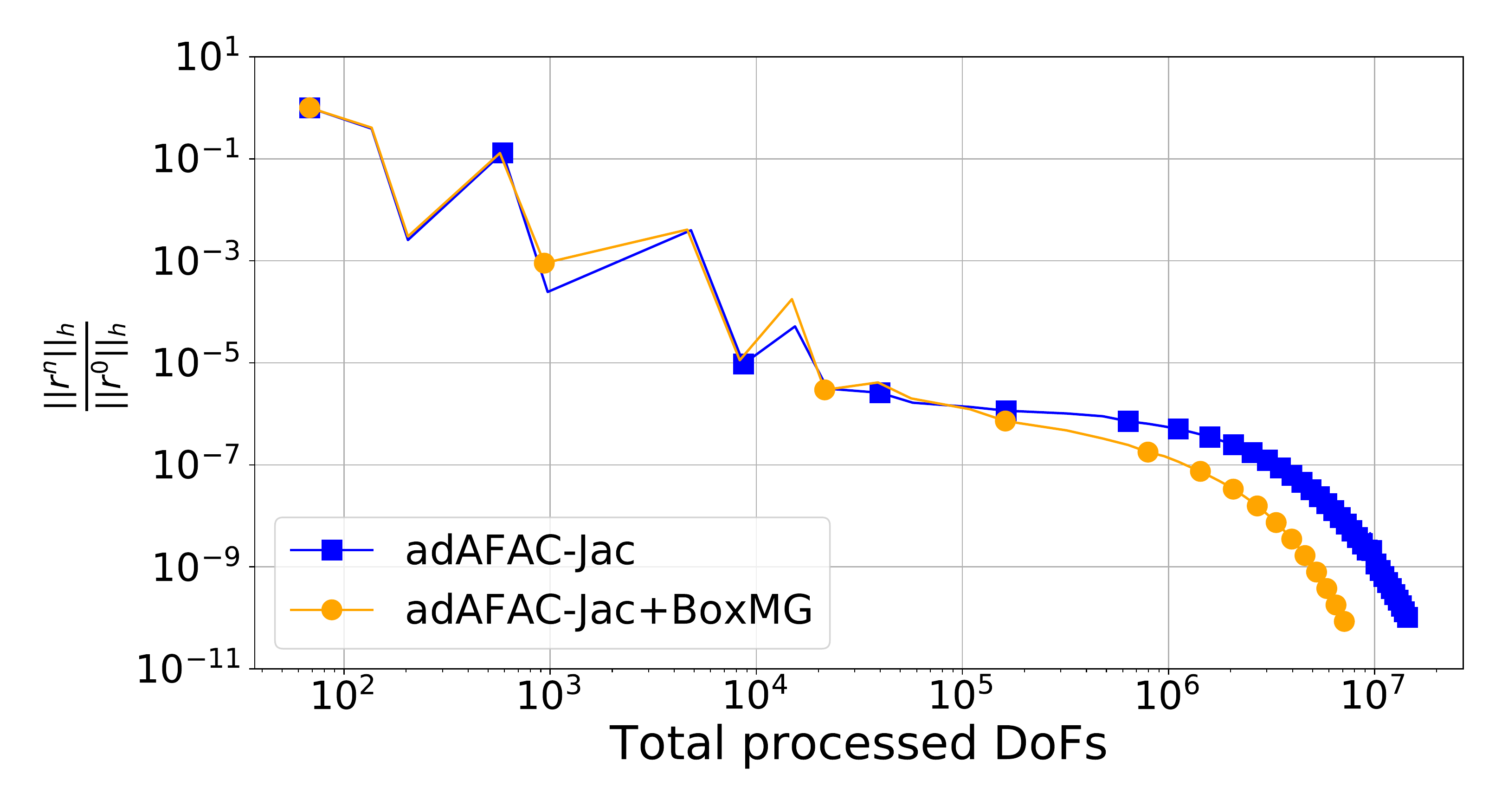}
   \includegraphics[width=0.39\textwidth]{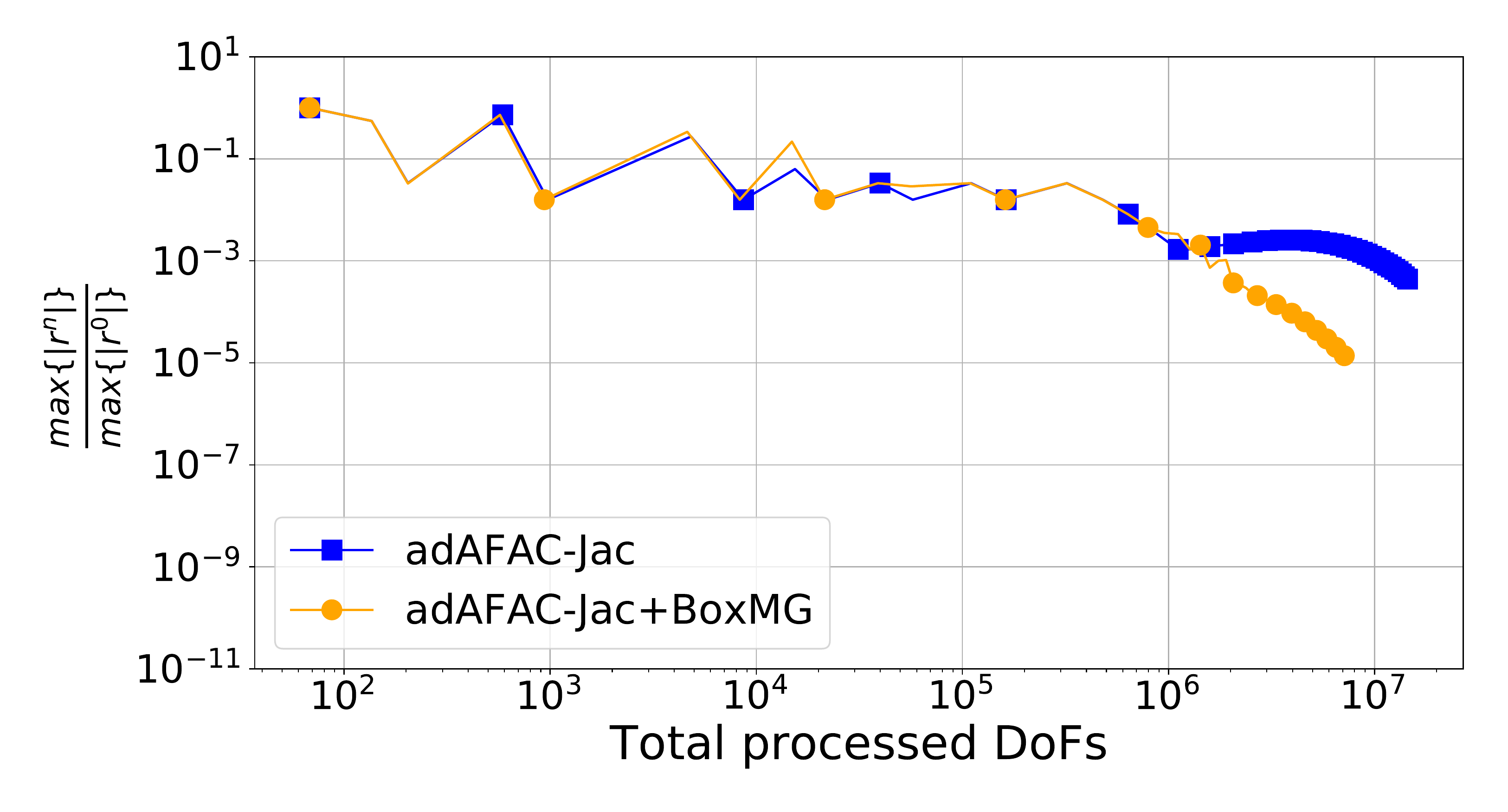}
 \end{center}
 \caption{
   Typical adaptive mesh for setup with a tiny, needle-like inclusion once the refinement criterion has
   stopped adding further elements (top left).
   The material inclusion either holds an $\epsilon $ which is bigger than its
   surrounding by a factor of ten (top row) or even by a factor of  1,000
   (bottom row).
   \label{figure:Poisson:adaptive-grid:inclusion}
 }
\end{figure}


Tiny, localised variations in $\epsilon $ are notoriously difficult to
handle for multigrid. 
The spike setup from our test suite yields a problem where diffusive
behaviour is ``broken'' along the inclusion.
The adaptivity criterion thus immediately refines along the tiny material spike
(Figure \ref{figure:Poisson:adaptive-grid:inclusion}) since the solution's
curvature and gradient there is very high.
We see diffusive behaviour around this refined area, but we know that there is
no long-range, smooth solution component overlapping the $\epsilon $ changes.

Again, a reasonable small variation in $\epsilon $ does not pose
major difficulties to either of our damped adAFAC solvers.
The strong localisation of the adaptivity ensures that the material transition
is reasonably handled, such that a sole geometric choice of inter-grid transfer
operators is totally sufficient.
However, this setup is challenging for additive multigrid
which fails to converge even with BoxMG.

Once we increase the material change by three orders of magnitude, we need an
explicit elimination of oscillations arising along the $\epsilon $ changes.
Solely employing algebraic BoxMG operators is insufficient.
They can mirror the solution behaviour to some degree but they are incapable to
compensate for the poor choice of our coarse grid points.
The present setup would require algebraic coarse grid
identification where the coarse grid aligns with the inclusion.

While adAFAC-PI with algebraic operators manages to obtain reasonable convergence for
a material variation of one order of magnitude nevertheless, it is unable to
converge for three orders of magnitude change even with algebraic inter-grid transfer operators.
adAFAC-Jac is able to handle the sharp, localised transition
which also can be read as extreme case of an anisotropic
$\epsilon $ choice in (\ref{equation:PDE}).
We see convergence for both its geometric variant and its algebraic
extension, though now the BoxMG variant is superior to its geometric
counterpart.

\begin{observation}
 adAFAC-Jac equips the geometric-algebraic BoxMG method with the opportunity to
 compensate, to some degree, for the lack of support of anisotropic refinement.
\end{observation}

\subsection{Non axis-aligned subdomains}

\begin{figure}
 \begin{center}
   \includegraphics[width=0.21\textwidth]{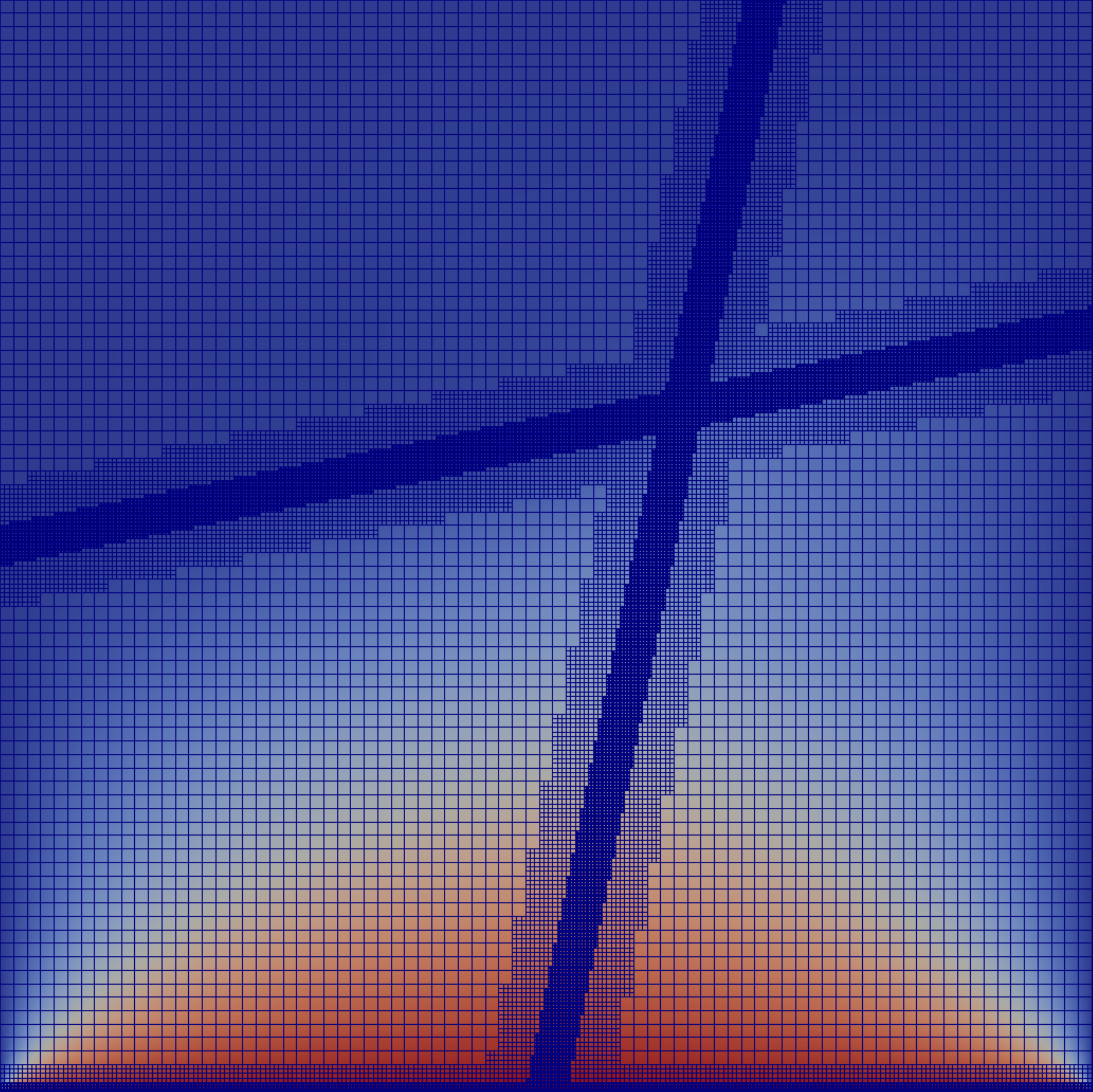}
   \includegraphics[width=0.39\textwidth]{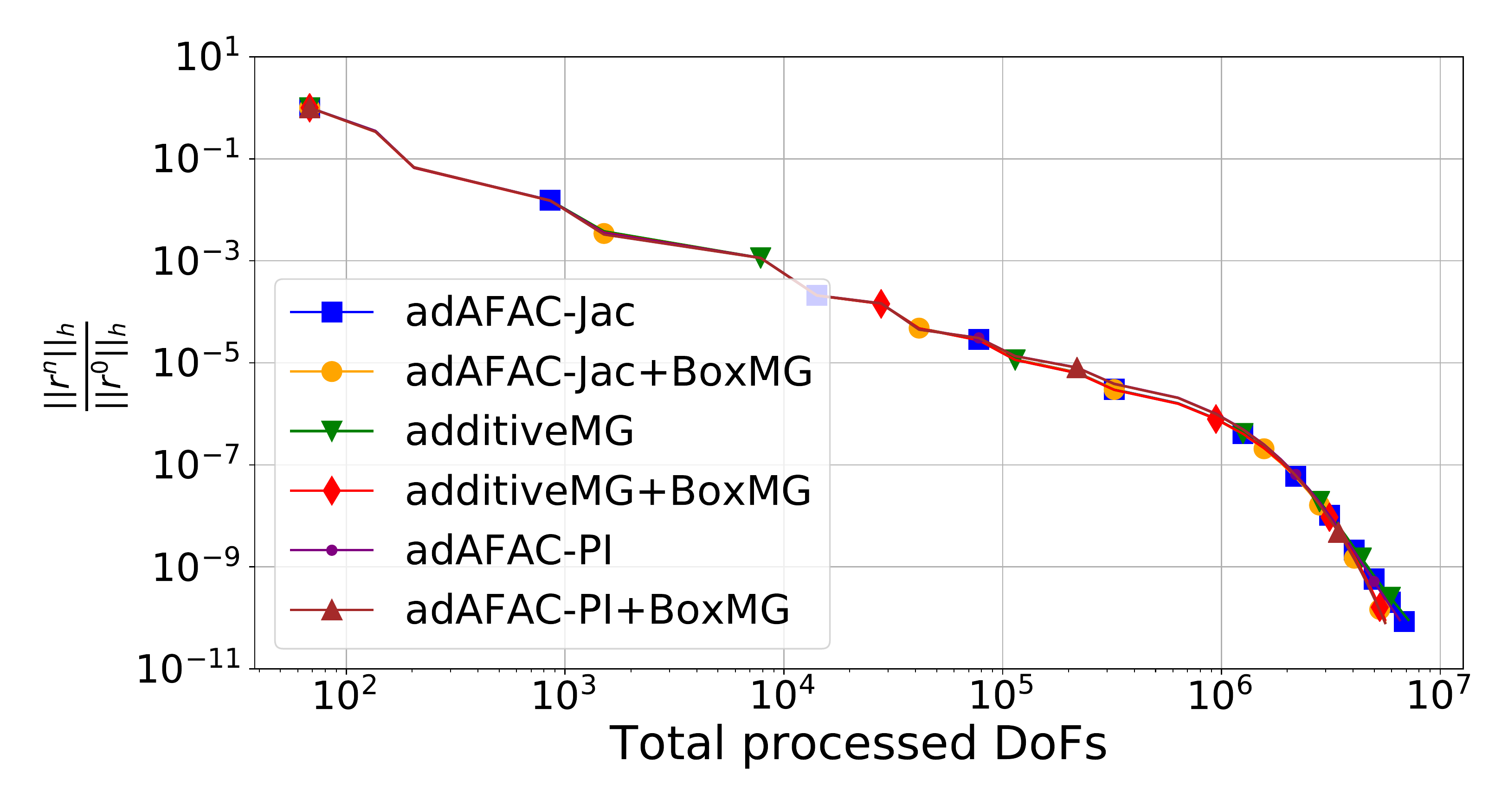}
   \includegraphics[width=0.39\textwidth]{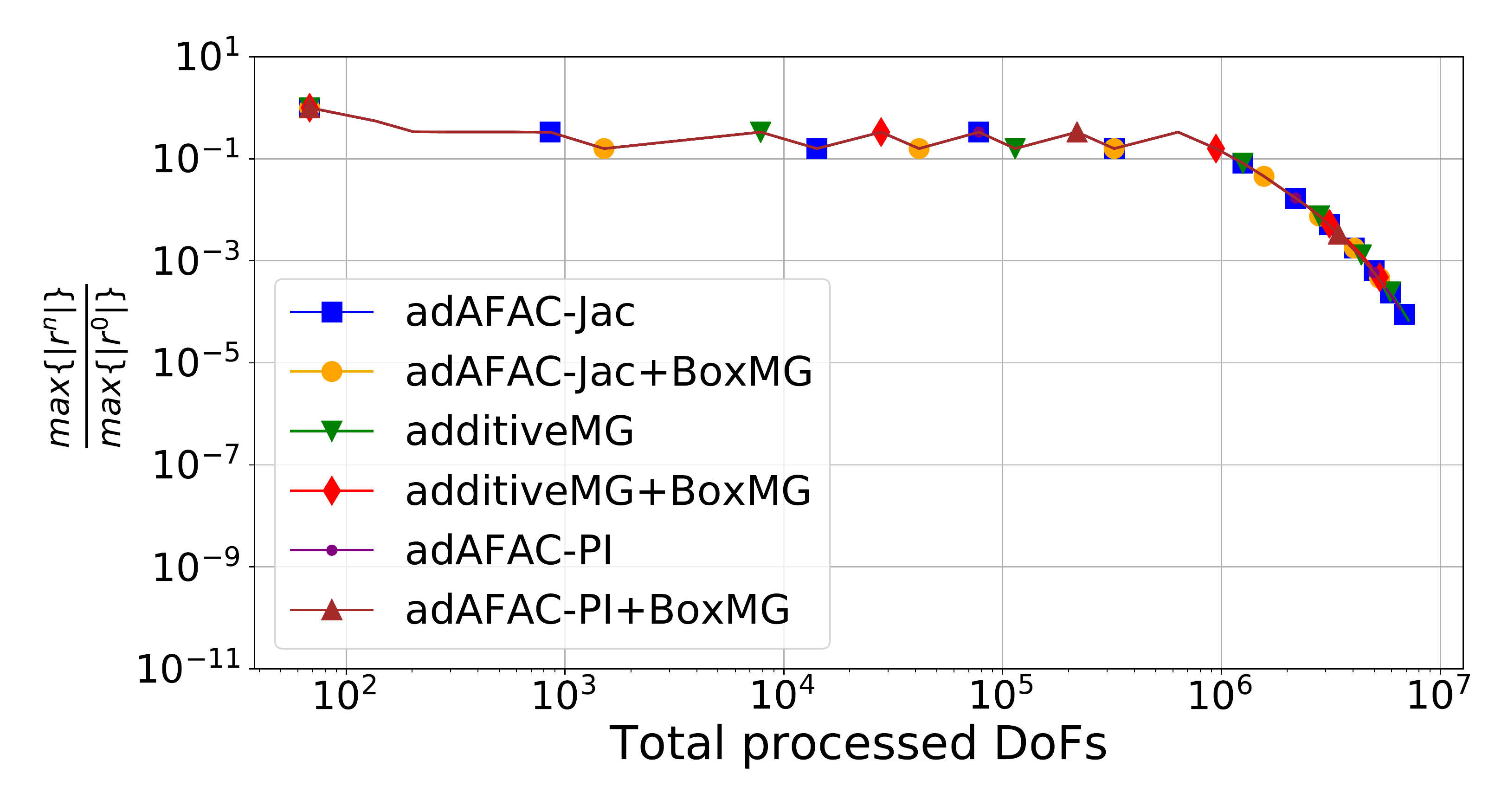}
   \\
   \includegraphics[width=0.39\textwidth]{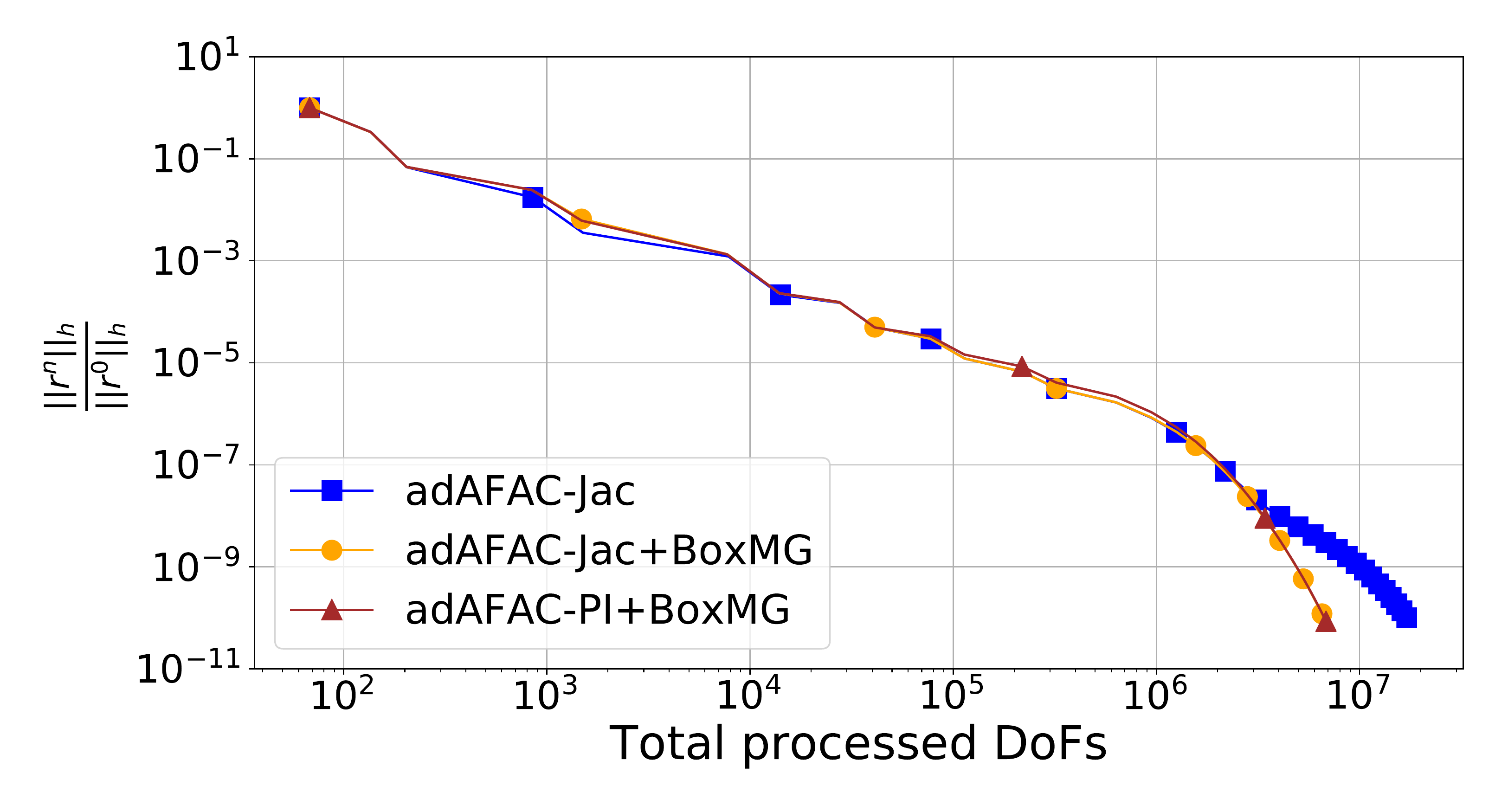}
   \includegraphics[width=0.39\textwidth]{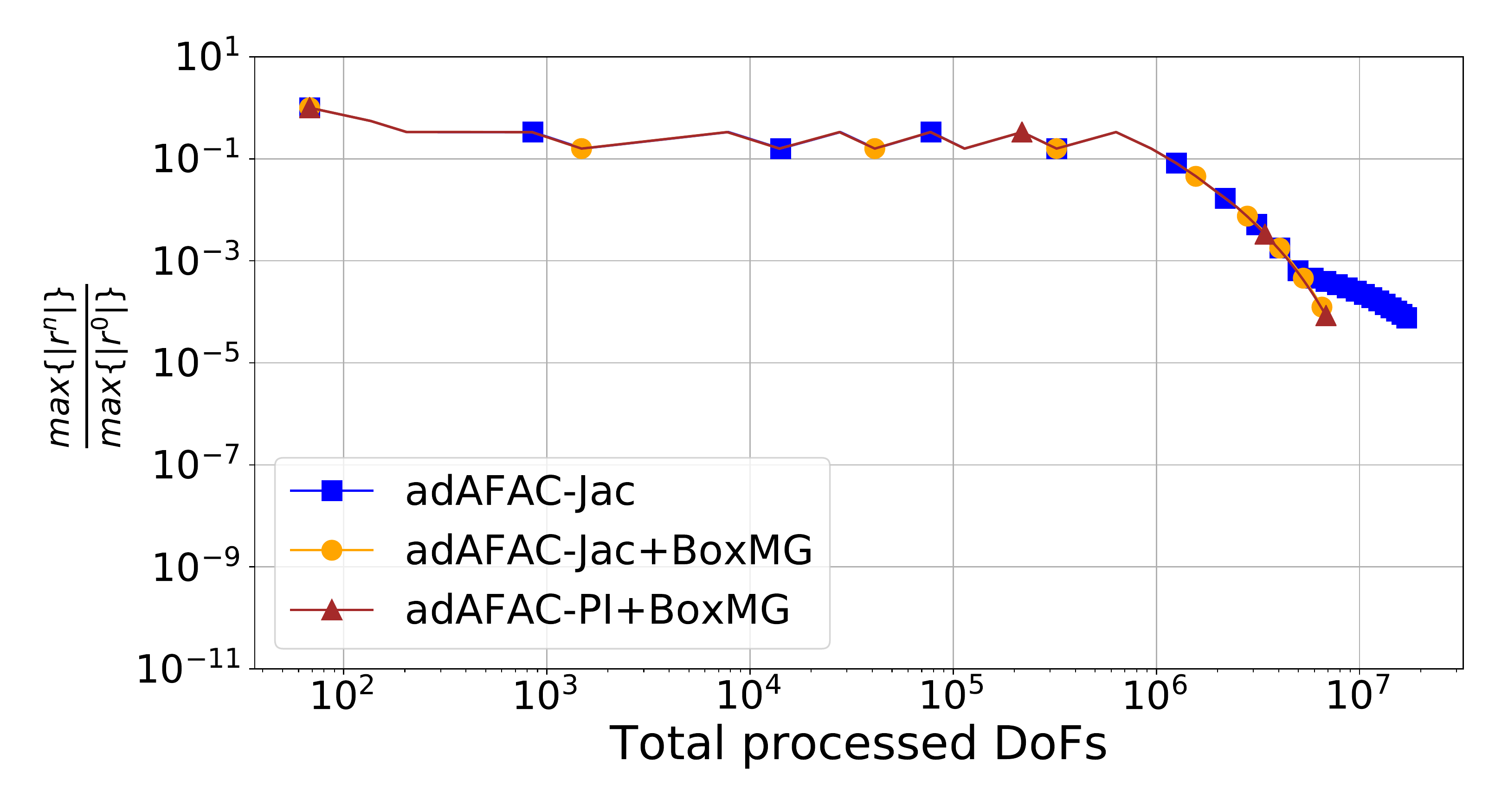}
 \end{center}
 \caption{
   Typical adaptive mesh for a setup where the regions with different
   material parameter $\epsilon $ are not axis-aligned.
   One order of magnitude differences in the material parameter (top)
   vs.~three orders of magnitude (bottom).
   \label{figure:Poisson:adaptive-grid:angled-lines-solvers}
 }
\end{figure}

%
%
We move on to our experimental setup with a deformed checkerboard setup (Figure
\ref{figure:Poisson:adaptive-grid:angled-lines-solvers}), where
the dynamic adaptivity criterion unfolds the mesh along the material
transitions.
The solution behaviour within the four subregions itself is smooth,
i.e.~diffusive, and the adaptivity around the material transitions thus is
wider, more balanced, than the hard-coded adaptivity directly at the bottom of
the domain.

With smallish variations in $\epsilon $, this setup does not pose a challenge to
any of our solvers, \replaced[id=R2]{irrespective of}{irrespectible} whether they work with algebraic or geometric
inter-grid transfer operators.
With increasing differences in $\epsilon $, we however observe that additive
multigrid starts to diverge. 
The smooth regions are still sufficiently dominant, and we suffer from
overcorrection.
adAFAC-PI performs better yet requires algebraic operators to remain robust up
to $\epsilon $ variations of three orders of magnitude.
adAFAC-Jac with geometric operators remains stable for all studied setups,
up to and including the five order of magnitude jump.
adAFAC-Jac with algebraic operators outperforms its
geometric cousin.
BoxMG's accurate handling of material transitions decouples the subdomains from
each other on the coarse correction levels.
Updates in one domain thus do not pollute the solution in a neighbouring domain.

\begin{observation}
 While the auxiliary equations can replace/exchange algebraic operators in some
 cases, they fail to tackle material transitions that are not grid-aligned.
\end{observation}

\subsection{Parallel adAFAC}

%
%
%
%
\added{
 We close our experiments with a scalability exercise.
 All data stem from a cluster hosting Intel E5-2650V4
 (Broadwell) nodes with 24 cores per node.
 They are connected via Omnipath.
 We use Intel's Threading Building Blocks (TBB) for the shared memory
 parallelisation and use Intel MPI for a distributed memory realisation.
}

\begin{figure}
 \begin{center}
  \includegraphics[width=0.45\textwidth]{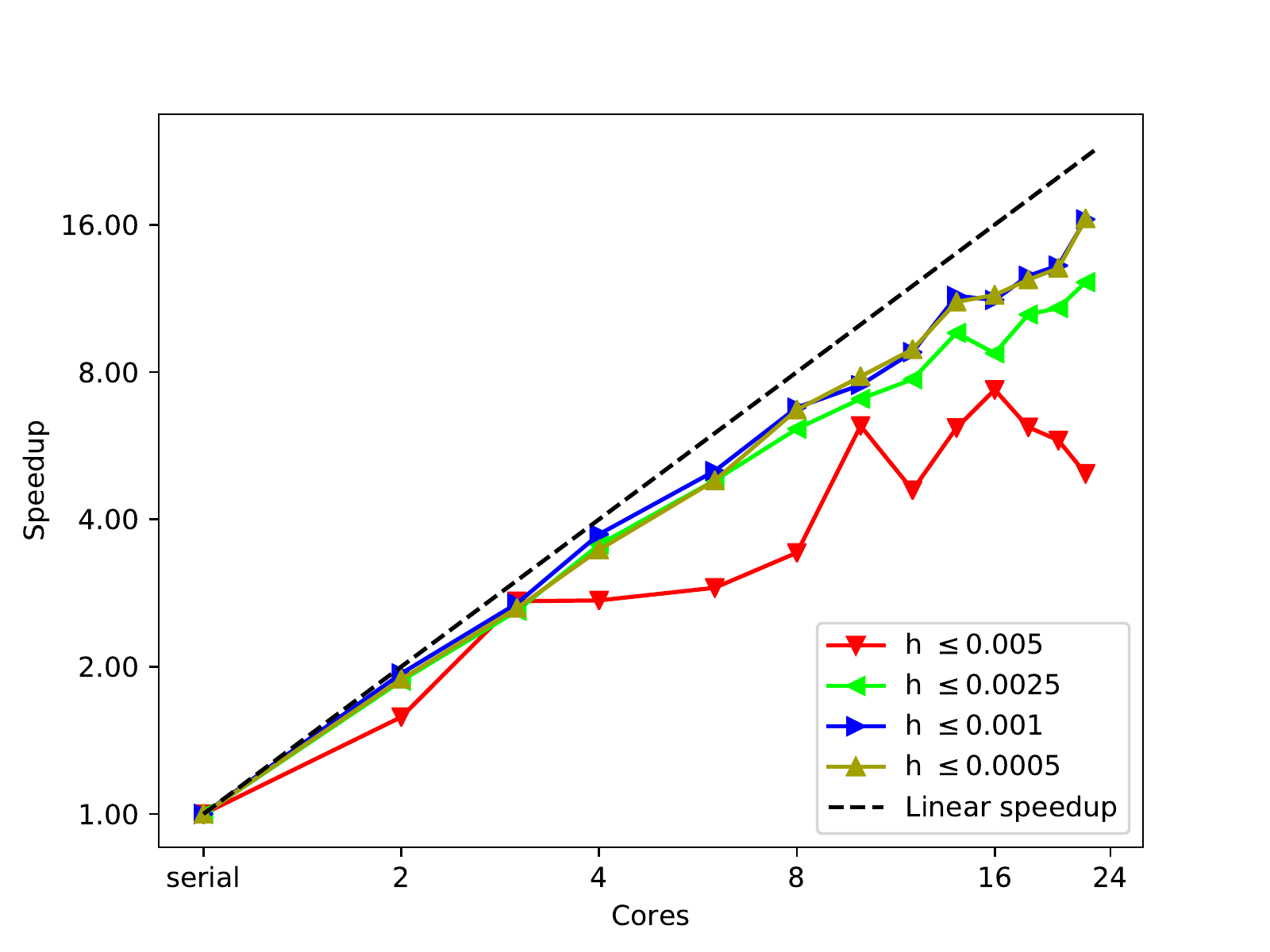}
  \includegraphics[width=0.45\textwidth]{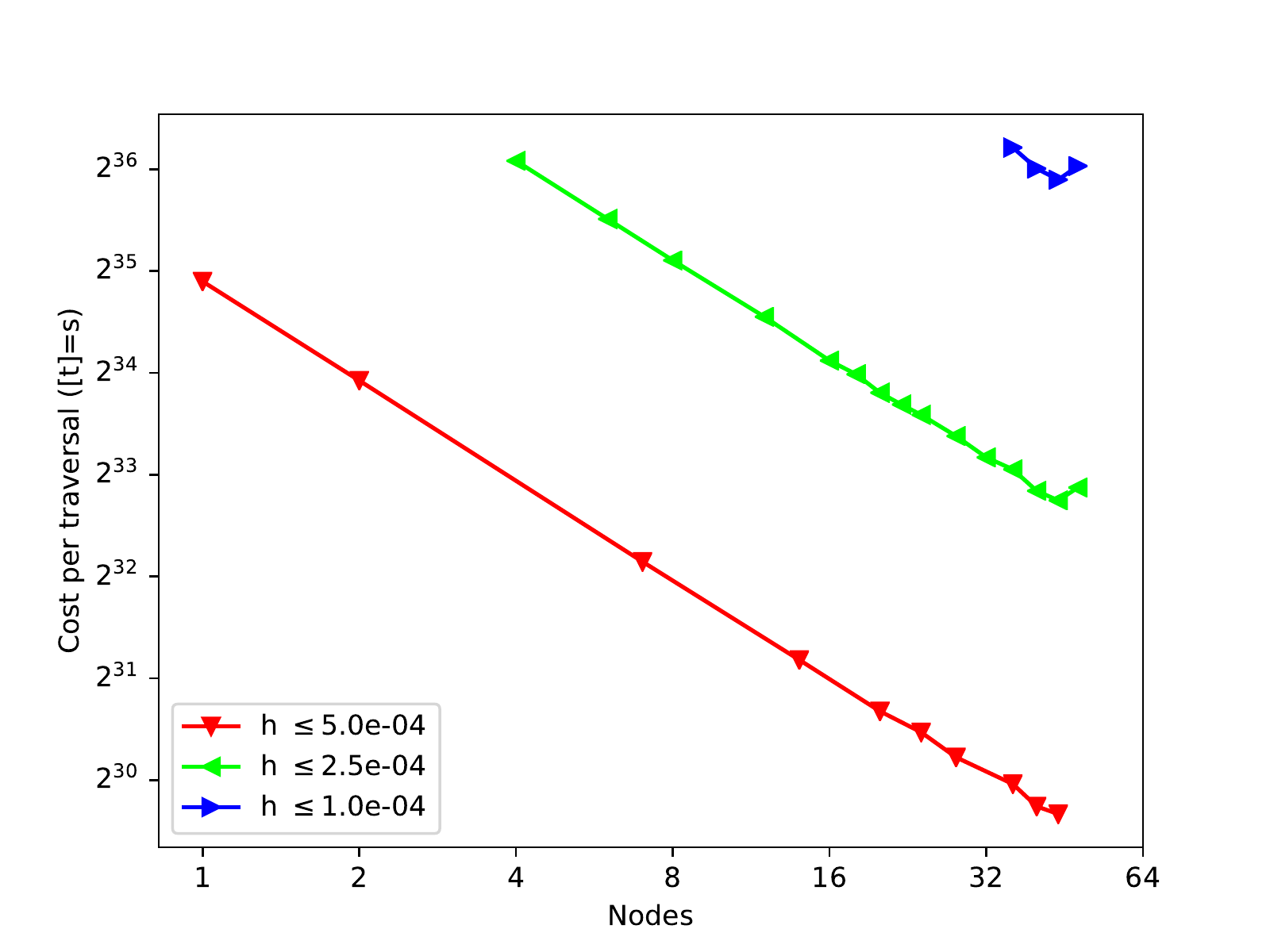}
 \end{center}
 \caption{
  \added{
    Left: Shared memory experiments with adAFAC.
    \added[id=R2]{All solver variants rely on the same code base, i.e.~exchange
    only operators, such that they all share the same performance
    characteristics.} 
    Right:
    Some distributed memory runtime results\added[id=R2]{ with the time for
    one multiscale grid sweep. This corresponds to one additive
    cycle as we realise single-touch semantics}.
  }
  \added[id=R3]{We study three different mesh sizes given via upper bounds on
  the $h$. Two ranks per node, i.e.~one rank per socket, are used.}
  \label{figure:multicore}
 }
\end{figure}

%
%
\added{
 Both the shared and the distributed memory parallelisation of our code use a
 multilevel space-filling curve approach.
 The fine grid cells are arranged along the Peano
 \replaced[id=R3]{space-filling}{space filling} curve and cut into curve
 segments of roughly the same number of cells.
 We use a
 non-overlapping domain decomposition on the finest mesh.
 Logically, our code does not distinguish between the code's shared and
 distributed memory strategy.
 They both decompose the data in the same way.
 The distributed memory variant however replaces memory copies along the
 boundary with MPI calls.
 All timings rely on runtimes for one cycle of a stationary mesh, i.e.~load
 ill-balances and overhead induced by adaptive mesh refinement are omitted.
}

\added[id=R3]{
 For all experiments, we start adAFAC and wait until the dynamic
 adaptivity has unfolded the grid completely such that it meets our prescribed
 $h$ as a maximum mesh size.
 We furthermore hardcode the domain decomposition such that the partitioning is
 close to optimal:
 We manually eliminate geometric ill-balancing, and we focus on the
 most computationally demanding cycles of a solve.
 Cycles before that, where the grid is not yet fully unfolded, yield performance
 which is similar to experiments with a bigger $h$.
}


%
%
\added{
 Our shared memory experiments \replaced[id=R3]{(Figure~\ref{figure:multicore})}{(Fig.~\ref{figure:multicore})} 
 show reasonable scalability up to eight cores if the 
 mesh is detailed.
 The curves are characteristic for both adAFAC-PI and adAFAC-Jac, i.e.~we have
 not been able to distinguish the runtime behaviour of these two approaches.
 If the mesh is too small, we see strong runtime variations. 
 Otherwise, the curves are reasonably smooth.
 Overall, the shared memory efficiency is limited by less than 70\% even if we
 make the mesh more detailed.
}

%
%
\added{
 Our code employs a very low order discretisation and thus exhibits a low
 arithmetic intensity.
 This intensity is increased by both adAFAC-PI and adAFAC-Jac, but the increase
 is lost behind other effects such as data transfer or the management of
 adaptive mesh refinement.
 The reason for the performance stagnation is not clear, but we may assume that
 NUMA effects play a role, and that communication overhead affects the 
 runtime, too.}
\added[id=R3]{With a distributed memory parallelisation, we can place two
 ranks onto each node. 
 NUMA then does not have further knock-on effects, and we obtain smooth
 curves until we run into too small partitions per node.
 With a low-order discretisation, our code is communication-bound---in line with
 most multigrid codes---yet primarily suffers from a strong synchronisation
 between cycles:
}

\replaced[id=R3]{Due to a non-overlapping domain decomposition on the finest
grid, all traversals through the individual grid partitions are synchronised with
each other.
}
{Besides architectural properties, the solvers however are intrinsically
 synchronised with each other.
}
\added{
 Our adAFAC implementation merges the coarse grid updates into the fine grid
 smoother, but each smoothing step requires a core to synchronise with all other cores.
 We eliminate strong scaling bottlenecks due to small system solves, but we have
 not yet eliminated scaling bottlenecks stemming from a tight synchronisation of
 the (fine grid) smoothing steps.
}

%
%
%
%
%
%

\begin{observation}
 \added{
 Despite adaFAC's slight increase of the arithmetic intensity, it seems to be
 mandatory to switch to higher order methods \cite{Gholami:16:SolverComparison}
 or approaches with higher asynchronicity
 \cite{Wolfson:19:AsynchronousMultigrid} to obtain better scalability.
 }
\end{observation}

\added{
 \noindent
 This is in line with other research \cite{Gholami:16:SolverComparison,Wolfson:19:AsynchronousMultigrid}.
}

\section{Conclusion and outlook}
\label{section:conclusion}

%
%
We introduce two additive multigrid variants which are conceptually close to
\replaced{asynchronous fast adaptive composite grid}{AFAC} solvers and Mult-additive.
An auxiliary term in the equation ensures that overshooting of plain additive
multigrid is \deleted{immediately} eliminated.
Our results validate that we obtain reasonable \added{multigrid} performance and
stability.
\added{
They confirm that adAFAC aligns with the
three key concepts from the introduction:
It relies solely on the geometric grid hierarchy,
it sticks to the asynchronous additive paradigm, and all new ideas
can be used in combination with advanced implementation patterns such 
as single-touch formulations or quasi matrix-free matrix storage.
Beyond that, the results} 
uncover \replaced{further}{three surprising} insights: 
\deleted{(i)} adAFAC \deleted{seems to be well-suited to
geometric construction. It }delivers
reasonable robustness \replaced{when}{while} solely using geometric inter-grid transfer
operators.
The construction of good inter-grid transfer operators \added{for non-trivial
$\epsilon $} is far from trivial and computationally cheap.
It is thus conceptually an interesting idea to give up on the idea of a
good operator and in turn to eliminate oscillations resulting from poor
operators within the correction equation.
We show that this is a valid strategy for some setups.
\replaced{In this context,}{(ii)} adAFAC
can be read as an antagonist to BPX.
BPX omits the system operator from the correction equations and ``solely''
relies on proper inter-grid transfer operators.
With our geometric adAFAC variants, we work without algebraic operators but a
problem-dependent auxiliary smoothing\added{, i.e.~a problem-dependent
operator.}
\deleted{We use additional operator information to improve the convergence
of an additive formulation
(iii) BoxMG is a powerful algebraic inter-grid transfer operator construction
scheme but cannot compensate for flaws that are introduced by problems that
require semicoarsening.
Our results do suggest that the auxiliary smoothing might be able to
improve the robustness of BoxMG in such scenarios.}

%
%
It is notoriously difficult to integrate multigrid ideas into existing
solvers.
Multigrid builds upon several sophisticated building blocks and needs
mature, advanced data structures. 
On the implementation side, an interesting contribution of our work is the
simplification and integration of the novel adAFAC idea into well-established concepts.
The fusion of three different solves (real solution, hierarchical solution
required for \replaced{the hierarchical transformation multigrid (HTMG)
implementation variant}{HTMG} and \added{the} damping equations) does not
introduce any additional implementational complexity compared to standard relaxation strategies.
However, it increases the arithmetic intensity.
adAFAC can be implemented as single-touch solver on dynamically adaptive grids.
This renders it an interesting idea for high performance codes relying on
dynamic, flexible meshes.

%
%
Studies from a high performance computing point of view are among our next
steps.
Interest in additive solvers has recently increased as they promise to
become a  seedcorn for asynchronous algorithms
\cite{Wolfson:19:AsynchronousMultigrid}.
Our algorithmic sketches integrate all levels' updates into one grid sweep and
thus fall into the class of vertically integrated solvers
\cite{Adams:16:SegmentalRefinement,Weinzierl:18:BoxMG}.
It will be interesting to study how desynchronisation interplays with the
present solver and single-touch ideas.
Further, we have to apply the scheme to more realistic\replaced[id=R2]{ and}{,} more challenging
scenarios.
Non-linear equations here are particularly attractive, as our adAFAC
implementation already offers a \replaced{full approximation storage}{FAS} data
representation.
\added{
 The multigrid community has a long tradition of fusing different ingredients:
 Geometric multigrid on very fine levels, direct solvers on very coarse levels,
 algebraic techniques in-between, \replaced[id=R2]{for example.}{e.g.}
 adAFAC is yet another solver variant within this array of options, and it will be
 interesting to see where and how it can work in conjunction with other multilevel solvers. 
}
On the \replaced{numerical method}{method} side, we expect further payoffs by improving 
\replaced{individual components of the solver---such as tailoring the smoother to our modified
restriction or modifying the prolongation in tandem}{the solver components}. 
Notably ideas following \cite{Yang:14:Reducing}
which mimic a $V(1,1)$-cycle
or even a $V$-cycle with more smoothing steps are worth investigating.

\section*{Acknowledgements}

The authors would like to thank Stephen F. McCormick.
Steve read through previous work of the authors \cite{Weinzierl:18:BoxMG} and
brought up the idea to apply the proposed concepts to AFACx. 
This kickstarted the present research into an AFAC variant.
The authors furthermore would like to thank Thomas Huckle who spotted the
elimination of the $\epsilon $-dependency in the auxiliary equation.
Finally, the authors thank Edmond Chow for his comments on BPX and for
the revitalising of the authors' research on FAC through his own work on
asynchronous multigrid.

\ifthenelse{\boolean{arxiv}}{
  \bibliographystyle{siam}
}{}
\ifthenelse{\boolean{sisc}}{
  \bibliographystyle{siamplain}
  \footnotesize
}{}
\ifthenelse{\boolean{linearalgebra}}{
}{}

\bibliography{paper}

\begin{thebibliography}{10}

\bibitem{Lin:18:Performance}
Lin PT, Shadid JN, Hu JJ, Pawlowski RP, and Cyr EC.
\newblock Performance of fully-coupled algebraic multigrid preconditioners for
  large-scale VMS resistive MHD.
\newblock Journal of Computational and Applied Mathematics. 2018;{\bf
  344}:782--793.

\bibitem{Weinzierl:18:BoxMG}
Weinzierl M, and Weinzierl T.
\newblock Quasi-matrix-free hybrid multigrid on dynamically adaptive Cartesian
  grids.
\newblock ACM Transactions on Mathematical Software. 2018;{\bf
  44}(3):32:1--32:44.

\bibitem{Dongarra:14:ApplMathExascaleComputing}
Dongarra J, Hittinger J, et~al.. {A}pplied {M}athematics {R}esearch for
  {E}xascale {C}omputing.
\newblock DOE ASCR Exascale Mathematics Working Group:
  http://www.netlib.org/utk/people/JackDongarra/PAPERS/doe-exascale-math-report.pdf;
  2014.

\bibitem{Gmeiner:14:Parallel}
Gmeiner B, K{\"o}stler H, St{\"u}rmer M, and R{\"u}de U.
\newblock Parallel multigrid on hierarchical hybrid grids: a performance study
  on current high performance computing clusters.
\newblock Concurrency and Computation: Practice and Experience. 2014;{\bf
  26}(1):217--240.

\bibitem{Lu:14:HybridMG}
Lu C, Jiao X, and Missirlis NM.
\newblock A hybrid geometric + algebraic multigrid method with semi-iterative
  smoothers.
\newblock Numerical Linear Algebra with Applications. 2014;{\bf
  21}(2):221--238.

\bibitem{May:15:Scalable}
May DA, Brown J, and Le~Pourhiet L.
\newblock A scalable, matrix-free multigrid preconditioner for finite element
  discretizations of heterogeneous Stokes flow.
\newblock Computer methods in applied mechanics and engineering. 2015;{\bf
  290}:496--523.

\bibitem{Sundar:12:ParallelMultigrid}
Sundar H, Biros G, Burstedde C, Rudi J, Ghattas O, and Stadler G.
\newblock Parallel Geometric-algebraic Multigrid on Unstructured Forests of
  Octrees.
\newblock In: Proceedings of the International Conference on High Performance
  Computing, Networking, Storage and Analysis. SC '12. Los Alamitos, CA, USA:
  IEEE Computer Society Press; 2012. p. 43:1--43:11.

\bibitem{Gholami:16:SolverComparison}
Gholami A, Malhotra D, Sundar H, and Biros G.
\newblock FFT, FMM, or Multigrid? A comparative Study of State-Of-the-Art
  Poisson Solvers for Uniform and Nonuniform Grids in the Unit Cube.
\newblock SIAM Journal on Scientific Computing. 2016;{\bf 38}(3):C280--C306.

\bibitem{Rudi:15:Extreme}
Rudi J, Malossi ACI, Isaac T, Stadler G, Gurnis M, Staar PW, et~al.
\newblock An extreme-scale implicit solver for complex PDEs: highly
  heterogeneous flow in Earth's mantle.
\newblock In: Proceedings of the international conference for high performance
  computing, networking, storage and analysis; 2015. p. 1--12.

\bibitem{Reps:17:Helmholtz}
Reps B, and Weinzierl T.
\newblock A Complex Additive Geometric Multigrid Solver for the Helmholtz
  Equations on Spacetrees.
\newblock ACM Transactions on Mathematical Software. 2017;{\bf
  44}(1):2:1--2:36.

\bibitem{Ghysels:12:PolynomialSmoothers}
Ghysels P, Klosiewicz P, and Vanroose W.
\newblock Improving the arithmetic intensity of multigrid with the help of
  polynomial smoothers.
\newblock Numerical Linear Algebra with Applications. 2012;{\bf
  19}(2):253--267.

\bibitem{Ghysels:13:ModelMG}
Ghysels P, and Vanroose W.
\newblock Modeling the performance of geometric multigrid stencils on multicore
  computer architectures.
\newblock SIAM Journal on Scientific Computing. 2015;{\bf 37}(2):C194--C216.

\bibitem{Gmeiner:15:HHG}
Gmeiner B, R\"ude U, Stengel H, Waluga C, and Wohlmuth B.
\newblock Towards Textbook Efficiency for Parallel Multigrid.
\newblock Numerical Mathematics: Theory, Methods and Applications. 2015;{\bf
  8}:22--46.

\bibitem{Weinzierl:11:Peano}
Weinzierl T, and Mehl M; SIAM.
\newblock {Peano -- A Traversal and Storage Scheme for Octree-Like Adaptive
  Cartesian Multiscale Grids}.
\newblock SIAM Journal on Scientific Computing. 2011;{\bf 33}(5):2732--2760.

\bibitem{Weinzierl:19:Peano}
Weinzierl T.
\newblock The {P}eano software - parallel, automaton-based, dynamically
  adaptive grid traversals.
\newblock ACM Transactions on Mathematical Software. 2019;{\bf
  45}(2):14:1--14:41.

\bibitem{Griebel:90:HTMG}
Griebel M.
\newblock Zur L{\"o}sung von Finite-Differenzen-und Finite-Element-Gleichungen
  mittels der Hierarchischen-Transformations-Mehrgitter-Methode [On the
  solution of the finite-difference and finite-element equations through the
  hierarchical-transformational-multigrid method].
\newblock Technische Universit{\"a}t M{\"u}nchen. Institut f{\"u}r Informatik;
  1990.

\bibitem{Trottenberg:01:Multigrid}
Trottenberg U, Oosterlee CW, and Sch{\"{u}}ller A.
\newblock {M}ultigrid.
\newblock Academic Press; 2001.

\bibitem{Dendy:82:BlackboxMG}
Dendy JE.
\newblock {B}lack {B}ox {M}ultigrid.
\newblock {J}ournal of {C}omputational {P}hysics. 1982;{\bf 48}(3):366--386.

\bibitem{Dendy:10:BoxMgBy3}
Dendy JE, and Moulton JD.
\newblock {B}lack {B}ox {M}ultigrid with {C}oarsening by a {F}actor of {T}hree.
\newblock Numerical Linear Algebra with Applications. 2010;{\bf 17}:577--598.

\bibitem{Bastian:98:Additive}
Bastian P, Wittum G, and Hackbusch W.
\newblock Additive and multiplicative multi-grid a comparison.
\newblock Computing. 1998;{\bf 60}(4):345--364.

\bibitem{Wolfson:19:AsynchronousMultigrid}
Wolfson-Pou J, and Chow E.
\newblock Asynchronous Multigrid Methods.
\newblock In: 2019 IEEE International Parallel and Distributed Processing
  Symposium (IPDPS). IEEE; 2019. p. 101--110.

\bibitem{McCormick:86:FAC}
McCormick SF, and Thomas J.
\newblock The fast adaptive composite grid ({FAC}) method for elliptic
  equations.
\newblock Mathematics of Computation. 1986;{\bf 46}(174):439--456.

\bibitem{McCormick:89:AFAC}
McCormick SF, and Quinlan DJ.
\newblock Asynchronous multilevel adaptive methods for solving partial
  differential equations on multiprocessors: Performance results.
\newblock Parallel Computing. 1989;{\bf 12}(2):145--156.

\bibitem{Lee:04:AFAC}
Lee B, McCormick SF, Philipp B, and Quinlan DJ.
\newblock Asynchronous Fast Adaptive Composite-Grid Methods for Elliptic
  Problems: Theoretical Foundations.
\newblock SIAM Journal Numerical Analysis. 2004;{\bf 42}:130--152.

\bibitem{Phillip:00:Elliptic}
Phillip B. 2000. \textit{Elliptic Solvers with Adaptive Mesh Refinement on
  Complex Geometries}.
\newblock Technical Report UCRL-JC-137372. Lawrence Livermore National Lab., CA
  (US).

\bibitem{Tuminaro:00:Parallel}
Tuminaro RS, and Tong C.
\newblock Parallel smoothed aggregation multigrid: Aggregation strategies on
  massively parallel machines.
\newblock In: Supercomputing, ACM/IEEE 2000 Conference. IEEE; 2000. p. 5--5.

\bibitem{Vanvek:96:Algebraic}
Van{\v{e}}k P, Mandel J, and Brezina M.
\newblock Algebraic multigrid by smoothed aggregation for second and fourth
  order elliptic problems.
\newblock Computing. 1996;{\bf 56}(3):179--196.

\bibitem{Vanvek:95:Fast}
Van{\v{e}}k P.
\newblock Fast multigrid solver.
\newblock Applications of Mathematics. 1995;{\bf 40}(1):1--20.

\bibitem{Yang:14:Reducing}
Vassilevski PS, and Yang UM.
\newblock Reducing communication in algebraic multigrid using additive
  variants.
\newblock Numerical Linear Algebra with Applications. 2014;{\bf
  21}(2):275--296.

\bibitem{Adams:16:SegmentalRefinement}
Adams MF, Brown J, Knepley MG, and Samtaney R.
\newblock Segmental Refinement: {A} Multigrid Technique for Data Locality.
\newblock {SIAM} J Scientific Computing. 2016;{\bf 38}(4).

\bibitem{Smith:96:DomainDecomposition}
Smith B, and P~Bj{\o}rstad WG.
\newblock Domain Decomposition---Parallel Multilevel Methods for Elliptic
  Differential Equations.
\newblock Cambridge University Press; 1996.

\bibitem{Hart:89:FAC}
Hart L, and McCormick SF.
\newblock Asynchronous multilevel adaptive methods for solving partial
  differential equations on multiprocessors: Basic ideas.
\newblock Parallel Computing. 1989;{\bf 12}:131--144.

\bibitem{Knuth:90:AttributeGrammar}
Knuth DE.
\newblock The genesis of attribute grammars.
\newblock In: Deransart P, and Jourdan M, editors. WAGA: Proceedings of the
  international conference on Attribute grammars and their applications.
  Springer-Verlag; 1990. p. 1--12.

\bibitem{Brandt:77:Multi}
Brandt A.
\newblock Multi-level adaptive solutions to boundary-value problems.
\newblock Mathematics of Computation. 1977;{\bf 31}(138):333--390.

\bibitem{Dubey:16:SAMR}
Dubey A, Almgren AS, Bell JB, Berzins M, Brandt SR, Bryan G, et~al.
\newblock A Survey of High Level Frameworks in Block-Structured Adaptive Mesh
  Refinement Packages.
\newblock Journal of Parallel and Distributed Computing. 2014;{\bf
  74}(12):3217--3227.

\bibitem{Hart:86:FAC}
Hart L, McCormick SF, and O'Gallagher A.
\newblock The Fast Adaptive Composite-Grid Method (FAC): Algorithms for
  Advanced Computers.
\newblock Applied Mathematics and Computation. 1986;p. 103--125.

\bibitem{Jimack:11:Asynchronous}
Jimack PK, and Walkley MA.
\newblock Asynchronous parallel solvers for linear systems arising in
  computational engineering.
\newblock Computational technology reviews. 2011;{\bf 3}:1--20.

\bibitem{Rude:87:Multiple}
R{\"u}de U.
\newblock Multiple tau-extrapolation for multigrid methods.
\newblock Bibliothek d. Fak. f{\"u}r Mathematik u. Informatik, TUM; 1987.

\bibitem{Yavneh:12:NonnsymBoxMg}
Yavneh I, and Weinzierl M.
\newblock {N}onsymmetric {B}lack {B}ox {M}ultigrid with {C}oarsening by
  {T}hree.
\newblock {N}umerical {L}inear {A}lgebra with {A}pplications. 2012;{\bf
  19}(2):246--262.

\bibitem{Press:91:Multigrid}
Press WH, and Teukolsky SA.
\newblock Multigrid Methods for Boundary Value Problems. I.
\newblock Computers in Physics. 1991;{\bf 5}(5):514--519.

\bibitem{Kouatchou:00:Optimal}
Kouatchou J, and Zhang J.
\newblock Optimal injection operator and high order schemes for multigrid
  solution of 3D Poisson equation.
\newblock International journal of computer mathematics. 2000;{\bf
  76}(2):173--190.

\bibitem{Bjorgen:17:Numerical}
Bj{\o}rgen J, and Leenaarts J.
\newblock Numerical non-LTE 3D radiative transfer using a multigrid method.
\newblock Astronomy \& Astrophysics. 2017;{\bf 599}:A118.

\bibitem{Bungartz:04:SparseGrids}
Bungartz HJ, and Griebel M.
\newblock Sparse grids.
\newblock Acta Numerica. 2004;{\bf 13}:147--269.

\bibitem{Mehl:06:MG}
Mehl M, Weinzierl T, and Zenger C.
\newblock A cache-oblivious self-adaptive full multigrid method.
\newblock Numerical Linear Algebra with Applications. 2006;{\bf
  13}(2--3):275--291.

\end{thebibliography}

\ifthenelse{\boolean{appendix}}{
 \appendix
 \section{Additive solvers in literature}

\begin{table}
 \caption{
   Overview of different additive solvers related to our approach.
   We write down one update resulting from a two-grid cycle with one smoothing
   step $M$.
   Prior to the start of any cycle, 
   we assume $u_\ell =0$ for all points spatially coinciding with points on 
   finer grids; unless stated differently. 
   Formula terms that make no sense (as levels do not exist) are treated as
   zero.
   \label{table:literature:overview}
 }
 \begin{center} 
 \begin{tabular}{p{4cm}p{7cm}p{5.5cm}}
   Name & Two-level formulation & Disadvantages \\
   \hline
   Additive MG 
     &
     $
     r_{\ell} = \left\{
      \begin{array}{ll}
       b_\ell - A_\ell u_\ell  
       & \mbox{on fine grid} 
       \\
       R r_{\ell +1} 
      \end{array}
     \right.
     $
     & Not stable as solver.\\
     &
     $
     u_{\ell} \gets u_{\ell} + \omega  M_{\ell } ^{-1} r_\ell + P u_{\ell
     -1} $ 
     \\
     \hline
   Additive MG (with exponential damping) 
     &
     $
     r_{\ell} = \left\{
      \begin{array}{ll}
       b_\ell - A_\ell u_\ell  
       & \mbox{on fine grid} 
       \\
       R r_{\ell +1} 
      \end{array}
     \right.
     $
     & No level-independent convergence rate.\\
     &
     $
     u_{\ell} \gets u_{\ell} + \omega  \left( M_{\ell } ^{-1} r_\ell + P u_{\ell
     -1} \right) $ 
     \\
     \hline
   BPX 
     &
     $
     r_{\ell} = \left\{
      \begin{array}{ll}
       b_\ell - A_\ell u_\ell  
       & \mbox{on fine grid} 
       \\
       R r_{\ell +1} 
      \end{array}
     \right.
     $
     & Usually used as preconditioner only. \\
     &
     $
     u_{\ell} \gets u_{\ell} + \left\{ 
      \begin{array}{ll}
      \omega M_{\ell } ^{-1} r_\ell + P u_{\ell -1}
       & \mbox{on fine grid} 
       \\
       \omega h_{\ell}^{d-2} r_\ell + P u_{\ell -1}
      \end{array}
     \right.
     $ 
     & Ignores coarse grid matrices but uses coarse mesh width $h_\ell$. \\
     \\
     \hline
   Hierarchical basis 
     &
     $
     r_{\ell} = (id - I^TI) \left( b_\ell - A_\ell u_\ell + R r_{\ell +1} \right) 
     $
     &
     Start with $u_{\ell -1} = I u_\ell$. Stable. Inferior to BPX as
     preconditionert \cite[p.~74]{Smith:96:DomainDecomposition}.
     \\
     &
     $
     u_{\ell} \gets u_{\ell} + \omega M_{\ell } ^{-1} r_\ell + P u_{\ell -1} -
     I^T I u_\ell $ 
     \\
     \hline
   AFACc 
     &
     $
     r_{\ell} = \left\{
      \begin{array}{ll}
       b_\ell - A_\ell u_\ell  
       & \mbox{on fine grid} 
       \\
       R (id-I^TI) r_{\ell +1} 
      \end{array}
     \right.
     $
     \\
     &
     $
     u_{\ell} \gets u_{\ell} + \omega M_{\ell } ^{-1} r_\ell + P u_{\ell -1}
     $ 
     \\
     \hline
   damped BPX-like
     & 
     $
     r_{\ell} = \left\{
      \begin{array}{ll}
       b_\ell - A_\ell u_\ell  
       & \mbox{on fine grid} 
       \\
       R r_{\ell +1} 
      \end{array}
     \right.
     $
     &
     As proposed in \cite{Reps:17:Helmholtz}.
     \\
     &
     $
     u_{\ell} \gets u_{\ell} + \omega M_{\ell } ^{-1} r_\ell + P u_{\ell -1}
     - P I \omega M_\ell ^{-1} r_\ell
     $ 
     \\
     &
     $
     \phantom{u_{\ell}} = u_{\ell} + \omega (id-P I) M_{\ell } ^{-1} r_\ell + P
     u_{\ell -1}
     $ 
 \end{tabular}
 \end{center}
\end{table}

We enumerate the levels top-down, i.e.~our fine grid carries level 
$\ell _{max}$ (the biggest index), the next coarser grid has level $\ell
_{max}-1$, and so forth.
Without loss of generality, we may assume that the coarsest multigrid level is
level 0.
Many papers in multigrid instead count bottom-up
\cite{Smith:96:DomainDecomposition,Trottenberg:01:Multigrid}. 
In those works, $\ell =0$ is the finest level, $\ell = 1$ the next coarser level, 
and so forth.

Classic (additive) multigrid starts from the finest level and successively
constructs coarser levels.
This leads to a cascade of grids where each coarser grid has fewer degrees of
freedom. 
The number of degrees of freedom
decreases monotonically.
As a result, coarse grids do preserve a certain spatial degree of freedom
concentration if the initial fine grid is strongly adaptive.
Additive multigrid then computes the residual, restricts it to all
levels, smoothes concurrently and finally sums up sequentially.
\begin{equation}
  u_{\ell_{max}} \gets u_{\ell_{max}} + \sum _{\ell = \ell _{max}} ^ {\tilde \ell _{min}}
   \tilde P ^{\ell_{max} - \ell} \tilde M_{\ell } ^{-1} \tilde R ^{\ell_{max} -
   \ell} \left( b_{\ell_{max}} - A_{\ell_{max}} u_{\ell_{max}} \right)
   \label{equation:literature:additive}
\end{equation}
if we use one smoothing step $M^{-1}$.
We use the tilde symbol to highlight that the underlying grids are not the ones
we use in the present paper.
In our work, they stem from a spacetree.
They are not real coarsened representations of the global fine grid mesh.
We use a formulation where the level hierarchy
corresponds directly to the grid construction rule, and write 
all operators without the tilde (Table
\ref{table:literature:overview}).
Furthermore, we ignore realisation difficulties arising from the fact that
vertices of one resolution might carry different semantics.
We finally conclude that a mesh of one resolution $\ell $ might carry fewer
degrees of freedom than the next coarser mesh $\ell -1$ if $\ell $ is only a
further refinement in a particular, small region of the domain.

Additive multigrid applied to grids that are embedded into each other yields a
very simple update rule if we ignore levels that are not populated on some
resolutions:  
we start with the finest resolution, likely not covering all of the domain, and
work our way up.
Additive multigrid is notoriously instable, but can be made stable if we
exponentially damp the grid levels. 
Unfortunately, we then loose multigrid, i.e.~level-independent, convergence
\cite{Bastian:98:Additive,Reps:17:Helmholtz,Smith:96:DomainDecomposition}.

BPX---classically used as preconditioner---removes the difficulty from additive
multigrid that we have to construct coarse grid operators.
Instead, it applies a simple diagonal coarse grid matrix which anticipates the
scaling of a ``real'' coarse grid matrix
\cite[pp.~73--74]{Smith:96:DomainDecomposition}.
It uses the inter-grid transfer operators but does not use the system matrix per
se on coarser levels.
We observe that for sole Poisson equations, cubic elements and a Jacobi
smoother, BPX yields exactly the vanilla additive multigrid scheme.
It thus obviously is not well-suited as solver unless we damp it artificially.

The hierarchical basis variant of the additive multigrid
\cite{Bastian:98:Additive} is typically introduced form a geometry point of
view: 
Assume we have a grid on level $\ell $.
When we add additional degrees of freedom on level $\ell +1$, we do not add
those shape functions whose support point spatially coincides with any point on
level $\ell $.
Writing it down within an additive FAC framework requires work:
We eliminate those vertices that the hierarchical basis grid construction does
not create a posteriori.
We first pretend they did exist.
However, as they exist only virtually, they do not carry a residual.
Combining the injection $I$ from level $\ell $ to $\ell -1$, i.e.~from coarse to
fine, times the transpose of the injection allows us to eliminate residuals from
the equation.
In return, we have to ensure that the coarse grid points hold valid $u$ values
plus valid right-hand sides and that these right-hand sides are taken into
account.
With zeroed residuals, the linear update through $M^{-1}$ does not yield an
update for virtual degrees of freedom, i.e.~dofs that are updated by coarser
levels.
In the update, we use the formula of BPX.
However, as a virtual vertex has already been injected to the coarser levels,
the value of this vertex is added again due to the prolongation $P$.
We therefore manually remove the current value of a virtual vertex when we
prolongate.

It is worth a note that the hierarchical basis variant has to be a hybrid
betweeen a multilevel solver and a multilevel correction scheme.
As we explicitly mask out fine level vertices that spatially coincide with
vertices holding coarse level degrees of freedom, we have to make these coarse
level vertices holds their right-hand side. 
They therefore combine discretisation terms plus a correction term; 
similar to our FAS.

Original work around additive FAC \cite{Hart:89:FAC} introduces the
AFACc scheme.
It relies on an operator which sets points manually to zero which spatially
coincide with points of a coarser grid.
It is a plain additive multigrid scheme which eliminates residual contributions
on those points which are updated on coarser levels, too, before it restricts.
The points continue to be updated on the fine level.
They just have no effect on the right-hand sides of the next coarser level.
This damps the right-hand sides.
In line with our hierarchical basis formula, we relise the AFACc masking
operator (called $P^k$ in \cite[p.~137]{Hart:89:FAC}.
) by a combination of the identity with injection and its transpose.
AFACc can be read as a merger of additive multigrid with the hierarchical basis
concept.

More sophisticated is the AFACf version of the same authors.
Here, they damp the correction equation's right-hand side by again going one
level up.
Different to AFACc, we thus may read AFACf as a damped AFACx scheme in the
sense of the present paper:
\begin{eqnarray*}
  u_{\ell} & \gets & u_{\ell} + 
   \left( P c_{\ell -1} + M_{\ell } ^{-1} \right) \left( b_\ell -
   A_\ell u_\ell \right) \quad \mbox{with} \\
  c_{\ell -1} & \gets & M_{\ell -1} ^{-1} (id - M_{\ell -1} ^{-1} \left(
   P M_{\ell -2} ^{-1} R \right)) R \left( b_\ell - A_\ell u_\ell \right).
\end{eqnarray*}

\noindent
In the above formula, we restrict ourselves to an iterative formulation.
The original paper relies on an explicit inversion of $A_{\ell -2}$
instead of a smoothing operator $M_{\ell -2}$ in the correction equation. 
This is problematic:
In an iterative world, the multiresolution scheme uses an additional multigrid
solve imbedded into the (top-down) run-through through the levels.
This is expensive \cite{Lee:04:AFAC}.

Finally, we introduced a damped BPX-like version for the Helmholtz
equations \cite{Reps:17:Helmholtz} within the FAC context:
When we update a vertex, we examine whether this same vertex has been updated on
a coarser level already.
If this is the case, we inject the update, prolongate it again and diminish the
fine grid update by this value.
We effectively multiply a fine grid correction with $id-PI$.
For diffusion operators, this manually eliminates the spreading of residual
impact and thus moves the plain additive scheme into the BPX regime.

 \section{An alternative adAFACx variant}

AFACx uses smoothing steps on the auxiliary grid levels to determine a start
datum for the actual additive multigrid solve
(Figure~\ref{figure:ingredients:adaFAC-scheme}).
Following the idea that a smoothing step followed by an inter-grid transfer
operator can be cast into a  smoothed operator, we observe that AFACx can be
read as standard additive scheme starting from a $\tilde P\tilde M^{-1}_{\ell
-1}RIu_{\ell +1}$ FAS solution representation.

Even though we approximate the smoothed operator by half-weighting---and thus
eliminate a manual smoothing---and stick to one Jacobi update for $M^{-1}_{\ell
-1}$, AFACx cannot
be written straightforwardly as single-touch algorithm. 
However, we can reapply our approximation trick: (i) Make AFACx's FAS start from
$Iu_{\ell +1}$ and run one additive smoothing step. (ii) Make AFACx's auxiliary
equation run completely assynchronously. (iii) Approximate the modified starting
point of the original AFACx algorithm by damping the Jacobi step outcome of (ii)
with $\tilde P\tilde M^{-1}_{\ell -1}RIu_{\ell +1}$.

We tested this modification.
Yet, modifying the restriction operators to mimic
pre-smoothing has been more effective in our experiments than modifying the
prolongation as sketched above.
We furthermore can motivate it clearly from $V(1,0)$-cycles whereas the sketched
AFACx variant seems to be closer to $V(0,1)$.
We therefore focus on the former approach.

 \section{Calculus (intermediate steps)}

The intermediate calculation for
(\ref{equation:integration:difference-add-mult}) reads as
\begin{eqnarray*}
  u_{\ell _{max, mult}} ^{(n+1)} - u_{\ell _{max, add}} ^{(n+1)} & = & 
  PA_{\ell _{max}-1}^{-1}R(
   b_{\ell _{max}} - A_{\ell _{max}}
   \left[ 
    u_{\ell _{max}} ^{(n)} + \omega _{{\ell_{max}}} M^{-1}_{{\ell_{max}}}  (
    b_{\ell _{max}} - A_{\ell _{max}} u_{\ell _{max}} ^{(n)} ) \right]
  )
  \nonumber 
  \\
  && + 
   \left[ 
    u_{\ell _{max}} ^{(n)} + \omega _{{\ell_{max}}} M^{-1}_{{\ell_{max}}}  (
    b_{\ell _{max}} - A_{\ell _{max}} u_{\ell _{max}} ^{(n)} ) \right]
  \\
  && -
  PA_{\ell _{max}-1}^{-1}R(
   b_{\ell _{max}} - A_{\ell _{max}}
    u_{\ell _{max}} ^{(n)}
  )
  \\
  && - 
   \left[ 
    u_{\ell _{max}} ^{(n)} + \omega _{\ell _{max}} M^{-1}_{\ell_{max}}  (
    b_{\ell _{max}} - A_{\ell _{max}} u_{\ell _{max}} ^{(n)} ) \right]
  \\
  & = & 
  PA_{\ell _{max}-1}^{-1}R(
   b_{\ell _{max}} - A_{\ell _{max}}
   \left[ 
    u_{\ell _{max}} ^{(n)} + \omega _{{\ell_{max}}} M^{-1}_{{\ell_{max}}}  (
    b_{\ell _{max}} - A_{\ell _{max}} u_{\ell _{max}} ^{(n)} ) \right]
  )
  \\
  && -
  PA_{\ell _{max}-1}^{-1}R(
   b_{\ell _{max}} - A_{\ell _{max}}
    u_{\ell _{max}} ^{(n)}
  ).
\end{eqnarray*}

}{}


\end{document}